\documentclass[a4paper,11pt,english]{article}
\usepackage{latexsym,amsfonts,amsopn}
\usepackage{amssymb}
\usepackage{babel}
\usepackage{amstext}
\usepackage{amscd}
\usepackage{amsthm}
\usepackage{amsmath}
\usepackage{graphicx}
\usepackage{epsfig}
\usepackage{rotating}
\usepackage{pb-diagram}
\swapnumbers
\newtheoremstyle{thomasdefinition}
     {.5cm}
     {.3cm}
     {}
     {}
     {\bf}
     {.}
     {.3em}
     {\thmnumber{({\bf\bfseries #2})} \thmname{\bf\bfseries #1}\thmnote{\bf\bfseries #3}}
\theoremstyle{thomasdefinition}
\newtheorem{dfn}{Definition}[section]
\newtheorem{exa}[dfn]{Example}
\newtheorem{re}[dfn]{Remark}

\newtheoremstyle{thomastheorem}
     {.5cm}
     {.3cm}
     {\it}
     {}
     {\bf}
     {.}
     {.3em}
     {\thmnumber{({\bf\bfseries #2})} \thmname{\bf\bfseries #1}\thmnote{\bf\bfseries #3}}
\theoremstyle{thomastheorem}
\newtheorem{cor}[dfn]{Corollary}
\newtheorem{lem}[dfn]{Lemma}
\newtheorem{thm}[dfn]{Theorem}
\newtheorem{prop}[dfn]{Proposition}

\newtheoremstyle{thomasempty}
     {.5cm}
     {.3cm}
     {\it}
     {}
     {\bf}
     {}
     {0cm}
     {\thmnumber{({\bf\bfseries #2})} \thmname{\bf\bfseries #1}\thmnote{\bf\bfseries #3}}
\theoremstyle{thomasempty}
\newtheorem*{empt}{}
\newcommand{\btab}{\begin{tabular}}
\newcommand{\etab}{\end{tabular}}
\newcommand{\be}{\begin{enumerate}}
\newcommand{\ee}{\end{enumerate}}
\newcommand{\bi}{\begin{itemize}}
\newcommand{\ei}{\end{itemize}}
\newcommand{\ba}{\begin{array}}
\newcommand{\ea}{\end{array}}
\newcommand{\bsm}{\begin{smallmatrix}}
\newcommand{\esm}{\end{smallmatrix}}
\newcommand{\ds}{\displaystyle}
\newcommand{\sst}{\scriptstyle}

\newcommand{\mb}[1]{{\mathversion{bold}{\bf #1}}}

\renewcommand{\cal}[1]{\ensuremath{\mathcal{#1}}}
\newcommand{\ub}[4]{\setlength{\unitlength}{1mm} \begin{picture}(0,0)
 \put(0,0){\raisebox{#2mm}{$\underbrace{\hspace{#1mm}}_{#4}$}} \end{picture} \raisebox{#3mm}{\rule{0mm}{#3mm}} }
\newcommand{\ob}[4]{\setlength{\unitlength}{1mm} \begin{picture}(0,0)
 \put(0,0){\raisebox{#2mm}{$\overbrace{\hspace{#1mm} \rule{0mm}{2ex}}^{#4}$}} \end{picture} \rule{0mm}{#3mm} }
\newcommand{\stru}[1]{\makebox[0cm]{$\ba{c}\rule{0mm}{#1mm}\ea$}} 
\newcommand{\nosp}[3]{\setlength{\unitlength}{1mm} \begin{picture}(0,0) \put(#1,#2){#3} \end{picture}}
\newcommand{\wt}{\widetilde}
\newcommand{\wh}{\widehat}
\newcommand{\ol}{\overline}
\newcommand{\ra}{\rightarrow}
\newcommand{\kov}{\nabla}
\newcommand{\lp}{\left(}
\newcommand{\rp}{\right)}
\newcommand{\menge}[2]{ \big\{\, #1 \,\big|\, #2 \big\} }
\newcommand{\Menge}[2]{ \Big\{\, #1 \;\Big|\; #2 \Big\} }
\newcommand{\rca}[2]{ \renewcommand{\arraystretch}{#1}{#2} }
\newcommand{\R}{\ensuremath{\mathbb{R}}}
\newcommand{\Z}{\ensuremath{\mathbb{Z}}}
\newcommand{\N}{\ensuremath{\mathbb{N}}}
\newcommand{\G}{\ensuremath{\mathcal{G}}}

\newcommand{\RR}{\ensuremath{\mathcal{R}}}

\renewcommand{\H}{\ensuremath{\mathcal{H}}}

\newcommand{\B}{\ensuremath{\mathcal{B}}}
\newcommand{\Zent}{\ensuremath{\mathcal{Z}}}
\newcommand{\g}{\ensuremath{\mathfrak{g}}}
\renewcommand{\k}{\ensuremath{\mathfrak{k}}}
\newcommand{\h}{\ensuremath{\mathfrak{h}}}

\newcommand{\m}{\ensuremath{\mathfrak{m}}}

\newcommand{\s}{\ensuremath{\mathfrak{s}}}
\newcommand{\z}{\ensuremath{\mathfrak{z}}}

\newcommand{\ad}{\ensuremath{\mathfrak{ad}}}
\newcommand{\gl}{\ensuremath{\mathfrak{gl}}}
\newcommand{\so}{\ensuremath{\mathfrak{so}}}
\renewcommand{\o}{\ensuremath{\mathfrak{o}}}

\renewcommand{\sl}{\ensuremath{\mathfrak{sl}}}

\newcommand{\rad}{\ensuremath{\mathfrak{rad}}}
\newcommand{\iso}{\ensuremath{\mathfrak{is\mspace{-2mu} o}}}
\newcommand{\hol}{\ensuremath{\mathfrak{hol}}}
\newcommand{\aut}{\ensuremath{\mathfrak{aut}}}
\newcommand{\Iso}{\ensuremath{I\mspace{-3mu}so}}
\newcommand{\Aff}{\ensuremath{A\mspace{-2mu}f\mspace{-5mu}f}}
\newcommand{\Diff}{\ensuremath{Di\mspace{-2mu}f\mspace{-5mu}f}}

\DeclareMathOperator{\Sig}{Sig}
\DeclareMathOperator{\LA}{LA}

\advance\oddsidemargin by -25mm
\advance\evensidemargin by -25mm
\advance\topmargin by -1.0cm 

\begin{document}

%

\date{}
\title{Solvable Pseudo-Riemannian Symmetric Spaces}
\author{Thomas Neukirchner\thanks{e-mail: neukirch@mathematik.hu-berlin.de - \today}}
\maketitle
\abstract{
We present an approach to solvable pseudo-Riemannian symmetric spaces based on 
papers of M.Cahen, M.Parker and N.Wallach. Thereby we
reproduce the classification of solvable symmetric triples of Lorentzian signature
$(1,n-1)$ and complete the case of signature $(2,n-2)$.
Moreover we discuss the topology of non-simply-connected symmetric spaces.
}
\tableofcontents
\pagebreak

\setlength{\parindent}{0cm}
\setlength{\parskip}{2ex}
\setlength{\hfuzz}{3pt}

\begin{flushright}
\parbox{7.5cm}{\small \it "Symmetry, as wide or as narrow you may define its meaning, is one idea by which man 
through the ages has tried to comprehend and create order, beauty, and perfection."
\begin{flushright}
\emph{Hermann Weyl}
\end{flushright}
}
\end{flushright}
 
\setcounter{section}{1}
\section*{Introduction}
\addcontentsline{toc}{section}{Introduction}

\begin{dfn}[Symmetric space]\label{dfn:SS}
A smooth connected pseudo-Riemannian manifold $(M,g)$ is a symmetric space if any point $x \in M$ 
is an isolated fixed point of an involutive isometry $s_x$, called the symmetry at $x$. 
\end{dfn}

\vspace{-1em}
This definition represents clearly the geometric meaning of the attribute ``symmetric''.
To emphasize the algebraic point of view consider the multiplication $x \bullet y := s_x y$  given by the symmetry. 
Then the following properties can be derived:
\begin{center}
{\renewcommand{\arraystretch}{1.3}
\btab{l  p{11cm}}
(M0) & $\bullet : M \times M \ra M$ is smooth \\
(M1) & $x \bullet x = x$ (fixed point property) \\
(M2) & $x \bullet (x \bullet y) = y$ (involutivity)\\
(M3) & $x \bullet (y \bullet z)=(x \bullet y) \bullet (x \bullet z)$ (left-multiplication is an automorphism) \\
(M4) & Any $x \in M$ has a neighborhood $U(x)$ s.t. $x \bullet y =y$ implies $y=x$ 
     if $y \in U(x)$ (isolation of the fixed point)
\etab
}
\end{center}
Conversely, O.Loos proved that (M0)-(M5) yields an affine symmetric space $(M,\kov)$ (i.e. $s_x$ as in \ref{dfn:SS} is an affine map) 
w.r.t. to a canonical connection $\kov$ (\cite[II.3.1]{Loos69}).
As example consider a Lie group with the symmetry product given by $x \bullet y = x \cdot y^{-1} \cdot x$. Although 
this product is different from the original group product many properties of Lie groups carry over to symmetric spaces
(cf. \cite[II.3.2.b, II.2.8.a, III.1.7/8, IV.2.2]{Loos69}). Thus symmetric spaces generalize Lie groups from an algebraical point of view.
Furthermore the natural notion of homomorphism between pairs $(M,\bullet)$ coincides with the concept of affine maps (\cite[II.2.6]{Loos69}).
In particular  $Aut(M, \bullet) = \Aff(M,\kov)$, i.e. (M3) extends to the useful identity
\begin{equation}\label{eq:Aut=Aff}
\Phi \circ s_x = s_{\Phi(x)} \circ \Phi, \quad \forall \, \Phi \in \Aff(M,\kov)
\end{equation}
On the other hand symmetry interdicts the existence of any nontrivial tensor field of odd degree which is invariant
under point reflections. In particular the covariant derivative of the curvature tensor field $\RR$ has to vanish.
Thus the differential geometric version of a symmetric space is obtained by requiring $\kov \RR=0$. 
About 1926 E.Cartan recognized that conversely any (pseudo-)Riemannian manifold with parallel curvature is at least
a {\bf locally symmetric space} i.e. locally isometric to a symmetric space.
Examples are again Lie groups equipped with a bi-invariant metric.
Thus symmetric spaces generalizes Lie groups also from a geometric point of view. 
An overview on the various concepts of symmetric spaces, 
their associated notions of morphism and their infinitesimal versions is given in\cite{Bertram00}.

Our interest in symmetric spaces arises from the fact that their prime representation as homogenous space 
by means of the transvection group yields their holonomy as linearization of the isotropy group. In the Riemannian case
this was used by M.Berger to determine all possible holonomies whereas in the pseudo-Riemannian setting occur 
non-completely-reducibly acting holonomy groups - corresponding to symmetric spaces with non-semi-simple transvection group - 
which both were sparcely handled in the past.
The first examples for such spaces were again groups with bi-invariant metric (e.g. \cite{Wu67}).
Efforts to understand their general structure were made by A.Medina and P.Revoy (\cite{MedinaR85}, see also \cite{Bordemann97}) 
using double-extensions of Lie algebras
or more recently by \cite{KathO02}, who describe a metric Lie algebras by a twofold extension, 
which promises to be a more canonical approach.
These techniques could be transfered to an arbitrary symmetric space taking the involutive automorphism
(induced by conjugation with a symmetry) 
on the transvection group with its canonical bi-invariant metric into account.
Due to the lack of such a theory the classification results
of M.Cahen, M.Parker and N.Wallach concerning symmetric spaces of Lorentzian type and signature $(2,n-2)$ were obtained 
by fixing a basis in a tangent space which is in a certain sense adapted to the metric and discussing the implications of the Jacobi-identity
for an ansatz of the structure coefficients of the transvection algebra
(\cite{CahenP70},\cite{CahenW70}, \cite{Cahen72} \cite{CahenP80}, \cite{Cahen98}). It is the main purpose of the present paper,
which is a shortened and updated version of my diploma thesis \cite{Neukirchner02} to revisit 
those results, complete them where necessary and present them in a uniform manner in order to make the sources
of potentially further mistakes more visible. 

Another motivation for this work was to use the explicite knowledge of the holonomy group to find examples of pseudo-Riemannian
manifolds which admit parallel spinors by application of the holonomy principle. Their geometries are classified only for irreducible
manifolds \cite{BaumK99}. Otherwise local normal forms of such metrics has been constructed in \cite{Bryant00}, \cite{Kath00}, \cite{Leistner02} and 
as geodesically complete examples Lie groups with bi-invariant metric has been studied under this point of view in \cite{BaumK02}
and Lorentzian symmetric spaces with solvable transvection group in \cite{Baum00}.
Similar calculations might apply to the results of the present paper, however this work is not done yet. 

In the fist chapter we recall some elementary facts on symmetric spaces based on the notion of the transvection group
which coincides only in the Riemannian case with the connected component of the isometry group. 
After introducing symmetric triples 
we explain their one-to-one-correspondence to simply connected symmetric spaces and set up a little dictionary
between their geometric properties and the respective algebraic features of symmetric triples.
In particular it is indicated that the complexity and variety of symmetric triples is growing while the
signature of the pseudo-Riemannian metric increases as can be seen from the following figure.
As mentioned above the semi-simple symmetric triples has been classified in \cite{Berger57}. 
If the holonomy acts non-irreducibly there are exactly two invariant, irreducible and totally isotropic subspaces
provided the manifold is indecomposable (cf. \cite{Koh65}), thus the signature of the metric is $(n,n)$.
The case of a transvection group which admits a proper Levi-decomposition is treated in \cite{CahenP80}.
 
\begin{figure}[!h]\label{fig:tree}
\begin{center}
\scalebox{.9}{\input{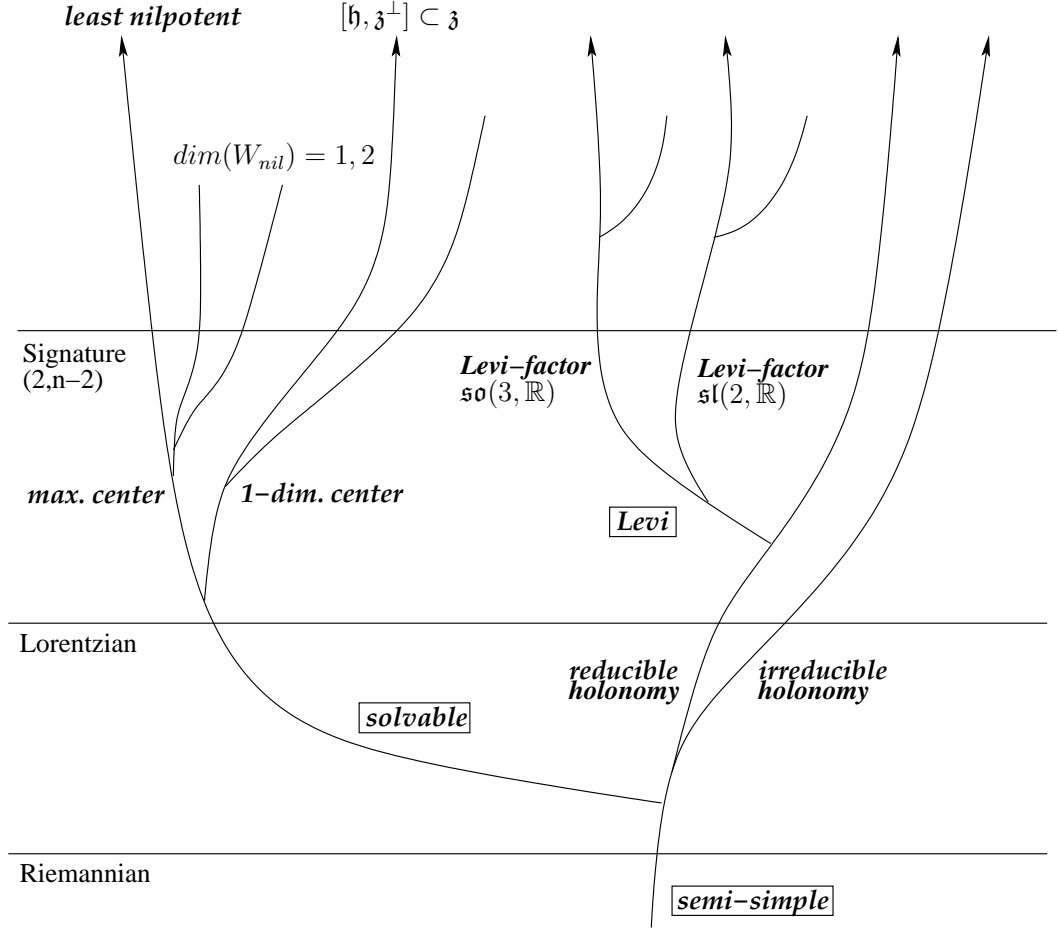}}
\caption{\it Indecomposable symmetric spaces arranged by increasing signature.
An arrow on the top of a branch denotes that the corresponding symmetric triples are classified for all
signatures. }
\end{center}
\end{figure}

The second and third part are dedicated to the solvable branch of the classification tree.
We develop an iterative  decomposition of the tangent space
which on the one hand is adapted to the metric (a generalized Witt-basis) and on the other hand 
gives the holonomy algebra the shape of certain upper triangular bloc-matrices.
Although the decomposition is not unique the dimensions of its subspaces are. These invariants yield 
the rough case differentiation made in \cite{CahenP70} while classifying solvable symmetric triples of signature $(2,n-2)$.
The adapted decomposition gives an ansatz for a solvable symmetric triple of arbitrary signature $(p,q)$
and it remains essentially to discuss the Jacobi-identity in order to get a symmetric triple.
If the dimension of the center of the transvection group is maximal
the situation simplifies. If additionally we restrict ourself to least nilpotent triples
(i.e. certain Jacobi-operators restricted to a subspace of the decomposition are invertible), their isomorphism classes can be determinated.
Although we can not guarantee that these triples are always indecomposable the class in which they are contained
is closed under decomposition. Moreover we give a sufficient criterion for indecomposability and describe those
components of the isometry group (more precisely its isotropy representation) which does not arise from the holonomy.
For the triples with maximal center which are \emph{not} least nilpotent the Jacobi-identity leads to
relations between the structure-coefficients of the transvection group which are unsolved in the general case.
For signature $(2,n-2)$ they can be managed. Only the nilpotent examples in these series has been determined in \cite{CahenP70}.
Further we cite a necessary condition for indecomposability in this class.
The case of non-maximal center is investigated in signature $(2,n-2)$. Essentially we reproduce the results of \cite{CahenP70}.
The most striking result is that the Jacobi-identity forces the dimension of a certain subspace occurring in the
iterated decomposition to be one. This might indicate that in the general case the Jacobi-identity strongly restricts the a-priory arbitrary
invariants derived from the iterated decomposition.

In the fourth chapter we return to the geometry and topology of symmetric spaces after a longer excursion to linear algebra
in the preceding parts. On the one hand it is known that the quotient manifold of a symmetric space induced by a
properly discontinuously acting group $D$ is again (globally) symmetric if and only if 
$D$ is invariant under conjugation with a symmetry $s_o$
and is contained in the centralizer $Z_{Iso(M,g)}\G(M,g)$
of the transvection group in the isometry group. This result is only partially satisfying since in general
it will be difficult to determine the whole isometry group. 
On the other hand the topological structure of affine symmetric spaces was described by S.Koh (\cite{Koh65}) without using
the representation of a symmetric space by its transvection group.
Combining both approaches we decompose a covering of symmetric spaces into one part induced by a discrete
subgroup of the center of the transvection group  (which is in principle even accessible from the corresponding Lie algebra)
and a finite covering. A refinement of the arguments shows that every solvable symmetric space is covered
finite times by a parallelizable manifold which is diffeomorphic to the product of a torus with $\R^k$.
If $k=0$, i.e. if the manifold is compact it is already flat. 
Finally we apply the above theory to calculate explicitely all solvable Lorentzian symmetric spaces up to a twofold covering.




\pagebreak

\setcounter{section}{0}
\section{Definitions and facts on symmetric spaces}
\subsection{Symmetric spaces - geometric features}

Let $\Iso=\Iso(M,g)$ denote the \mb{isometry group} of 
a symmetric space. The subgroup generated by pairs of symmetries is called the
\mb{transvection group} (also known as group of displacements):
\[
\G=\G(M,g) := \big\langle s_x \circ s_y \,\big| \, x,y, \in M \big\rangle 
\]
This name is justified by the fact that for points $x,y \in M$ lying on a smooth geodesic $\gamma$ 
the parallel displacement along $\gamma$ is realized by  $d(s_x \circ s_y)|_{\gamma}$ \cite[8.30]{O'Neill83}. 
For arbitrary points $x,y \in M$ exists at least a broken geodesic joining them and $s_x \circ s_y$ turns out to be a composition 
of displacements along the smooth extensions of its  geodesic segments \cite[I.1.4]{CahenP80}.
Moreover the composition of symmetries at the midpoints of the geodesic segments is an  isometry which maps
$x$ onto $y$. Hence $\G$ and thus $\Iso$ acts
transitively on $M$ which shows that every symmetric space is homogenous. Conjugation by the symmetry at a fixed base
point $o$ yields an automorphism of $\Iso$ which will be denoted by $\Sigma:=l_{s_o} \circ r_{s_o}$. Then for any Lie subgroup
$L \subset \Iso$ which acts transitively on $M$ and is invariant under $\Sigma$ its isotropy group $K:=\menge{g \in L}{go=o}$
satisfies
$L^{\Sigma}_0 \subset K \subset L^{\Sigma}$
where $L^{\Sigma}:=\menge{g \in L}{\Sigma(g)=g}$ and $L^{\Sigma}_0$ is its connected component of the identity
\cite[11.28]{O'Neill83}.
Differentiation leads to a \mb{symmetric decomposition} of the corresponding Lie algebra
\begin{equation}\label{eq:SymDec}
\LA(L)= \underbrace{Eig^{\sigma}(+1)}_{=\LA(K)} \oplus Eig^{\sigma}(-1) =: \k \oplus \m
\end{equation}
i.e. $[\k,\k] \subset \k$, $[\k,\m] \subset \m$ and $[\m,\m] \subset \k$
where $Eig^{\sigma}(\pm 1)$ denotes the eigenspace of $\sigma := d\Sigma_e$ w.r.t the Eigenvalues $\pm 1$. 
In many aspects it is optimal to choose as $L$ the transvection group. This is justified by the following facts:

\begin{prop}\label{prop:CharTrans}
Let $(M,g)$ be a symmetric space and $L$ a Lie subgroup of $\Iso(M,g)$ which is invariant under $\Sigma$ and acts transitively
on $M$. Then 
 \be
\item[(i)] $\G(M,g)$ is the  Lie subgroup of $L$ corresponding to the Lie subalgebra 
$\g=[\m,\m] \oplus \m$. 

\item[(ii)] $\G$ is the smallest subgroup of $\Iso(M,g)$ which acts transitively on 
$M$ and is invariant under $\Sigma$. 

\item[(iii)] The isotropy subgroup $\H$ w.r.t. $\G$ linearized at the basepoint $o$ coincides with the holonomy group 
$Hol(o) \subset Gl(T_o M)$ at $o$.
\ee
\end{prop}
 
Thus every symmetric space admits a unique representation as homogenous space 
by means of its transvection group. This representation is called \mb{prime} due to 
the property (ii) of the preceeding proposition.

\begin{proof}
(i) 
$exp(2X)=exp(X) \circ \Sigma(exp(-X))= exp(X) \circ s_o \circ exp(X)^{-1} \circ s_o=s_{exp(X)o} \circ s_o, \; X \in \m$ due to (\ref{eq:Aut=Aff}). 
Hence $\langle exp(\m) \rangle \subset \G(M,g)$. The opposite inclusion is achieved by an induction argument on the number
of smooth geodesic segments needed to connect an arbitrary  point to the basepoint $o$ (cf. \cite[I.1.4]{CahenP80}).
Thus $\G(M,g)$ is the connected Lie subgroup corresponding to the Lie subalgebra $[\m,\m] \oplus \m$ generated by $\m$.

(ii) For Lie subgroups $G \subset Isom(M,g)$ the assertion follows immediately from (i). 
In general $\forall p \in M,\;\exists g \in G$ with $g(o)=p$. Then (\ref{eq:Aut=Aff}) yields 
$s_{p} \circ s_{o}=s_{g(o)} \circ s_{o}=g \circ s_o \circ {g}^{-1} \circ s_o = g \circ \Sigma(g)^{-1} \in G$,
thus by definition  $\G(M,g) \subset G$.

(iii) is based on the fact that parallel translation along any smooth path $\gamma(t) \subset M \simeq \G/\H$ is explicitely given
by the differential of its horizontal lift to $\G$ where the horizontal subspaces are determined by left-translation of $\m \subset T_e\G$
(\cite[Thm.II.3.3]{Loos69}).
Note that the linear isotropy representation is faithful since $\G$ acts effectively on $M$. 
An alternative proof of the infinitesimal statement $\hol(o)=\h$ will be given below.
\end{proof}

\subsection{Symmetric triples - algebraic features}

One of the most striking features in the theory of symmetric spaces is that their local structure is completely 
encoded in a linearized form of their appearance. The linearization derives naturally from the chosen notion of symmetric
spaces and as the letter they are more or less equivalent. 
In view of the prime representation of a symmetric space we are lead to the following

\begin{dfn}[Symmetric triple]
\be
\item[(i)] A symmetric triple $(\g,\sigma,B)$ consist of a finite dimensional real
Lie algebra $\g$, an involutive automorphism $\sigma$ of $\g$ s.t. the symmetric decomposition $\g=\h \oplus \m$
induced by the eigenspace decomposition of $\sigma$ 
satisfies $\h=[\m,\m]$ and a non-degenerate symmetric bilinear form $B$ on $\g$ which is invariant under $\sigma$ and $\ad(\g)$, i.e.
$\sigma \in SO(B,\g)$ and $\ad(\g) \subset \so(B,\g) $.  
\item[(ii)] $\Sig(B|_{\m \times \m})=(p,q)$ is called the signature of the symmetric triple. 
\item[(iii)] An isomorphism of symmetric triples $(\g_i,\sigma_i,B_i)$ ($i=1,2$)
is a Lie algebra isomorphism $\Phi:\g_1 \rightarrow \g_2$ s.t. $\Phi \circ \sigma_1 = \sigma_2 \circ \Phi$ and 
$B_1 = \Phi^* B_2$.
\ee
\end{dfn}

Since $\sigma$ is uniquely determined by predicting its eigenspace decomposition 
$\g = \h \oplus \m$ we will refer to $(\h \oplus \m, B)$ as a symmetric triple as well. Some properties of $B$
can be reformulated as follows:

\begin{lem}\label{lem:InvBi}
\be
\item[(i)] The $\sigma$-invariance of $B$ is equivalent to $\m \perp_B \h$.
\item[(ii)] Let $(\g,\sigma)$ be a pair as in the definition above. 
Then every non-degenerate $\ad(\h)|_{\m}$-invariant symmetric
bilinear-form $B_{\m}$ on $\m$ has a unique $\ad(\g)$- and 
$\sigma$-invariant extension $B$ to $\g$.
\item[(iii)] In the situation of (ii) the following is equivalent:
\be
\item[(a)] The extension $B$ is non-degenerate, i.e. $(\g,\sigma,B)$ is a symmetric triple.
\item[(b)] The representation $\ad_{\g}(\cdot)|_{\m}:\h \rightarrow \gl(\m)$ is  faithful.
\item[(c)] $\h$ contains no non-trivial ideal of $\g$.
\item[(d)] If $G$ is a connected Lie group with $LA(G)=\g$ and $H$ a closed Lie subgroup
with $LA(H)=\h$ then the natural action of $G$ on $G/H$ is almost effective i.e. 
$N=\ker\{G \rightarrow \Diff(G/H)\}$ is discrete.
\ee  
\ee  
\end{lem}

\begin{proof} 
(i) is clear. For (ii) We define a $(4,0)$-tensor by $Q(X,Y,U,V):=B_{\m}\big( \big[[X,Y],U \big],V \big)$.
It turns out that $Q$ is a curvature tensor and as usual it can be identified
with $\wt{Q} \in S^2(\Lambda^2(\m))$ by
$\wt{Q}(X \wedge Y,U \wedge V)=Q(X,Y,U,V)$.
Then $\wt{Q}$ projects to a bilinear form $B_{\h}$ on $h$ under the map 
$\Pi: \Lambda^2{\m} \rightarrow \h$, $\sum_i X_i \wedge Y_i \mapsto \sum_i [X_i,Y_i]$
since $\ker(\Pi) \subset \ker(\wt{Q})$ and $B=B_{\h} \oplus B_{\m}$ is the required bilinear form.\\
By the $\ad(\m)$-invariance there is no other way to define $B_{\h}$ and
by the $\sigma$-invariance $\h \perp \m$, thus $B$ is unique.

(iii) (a) $\Leftrightarrow$ (b): $B$ is non-degenerate iff $B_{\h}$ is non-degenerate iff
$ B_{\h}(h,\h)=B_{\m}([h,\m],\m)=0, \; h \in \h$ implies $h=0$ iff
$ \ad(h)|_{\m}=0$ implies  $h=0$ iff $\ad_{\g}:\h \rightarrow \gl(\m)$ is faithful.

(b) $\Leftrightarrow$ (c): Let $\mathfrak{a} \subset \h$ be an ideal of $\g$. Then
$[\mathfrak{a},\m] \subset \mathfrak{a} \cap \m=\{0\}$. By the assumption b) follows $\mathfrak{a}=\{0\}$.
Conversely, assume there exists $0 \ne a \in \h$ with $[a,\m]=0$. Then
$[a,[X,Y]]=[[Y,a],X]+[[a,X],Y]=0$, hence [a,\h]=0. Thus $\R \cdot a$ would be an ideal of $\g$
contained in $\h$ which contradicts c).

(c) $\Leftrightarrow$ (d): Every connected normal Lie subgroup of $G$ contained in
$H$ corresponds to an ideal of $\g$ contained in $\h$. Thus condition (c) holds iff every normal
Lie subgroup of $G$ in $H$ is discrete. On the other hand it is easy to show that $N$ is the largest normal subgroup
of $G$ contained in $H$.
\end{proof}


The prime representation of a given symmetric space $(M,g)$ as a quotient $\G/\H$ leads naturally 
to the symmetric triple
\mb{$\tau(M,g):= \big( LA(\G) , d\Sigma_e , (d\pi|_{\m})^* g_o \big)$}
where $\pi: \G \rightarrow \G/\H \simeq M$ denotes the canonical submersion.  
By proposition \ref{prop:CharTrans} the identity $\h=[\m,\m]$ is satisfied and $(d\pi|_{\m})^* g_o $ 
extends to a non-degenerate bilinear form on whole $\g$ according to lemma \ref{lem:InvBi}.
Note that $\tau(M,g)$ actually depends on the chosen base point in $M$.

Conversely, let $(\g,\sigma,B)=(\h \oplus \m,B)$ be a symmetric triple. 
Consider the simply connected Lie group $\wt{\G}$ with $ LA(\wt{\G})=\g $. 
Hence $\sigma \in Aut(\g)$ admits a lift to an involutive Lie group automorphism
$\wt{\Sigma}\in Aut(\wt{G})$. Furthermore $\wt{\H}:=\wt{\G}^{\wt{\Sigma}}$ is closed
and connected due to \cite[p.293]{Koh65}\footnote{
Alternatively one could define $\wt{\H}:=\wt{\G}^{\wt{\Sigma}}_0$.}
and thus $\wt{M}=\wt{\G}/\wt{\H}$ is simply connected due to the exact homotopy sequence
$1=\pi_1(\wt{G}) \rightarrow \pi_1(\wt{\G} / \wt{\H}) \rightarrow \pi_0(\wt{\H})=1$.
Finally one observes that $B \in S^2(\m)$ is $Ad(\wt{H})$-invariant due to the connectedness of $\wt{H}$.
Then $\wt{\Sigma}$ equips $\wt{\G} / \wt{\H}$ with symmetries w.r.t. the pseudo-Riemannian
metric $g^B$ gained by left-translation of $(d\pi|_{\m})_* B$
(see \cite[11.22/29]{O'Neill83}). The so constructed simply connected symmetric spaces
will be denoted by 
\mb{$\theta(\g,\sigma,B) := \big( \wt{\G}/\wt{\H},g^B \big)$}.

\begin{thm}\label{thm:Cor}
The above construction establishes  a bijection between isometry classes of simply connected symmetric spaces and
isomorphism classes of symmetric triples. In particular the maps $\tau$ and $\theta$ are compatible with isomorphisms resp. isometries.
\end{thm}

This result is well known (see e.g. \cite[Ch.I \S2]{CahenP80}, \cite[II.4.12]{Loos69} or \cite[I.2.9]{Bertram00}). 
As will be shown  below the map $(M,g) \mapsto \tau(M,g)$ can be understood as the evaluation of the curvature tensor.
Thus it is possible to verify that the symmetric triples $(\g,\sigma,B)$ and $\tau\big(\theta(\g,\sigma,B)\big)$
are isomorphic. However it is worth to mention that the group $\wt{\G}/N$ is indeed the transvection group 
of the symmetric space $\theta(\g,\sigma,B)$ with the discrete group $N$ as in lemma \ref{lem:InvBi}.
In general $N \neq \{e\}$ as the example of the sphere $S^2$ with the standard metric shows.    

The differential geometric manifestation of symmetric spaces as manifolds with parallel curvature provides
a slightly different object as linearization, the so-called {\bf Lie triple system (Lts)}. It associates
to a locally symmetric space the trilinear map $[\cdot,\cdot,\cdot]: T_oM^3 \ra T_oM, \; [X,Y,Z]:=\RR_o(X,Y)Z$
given by its curvature tensor and the metric $g_o$ on a fixed tangent space $T_oM$ 
(for an abstract definition of a Lts see e.g. \cite[p.78]{Loos69}). 

Now any Lts admits an embedding into a Lie algebra s.t. its image carries the structure of an symmetric triple
and conversely the Lts is regained by the trilinear map $[[\cdot,\cdot],\cdot]:\m^3 \ra \m$. More precisely consider
\begin{align*}
Aut(T_oM,\RR_o,g_o):=&\menge{ A \in O(T_oM,g_o) }{ \RR_o\big( A(\cdot), A(\cdot) \big) \circ A = A \circ \RR_o (\cdot,\cdot) } \\
\Rightarrow \quad \aut:=\aut(T_oM,\RR_o,g_o) =& 
\menge{ a \in \o(T_oM,g_o)}{ \big[a,\RR_o(\cdot,\cdot) \big] - \RR_o \big( a(\cdot),\cdot \big) - \RR_o \big( \cdot, a(\cdot)\big) = 0 }
\end{align*}
Then $\RR$ being parallel enables the definition of the following Lie bracket on the vector space $\aut \oplus T_oM$
\begin{align}\label{eq:Bracket}
& [\cdot,\cdot]:T_o \times T_o \rightarrow \aut : &&[X,Y]=\RR_o(X,Y) \notag\\
& [\cdot,\cdot]:\aut \times T_o \rightarrow T_o: && [R,Z]=R(Z) \\
& [\cdot,\cdot]:\aut \times \aut \rightarrow \aut: && [R,S]=R \circ S - S \circ \notag R
\end{align}
(cf. \cite[p.229]{KobayashiN1}) and $\Re_o \oplus T_oM = \h \oplus \m$ becomes a symmetric triple called the {\bf standard embedding}
of the Lts $(T_oM,\RR_o,B_o)$, where $\Re_o := span\menge{\RR_o(X,Y)}{X,Y \in T_oM} \subset \aut$. 

On the other hand the theory of natural reductive homogenous spaces applies to the prime representation of a symmetric space
and yields under the identification $\m \simeq T_oM$ the curvature formula $\RR_o(X,Y)=\ad([X,Y])|_{\m}$
(\cite[prop.11.31]{O'Neill83}). Hence the symmetric
triple $\tau(M,g)$ and the standard embedding of the Lts $(T_oM,\RR_o,g_o)$ are isomorphic for a (globally) symmetric space $(M,g)$. 
Analogous any symmetric decomposition $LA(L) = \k \oplus \m$ w.r.t a Lie subgroup $L \subset \Iso(M,g)$ can be comprehended as 
subalgebra of $\aut \oplus T_oM$ with the Lie bracket given by (\ref{eq:Bracket}).

Note that the infinitesimal version of proposition \ref{prop:CharTrans}(iii) derives now alternatively from the Ambrose-Singer-theorem 
since $\hol(o)=span \menge{\R_o(X,Y)}{X,Y \in \m} = \ad([\m,\m])=\ad(\h) \simeq \h$ by the faithfulness of 
the isotropy representation $\ad(\cdot)|_{\m}:\h \ra \gl(\m)$.

\begin{prop}\label{prop:Aut} 
Let $(M,g)$ be a simply connected symmetric space and $(T_oM,\RR_o,g_o)$  its associated Lie triple system.
\be
\item[(i)] The symmetric decomposition of the isometry algebra is given by
\[ 
 \iso(M,g) \simeq \aut(T_oM,\RR_o,g_o) \oplus T_oM
\]
w.r.t the Lie algebra structure of (\ref{eq:Bracket}). 

\item[(ii)] Any automorphism of the Lts $(T_oM,\RR_o,g_o)$ extends uniquely to a Lie algebra automorphism 
of its standard embedding $\Re_o \oplus T_oM \simeq \g(M,g)=\g$ which commutes with $\sigma$, i.e
\[
Aut(T_oM,\RR_o,g_o) \simeq \Zent_{\ds Aut(\g)} \langle \sigma \rangle
\]
\ee
\end{prop}

\begin{proof}
(i) By virtue of proposition \ref{prop:CharTrans}(i) it suffices to determine the isotropy part of the isometry algebra. 
By \cite[Thm.8.17]{O'Neill83} 
it is maximal in the sense that the necessary condition of 
leaving the metric and curvature invariant is already sufficient.

It should be emphasized that the assumption $M$ being simply connected cannot be omitted although the assertion is formulated on the
Lie algebra level. An easy counterexample is the flat 2-dimensional cylinder. 
Thus the dimension of the isometry group can be diminished by a covering whereas
every covering of symmetric spaces yields a covering of the transvection groups (cf. proposition \ref{prop:SymCov}).

(ii) Let $A \in Aut(T_oM,\RR_o,g_o)$ and define its extension
to $\Re_o=[T_oM,T_oM]$ by 
$\wt{A}\big( \sum_i [X_i,Y_i] \big) := \sum_i \big[ A(X_i),A(Y_i) \big]$, $X_i,Y_i \in T_oM$.
$\wt{A}$ is well defined iff $R=\sum_i [X_i,Y_i]=0$ implies $\wt{A}(R)=0$:
\[
\big[ \wt{A}(R),A(Z) \big]= \Big[ {\sum}_i [A(X_i),A(Y_i)] , A(Z) \Big] = A \big( \Big[ {\sum}_i [X_i,Y_i] , Z \Big] \big)= A\big([R,Z]\big)=0
\]
Since the adjoint representation of $\Re_o$ on $T_oM$ is faithful (\ref{lem:InvBi}(iii)) and $A(T_oM)=T_oM$ it follows $\wt{A}(R)=0$.  
The proof also shows that the assertion actually can be generalized to surjective homomorphisms of Lts's or homomorphism between
semi simple Lts's (cf. \cite[V.1.9]{Bertram00}). 
\end{proof}

To conclude this paragraph we recall further relations between symmetric spaces and symmetric triples.

\begin{dfn}[Decomposability] 
\be
\item[(i)] A connected pseudo-Riemannian manifold $(M,g)$ is said to be decomposable if there
exists a proper subspace $V \subset T_oM$  which is non-singular w.r.t. $g_o$ (i.e. $g_o|_{V \times V}$ is non-degenerate) 
and invariant under $Hol(o)$. Otherwise $(M,g)$ is indecomposable\footnote{
In \cite{CahenP80},\cite{Wu67} the equivalent notion is \emph{weakly irreducible}.}.
\item[(ii)] A symmetric triple $\tau=(\h \oplus \m,B)$ is said to be decomposable if 
$\m$ contains a proper subspace $\m'$ which is non-singular w.r.t. $B$ 
and invariant under $\ad(\h)$. Otherwise $\tau$ is indecomposable.

A symmetric triple $(\h \oplus \m,B)$ is said to be {\bf pseudo-Euclidian} if $\h=0$,
hence $\g=\m$ is abelian.
\ee
\end{dfn}

By the theorem of {\bf de Rham-Wu} (\cite{Wu64}) a decomposable pseudo-Riemannian manifold which is additionally simply connected and complete
is isometric to a pseudo-Riemannian product manifold. Furthermore a pseudo-Riemannian product is symmetric iff its factors
are symmetric. Hence for the purpose of a classification of simply connected symmetric spaces we can restrict ourself to indecomposable ones.
(for a discussion of ambiguity in the de Rham-Wu decomposition see \cite[Appendix I, example 3]{Wu67}, \cite[I.4.8]{CahenP80}).
In order to relate the decomposition of symmetric spaces to those of symmetric triples via the correspondence of theorem \ref{thm:Cor} we make the
following preliminary remark: 

Consider two symmetric triples $\tau_i=(\g_i,\sigma_i,B_i),\;i=1,2$. Then their direct sum is defined as
$\tau_1 \oplus \tau_2 := (\g_1 \oplus \g_2, \sigma_1 \oplus \sigma_2, B_1 \oplus B_2)$
where $\g_1 \oplus \g_2$ is the direct sum of Lie algebras, $(\sigma_1 \oplus \sigma_2)(X_1,X_2) =
(\sigma_1(X_1),\sigma_2(X_2))$ and $ B = B_1 \oplus B_2$ means $B|_{\g_i \times \g_i}=B_i$ and
$\g_1 \perp \g_2$. Clearly $\tau_1 \oplus \tau_2$ is decomposable. 
Conversely, if $(\h \oplus \m,B)$ is decomposable there exists $\m_i \subset \m$ s.t. 
$\m=\m_1 \oplus \m_2$, $\m_1 \perp_B \m_2$, $\ad(\h)\m_i \subset \m_i$ and thus
$(\h \oplus \m) \simeq ([\m_1,\m_1]\oplus\m_1,B_1) \oplus ([\m_2,\m_2] \oplus \m_2,B_2)$ where
$B_i$ denotes the restriction of $B$ to $[\m_i,\m_i] \oplus \m_i$.   

\begin{prop}\label{prop:Decomposability}
Let $(M,g,o)$ be a simply connected symmetric space. Then $M$ is indecomposable if and
only if $\tau(M,g)$ is indecomposable. Furthermore
\begin{align*}
\tau(M_0 \times \ldots \times M_r,g_0 \times \ldots \times g_r) \;&\underset{\text{isomorph}}{\simeq}\:
\tau(M_0,g_0) \oplus \ldots \oplus \tau(M_r,g_r) \\
\theta(\tau_0) \times \ldots \times \theta(\tau_r) \; &\underset{\text{isometric}}{\simeq} \;
\theta(\tau_0 \oplus \ldots \oplus \tau_r)
\end{align*}
where all $(M_i,g_i)$ are simply connected symmetric spaces and $\tau_i$ are symmetric
triples.
\end{prop}

The following lemma describes algebraically the space of parallel vector fields ($M$ simply connected), 
i.e. those elements in $\m$ which are in the kernel of the adjoint representation $\ad(\h)|_{\m}$
due to the holonomy principle.

\begin{lem}\label{lem:Center}
Let $\tau = (\h \oplus \m,B)$ be a symmetric triple. Then:
\begin{enumerate}
\item[(i)]
 $\big \{ X \in \m \,|\, [\h,X]=0 \big\} = \z(\g) = [\g,\g]^{\perp} = [\h,\m]^{\perp} \cap \m$
where $\z(\g) $
is the center of the Lie algebra $\g$.
\item[(ii)]If $\g$ is solvable its center $\z(\g)$ is non-zero.
\item[(iii)]If $\tau$ is indecomposable and $\dim(\m) > 1$ then $\z=\z(\g)$ is totally isotropic
(i.e. $B|_{\z \times \z}=0$).
\end{enumerate}
\end{lem}

\begin{proof}
(i) Using the $\ad(\g)$-invariance of $B$ we have the following equivalences:
\[
X \in \z(\g) \quad \Leftrightarrow  \quad  B([X,U],V)= B(X,[U,V])=0,\; \forall U,V \in \g \quad \Leftrightarrow  \quad X \in [\g,\g]^{\perp}
\]
This results holds in general for a metric algebra (cf. \cite[prop.2.1.]{Bordemann97}). Furthermore
\begin{align*}
&[\g,\g]=[\h \oplus \m,\h \oplus \m]= \big([\h,\h] + [\m,\m] \big) \oplus [\h,\m]=\h \oplus [\h,\m]\\
\Rightarrow \quad  &[\g,\g]^{\perp}= \big( \h \oplus [\h,\m] \big)^{\perp} = \h^{\perp} \cap [\h,\m]^{\perp}
=\m \cap [\h,\m]^{\perp}
\end{align*}
and thus $B|_{\m \times \m}$ being non-degenerate yields
\[
X \in \m \cap [\h,\m]^{\perp}  \quad \Leftrightarrow  \quad X \in \m \text{ and } B(X,[\h,\m])= B([X,\h],\m)=0
    \quad \Leftrightarrow  \quad X \in \m \text{ and } [\h,X]=0
\]
(ii) $\g$ being solvable implies $[\g,\g] \subsetneq \g$. Hence, by (i) $\z(\g)=[\g,\g]^{\perp} \ne \{0\}$.

(iii) Assume $\z$ is not totally isotropic. Then there exists a vector $Z  \in \z$ with
$B(Z,Z) \ne 0$. Thus $\R \cdot Z$ is a non-zero, non-singular and holonomy-invariant subspace since
$\ad(\h)Z=0$. From the indecomposability of $\tau$ follows $\m=\R \cdot Z$ and thus $\dim(\m)=1$ which 
has been excluded.  
\end{proof}

Finally the algebraic distinction of Lie algebras is reflected in certain geometric features:

\begin{prop}\label{prop:Levi} 
Consider an indecomposable symmetric space $(M,g)$ with $\dim(M) >1$.  
\be
\item[(i)] The holonomy $\hol(o) = \ad(\h)|_{\m}$ acts completely reducible (resp. nilpotent) if and only if $\g$ is semi-simple (resp. solvable).
\item[(ii)] The Ricci-tensor $\RR ic$ is non-degenerate if and only if $\g$ is semi-simple. Otherwise it is two-step nilpotent i.e. $\RR ic^2=0$.
\item[(iii)] If $\g$ admits a proper Levi-decomposition $\g = \rad(\g) \rtimes \s$ the signature $\Sig(g)=(p,q), \;(p \leq q)$ 
restricts the following dimensions
\[
\dim(\s) \leq \frac{p(p-2)}{2} \qquad \dim(\z(\g)) = p-2
\]
\ee
\end{prop}

\begin{proof} For (i) we refer to \cite[prop.I.5.4]{CahenP80}. The first part of (ii) is well known (cf. \cite[prop.VI.1.3]{Loos69}). Note
that for a pseudo-Riemannian manifold the Ricci-tensor is always symmetric. Further it has been shown in \cite{BoubelB01} that a degenerate
Ricci-tensor of an indecomposable Ricci-parallel  pseudo-Riemannian manifold is already two-step nilpotent.
(iii) is a consequence of a rigorous investigation
of symmetric triples with a non-trivial Levi-decomposition achieved in \cite[II.\S1]{CahenP80} (In particular cor.1.7).  
\end{proof}

\begin{re} A necessary condition for the existence of a parallel spinor-field is that the Ricci-tensor is 2-step nilpotent
(cf. \cite[Satz 4.3]{Baum81}).
\end{re}

\begin{cor}\label{cor:SeperationProperty}
Let $\tau=(\g,\sigma,B)$ be an indecomposable symmetric triple.
\begin{enumerate}
\item[(i)]If $\tau$ is Riemannian, $\g$ is semi-simple or one-dimensional.
\item[(ii)]If $\tau$ is Lorentzian, $\g$ is either semi-simple or solvable.
\item[(iii)]If $\tau$ has signature $(2,n-2)$, $\g$ can be semi-simple, solvable, or
proper Levi with Levi factor $\s \in \{\so(3,\R),\sl(2,\R)\}$.

\item[(iv)]If $\tau$ is Lorentzian or of signature $(2,n-2)$ or has maximal center i.e $\dim(\z(\g))=p$ then:
\[
\g \text{ is solvable } \quad \Leftrightarrow \quad \z(\g) \ne \{0\}
\]
\end{enumerate}
\end{cor}

The semi simplicity of the transvection algebra of a Riemannian symmetric space has been used already
by Cartan for their classification. The ``separation property'' in the Lorentzian case,
i.e. the division into solvable and semi-simple transvection algebras was found in \cite[Theorem 3]{CahenW70}. 
In the same paper it was already observed that this result
does not hold for higher signatures. A sketch of the proof for assertion (iii) can be found in
\cite[Lemma 5.1.]{CahenP70}. Finally all statements can be derived from the inequalities stated in proposition \ref{prop:Levi}(iii)
and lemma \ref{lem:Center}.

\pagebreak
\addtolength{\evensidemargin}{-.5cm}
\addtolength{\oddsidemargin}{-.5cm}
\addtolength{\textwidth}{1cm}
\addtolength{\textheight}{1cm}
\addtolength{\footskip}{1cm}

\begin{setlength}{\hfuzz}{100pt}
\begin{samepage}
\begin{center}
\subsection[Tabular summary]{\sloppy Correspondence between geometrical and algebraical facts of symmetric spaces and symmetric triples}
{\small
\fbox{\fbox{
\rca{1.8}{
\btab{c}
\btab{  p{70mm} p{5mm}  p{70mm} }
\mb{$\Iso$} - isometry group   && \mb{\iso = LA(\Iso) } \\
\mb{$\G$} - transvection group  && \mb{$\g = LA(\G)$}  \\
\mb{$\Sigma = C_{s_o}$} conjugation by the symmetry $s_o$ && \mb{$\sigma = d\Sigma_{id}$} \\
\hline
\mb{$\H=Stab_o \G$} and $\G^{\Sigma}_0 \subset \H \subset \G^{\Sigma}$ && 
\mb{$\g = \h \oplus \m$} with $\h = LA(H)=Eig^{\sigma}(+1)$ and $\m = Eig^{\sigma}(-1)$\\ 
\mb{$Ad(\H)|_{\m}$} linear isotropy representation && \mb{$\ad(\h)|_{\m}$} the corresponding Lie algebra repr. \\
\mb{$\RR_o(X,Y)Z=\kov_{[X,Y]}Z-[\kov_X,\kov_Y]Z$} &\mb{$=$}& \mb{$[[X,Y],Z]$} \\
\mb{$\RR ic_o(X,Y)=tr_{\m}(Z \mapsto \RR_o(X,Z)Y)$} &\mb{$=$}& \mb{$-\frac{1}{2} \cdot \B^{\g} $} Cartan-Killing-form\\
\hline
\multicolumn{3}{c}{\mb{one-to-one correspondence:}} \\
$(M,g)$ a symmetric space with its prime representation $(\G/\H,\Sigma,g)$ 
& $\rightarrow$ & \mb{$\tau(M,g)= \big( LA(\G),d\Sigma_e,g_o \big) $} the associated symmetric triple \rule{0mm}{-4mm}\\ 
\mb{$\theta(\g,\sigma,B)=\big( \wt{\G}/\wt{H},\wt{\Sigma},g^B \big)$}  the associated simply connected symmetric space
& $\leftarrow$ & $\tau=(\g,\sigma,B)=(\h \oplus \m , B)$ a symmetric triple\\ \hline
$\G$  is the smallest $\Sigma$-stable subgroup of $\Iso(M,g)$ acting transitively on $M$ &&
$\g=[\m,\m] \oplus \m$ w.r.t. any symmetric decomposition $\k \oplus \m$\\ \hline
\mb{ $\H \simeq Ad(\H)|_{\m} = Hol(o)$} the holonomy group at $o$
&& \mb{$\h \simeq \ad(\h)|_{\m}=\hol(o)$} due to $\h = [\m,\m]$ and the theorem of Ambrose-Singer\\ \hline 
%
\mb{$LA \big( Stab_o(Iso)  \big) \simeq \aut(T_oM,\RR_o,g_o)$} for $M$ simply connected. && 
\raisebox{-1mm}{\mb{$ \Zent_{Aut(\g)} \langle \sigma \rangle = Aut \big(\m,[[\cdot,\cdot],\cdot],B \big) $}} \\ \hline
%
%
%
\multicolumn{3}{c}{\mb{de Rham-Wu decomposition:}}\\
$\theta(\tau_1 \oplus \tau_2) \underset{isometric}{\simeq} \theta(\tau_1) \times \theta(\tau_2) $&&
{\scriptsize  $\tau(M_1 \times M_2,g_1 \times g_2)  \underset{isomorph}{\simeq} \tau(M_1,g_1) \oplus  \tau(M_1,g_1)$ \stru{6}} \\ \hline
space of \mb{parallel vector fields} && \mb{$\z(\g)=[\g,\g]^{\perp}=[\h,\m]^{\perp}\cap \m$} \\ \hline
\etab
\\
\btab{p{8mm}  p{62mm} p{5mm}  p{65mm} }
& $\hol(o)$  acts \mb{$\begin{cases} \text{completely reducible} \\ \text{nilpotent} \end{cases}$} \stru{11}& \mb{$\Leftrightarrow$} & 
$\g$ is \mb{$\begin{cases} \text{semi-simple} \\ \text{solvable} \end{cases}$} \stru{11}\\ \cline{2-4}
& \mb{$\Sig(g)=(p,q)$}, $(p \leq q)$ signature
of the pseudo-Riemannian metric $g$ && If $\g=\rad(\g)\rtimes \s$ is a proper Levi decomp.
 \mb{$\dim(s)\leq \frac{p(p-2)}{2}$} and \mb{$\dim(\z(\g))\leq p-2$}\stru{7}  \\ \cline{2-4}
& $(M,g)$ Riemannian $\Rightarrow$ $\G=\Iso_0$ && 
 $\g$ semi-simple $\Rightarrow$ any derivation is inner\\ \cline{2-4}
& $\RR ic$ is \mb{non-degenerate} i.e. $(M,g)$ is Einstein w.r.t. $g \sim \RR ic$ & $\Leftrightarrow$ & $\g$ is \mb{semi-simple} \\
& \mb{$\RR ic^2=0$} i.e.  $(M,g)$ is Einstein iff it is Ricci-flat  & $\Leftrightarrow $&
 \mb{$\rad(\g) \neq \{0\}$} 
\etab
\etab
}}}}
\begin{picture}(0,0) \setlength{\unitlength}{1mm}
 \put(-61,33){\rotatebox{90}{\raisebox{0mm}[0mm][0mm]{\makebox[0mm]{$M$ indecomposable and $\dim(M)>1$ }}}} \end{picture}
\end{center}
\end{samepage}

\end{setlength}

\addtolength{\textheight}{-1cm}

\pagebreak
\addtolength{\textwidth}{-1cm}
\addtolength{\evensidemargin}{.5cm}
\addtolength{\oddsidemargin}{.5cm}
\addtolength{\footskip}{-1cm}


\section{ \sloppy Solvable Symmetric Triples}

In this part we present an approach to uncover some structure and invariants of indecomposable solvable symmetric triples.
The initial position differs completely from the semi-simple case. 
On the one hand a structure theory for solvable Lie algebras comparable to that of semi-simple ones is not available 
- not to mention a classification. On the other hand there is no intrinsic $\ad$-invariant scalar product on a solvable Lie algebra as
the Cartan-Killing form in the semi-simple case - there might be even no one.
Thus it seems to be advantageous to treat the algebraic and metric structure of a symmetric triple simultaneously as pointed out in the introduction.

Our first step in this direction is to decompose a  pseudo-Euclidian vector space adapted to a predicted (singular) subspace
which will be in the following the image of the action of the holonomy algebra $[\h,\m] = [\g,\g] \cap \m$
or algebraically spoken the first derivative of the transvection algebra. 
The stabilizer $\Gamma \subset O(\m,B)$ of $[\h,\m]$ and $B$ contains the automorphism group $Aut=Aut\big(\m,[[\cdot,\cdot],\cdot],B\big)$ 
(\ref{lem:SkewInIso}) which is interesting from a geometric point of view (\ref{prop:Aut}).
In order to find the invariants and a good normal form of a symmetric triple the action of $\Gamma / Aut$ will be studied
parameterizing all adapted decompositions as graphs over a given one. 

In general this procedure has to be iterated considering the image of the action of $\h$ on 
the quotient $[\h,\m]/[\h,\m]^{\perp}=\z^{\perp}/\z$
and so forth. The depths of the iteration is bounded by the signature of the symmetric triple.
This construction is also related to double extensions (cf. \cite[3.3]{Neukirchner02}).

\subsection[\sloppy Decomposition of a pseudo-Euclidian vector space adapted to a subspace]
{\sloppy Decomposition of a pseudo-Euclidian vector space adapted to a predicted subspace}\label{ss:AdapDecomp}

Given a real (finite dimensional) vector space $V$ equipped with a non-degenerate symmetric bilinear-form
$B$ and a fixed subspace $F \subset V$, it is natural to ask for a decomposition of V into a direct
sum of subspaces which is adapted to the predicted data $F$ and $B$. The freedom in choosing such a decomposition
and moreover an adapted basis in such a decomposition is measured by the group $\Gamma \subset GL(V)$ which preserves $B$ and $F$:
\[
\Gamma = stab(F) \cap O(V,B) = \big\{Q \in O(V,B)\;|\;Q(F)=F \big\}
\]

\begin{lem}[generalized Witt-decomposition]\label{lem:WittDecomp}
Let $V$ be a vector space, $B$ a symmetric non-degenerate bilinear form on $V$ and $F \subset V$ a linear subspace. 
Then $V$ admits a decomposition into $V = E \oplus W \oplus U \oplus E^*$ s.t.
\begin{align*}
&E =F \cap F^{\perp}   && B|_{E \times E}=B|_{E^* \times E^*}=0\\
&E \oplus W =F          && U \perp W\\
&E \oplus U =F^{\perp}  && (E \oplus E^*) \perp (U \oplus W)
\end{align*}
Every such decomposition is said to be \mb{adapted to $F$ and $B$} - abbreviated by \mb{$(F,V,B)$-adapted}.
The family of adapted decompositions will be denoted by \mb{ $\mathcal{D}=\mathcal{D}(F,V,B)$}.
\end{lem}

\begin{proof}
Consider the following diagram of natural inclusions:
\[
\rca{1.3}{\begin{array}{ccc}
F \cap F^{\perp} &\subset &F\\
\cap & {} & \cap\\
F^{\perp} &\subset &V
\end{array}}
\]

\pagebreak

First choose supplementary subspaces $W$ (resp. $U$) of $E:=F \cap F^{\perp}$ in $F$ (resp. $F^{\perp}$). Then:
\begin{itemize}
\item $W \perp U$, since $W \subset F \perp F^{\perp} \supset U$.
\item $W,U$ are non-singular: e.g. $W \cap W^{\perp}=W \cap F \cap W^{\perp}=W \cap (F^{\perp} + W)^{\perp} \subset W \cap F \cap F^{\perp}=\{0\}$.
\item $W \cap U =\{0\}$ due to the above properties.
\end{itemize}
Thus $V = (W \oplus U) \oplus  (W \oplus U)^{\perp}$ and by the dimension formula $\dim(W \oplus U )^{\perp} = 2 \cdot \dim(E)$ since
$\dim(W \oplus U )= \Big( \dim(F)-\dim(E) \Big)+ \Big( \dim(F^{\perp}) -\dim(E) \Big) = n-2 \cdot \dim(E)$. Set $\dim(E)=:p$.

$E=F \cap F^{\perp}=(F + F^{\perp})^{\perp} \subset (W \oplus U )^{\perp}$ is a totally isotropic subspace.
Hence $\Sig \big( (W \oplus U )^{\perp} \big)= (p,p)$. Thus by the theorem of Sylvester
$ \big( (W \oplus U )^{\perp},B|_{(W \oplus U )^{\perp} \times (W \oplus U )^{\perp}} \big)$ is isometric to
$\big( \R^{2p},\langle \cdot , \cdot \rangle_p \big)$. On the other hand $\R^{2p}=E^+ \oplus E^-$, where
$E^{\pm}=span\{f_i^{\pm} \;|\; f_i^{\pm} =e_i \pm e_{p+i},\;1 \le i \le n \}$ are two totally isotropic subspaces of $\R^{2p}$.
Obviously $E^+$ can be mapped isometrically onto $E \subset V$ and according to the theorem of Witt it admits an extension
to an isometry from $\R^{2p}$ into $V$. The image of $E^-$ under this isometry obviously satisfies the required properties for $E^*$.
\end{proof}

\begin{re}\label{re:AdapBasis}
\be
\item[(i)] $E,E^*$ are said to be \mb{singular paired w.r.t. $B$}, i.e. $E,E^*$ are both totally isotropic and 
$B|_{E \times E^*}$ is non-degenerate.
\item[(ii)] If $F^{\perp}$ is totally isotropic $E=F^{\perp} \subset F$, hence $U=\{0\}$. If $F^{\perp}$ is even maximal
isotropic, $B_{W \times W}$ is definite (cf. (iii)) and the decomposition $V=E \oplus W \oplus E^*$ is known as \mb{Witt decomposition}.
\item[(iii)] $\Sig(W)+\Sig(U)+ \big( \dim(E),\dim(E) \big)=\Sig(V)$
\ee
\end{re}

\begin{dfn}[Adapted Basis]\label{dfn:AdapBasis} 
Choose a pseudo-orthonormal basis $\beta_U$ (resp. $\beta_W$) of $U$ (resp. $W$) and an arbitrary
basis $\beta_E$ of $E$. Further associate to $\beta_E=\{b_i\}$ its dual basis $\beta_E^*=\{b_i^*\}$ in $E^*$ i.e.
$B(b_i,b_j^*)=\delta_{ij}$. Then the matrix of $B$ w.r.t. the basis $\beta=\{\beta_E,\beta_W,\beta_U,\beta_E^*\}$ receives the form
\vspace{.5cm}
\begin{equation}\label{mat:Adapted}
Mat_{\beta}(B)=
\rca{1.3}{\left( \ba{c|cc|c} 
\nosp{0}{5}{$E$} & \nosp{0}{5}{$W$} & \nosp{0}{5}{$U$} & \nosp{0}{5}{$E^*$} I_p \\ \hline 
&I_{\Sig(W)} && \\
&& I_{\Sig(U)}&\\ \hline
I_p &&&     
\ea \right)}
\end{equation}
where $p=\dim(E)$, $I_p$ is the $p$-dimensional identity matrix and
$I_{(p,q)}=\left( \begin{smallmatrix}-I_p & 0 \\ 0 & I_q \end{smallmatrix} \right)$.
Every basis $\beta$ of $V$, such that $Mat_{\beta}(B)$ is of the form (\ref{mat:Adapted}) w.r.t.
an adapted decomposition $D=(E,W,U,E^*)\in \mathcal{D}$
is said to be an \mb{$(V,B,D)$-adapted basis}. 
The family of such basis will be denoted by \mb{ $\mathcal{B}_D$}.
\end{dfn}

\begin{re}\label{re:AdapBasisII}
\be
\item[(i)] $\Gamma = stab(F) \cap O(V,B)$ acts simply transitively on $\{ \mathcal{B}_D \}_{D \in \mathcal{D}}$. 
\item[(ii)] Since $E = \ker(B|_{F \times F})=\ker(B|_{ F^{\perp} \times  F^{\perp}})$, $B$ transfers onto the quotients
$F/E$ and $ F^{\perp}/E$. Thus there exist canonical isometries $(W,B|_{W \times W}) \simeq (F/E,B_{F/E})$ and
$(U,B|_{U \times U}) \simeq (F^{\perp}/E,B_{F^{\perp}/E})$. 
Hence any basis
$\beta=\{\beta_E,\beta_W,\beta_U,\beta_E^*\}\in \mathcal{B}_D$ is uniquely determined by the basis
$\beta_E$ of $E$ (which fixes $\beta_E^*$ by duality) and the pseudo-orthonormal basis $\beta_{F/E}$
(resp.$\beta_{F^{\perp}/E}$) in $(F/E,B_{F/E})$ (resp. $(F^{\perp}/E,B_{F^{\perp}/E})$) corresponding
to $\beta_W,\beta_U$ under the above isometries.
\item[(iii)] There is a bijection between $(V,B,F)$-adapted decompositions and the set \mb{ $\mathcal{E}^*$} of
totally isotropic supplementary subspaces $E^*$ of $F + F^{\perp}$ in $V$ (cf. \cite[page 69]{Boubel00}).
\end{enumerate}
\end{re}

\pagebreak

For computational purposes later on we need a concrete realization of all adapted decompositions and basis.
It will be obtained by investigating the degrees of freedom in the construction considered in the proof of lemma \ref{lem:WittDecomp}.
Let $D =(E,W,U,E^*) \in \cal{D}$ be a fixed adapted decomposition.

\begin{description}\label{dfn:Rho}
\item[\mb{ \raisebox{1mm}{$\boxed{\rho_{\wt{U},\wt{W}}}$} }]
The first freedom lies in the choice of the complements $W,U$. 
We denote the family of supplementary subspaces of $E$ in $F$ (resp. $F^{\perp}$) by \mb{$\mathcal{W}$} (resp. \mb{$\mathcal{U}$}).
Any such $\wt{W}\in \cal{W}$ may be  considered as graph of a linear map $L:{W} \rightarrow E$ which establishes the bijection
\begin{align*}
Hom({W},E)\quad &\rightarrow \quad \mathcal{W}\\
L_{\wt{W}} \quad &\mapsto \quad \wt{W}=\Gamma(L_{\wt{W}}) :=\menge{\big(w,L_{\wt{W}}(w) \big) \in W \oplus E=F }{ w \in W}
\end{align*}
The inverse of this map is obviously given by $\mathcal{W} \ni \wt{W} \mapsto -pr^{\wt{W} \oplus E}_2|_W$. 
Analogous we construct a bijection $Hom({U},E) \leftrightarrow \mathcal{U}$.
We assign to each $\wt{W} \in \mathcal{W},\;\wt{U} \in \mathcal{U}$ the map
$\wt{\rho}_{\wt{U},\wt{W}} \in End(F +F^{\perp})=
Hom(E \oplus {W} \oplus {U},E \oplus \wt{W} \oplus \wt{U})$ defined by
$\wt{\rho}_{\wt{U},\wt{W}}(e+w+u)= \big( e+w+L_{\wt{W}}(w)+u+L_{\wt{U}}(u) \big)$
which is an isometry on $F+F^{\perp}$
since $E \in \ker(B|_{(F+F^{\perp}) \times(F+F^{\perp})})$. 
Thus there exists an extension to an isometry  \mb{$\rho_{\wt{U},\wt{W}} \in O(V,B)$} due to the theorem of Witt.

\item[\mb{ \raisebox{1mm}{$\boxed{\rho_{ \wt{E^*} } }$ } }]
The second ambiguity concerns the choice of a Witt complement $E^*$ of $E$ in $({W} \oplus {U})^{\perp}$.
Again every such complement $\wt{E^*}$ can be comprehended as a graph of an element in $Hom(E^*,E)$
\begin{align*}
Hom({E^*},E)\quad &\rightarrow \quad \mathcal{E}^*\\
L_{\wt{E^*}} \quad &\mapsto \quad
\wt{E^*}=\Gamma(L_{\wt{E^*}}):= \menge{ \big( e^*,L_{\wt{E^*}}(e^*) \big) \in E^* \oplus E=(W \oplus U)^{\perp} }{e^* \in E^*}
\end{align*}
$\wt{E^*}$ being totally isotropic imposes the following condition on $L_{\wt{E^*}}=L$:
\begin{equation}\label{eq:Skew}
B ( e^*+Le^*,f^*+Lf^* ) =  0 \quad \Leftrightarrow \quad  B(e^*,Lf^*)+B(Le^*,f^*)=0 \quad \forall e^*,f^* \in {E^*} 
\end{equation}
The subset in $Hom({E^*},E)$ satisfying condition (\ref{eq:Skew}) will be denoted by \mb{$Skew({E^*},E,B)$}.
Any $L_{\wt{E^*}} \in Skew({E^*},E,B)$ yields the isometry $\wt{\rho}_{\wt{E^*}}(e+e^*)=\big(e+e^*+L_{\wt{E^*}}(e^*)\big)$ on $E \oplus E^*$.
The isometric extension to $V$ that acts identically on $(E \oplus {E^*})^{\perp}$ will be denoted by \mb{${\rho}_{\wt{E^*}} \in O(V,B)$}.
\end{description}

Summarizing, the following bijection has been obtained:
\begin{align*}
&Hom({W},E) \times Hom({U},E) \times Skew({E^*},E,B) \quad \longrightarrow \quad \mathcal{D}\\
&\quad \Big( L_{\wt{W}},L_{\wt{U}},L_{\wt{E^*}} \Big) \quad 
\longmapsto \quad \Big( E,\;\Gamma(L_{\wt{W}}),\;\Gamma(L_{\wt{U}}),\;{\rho}_{\wt{W},\wt{U}} \,\Gamma(L_{\wt{E^*}}) \Big)
\end{align*}

In other words  $\mathcal{D}$ is generated by the isometries ${\rho}_{\wt{U},\wt{W}} \circ {\rho}_{\wt{E^*}}$ acting on the
initial decomposition $D=(E,{W},{U},{E^*}) \in \mathcal{D}$. These isometries may be characterized as those in $\Gamma$
which act identically on the vector spaces $E$, $F/E$ and $F^{\perp}/E$. Consequently they form a subgroup of $\Gamma$:
%
\begin{align*}
\text{\mb{$\Gamma_{decomp.}$}} & := \big\{\gamma \in \Gamma \;\big| \; 
  \gamma|_E=id_E,\;\gamma|_{F/E}=id|_{F/E},\; \gamma|_{F^{\perp}/E}=id|_{F^{\perp}/E} \big\} \\
&= \big\{ {\rho}_{\wt{U},\wt{W}} \circ {\rho}_{\wt{E^*}} \;\big| \; 
   \wt{U} \in \mathcal{U},\;\wt{W} \in \mathcal{W},\;\wt{E^*} \in \mathcal{E}^* \big\}
\end{align*}
Moreover $\Gamma_{decomp.}$ is a normal subgroup of $\Gamma$ since for $\gamma_d \in \Gamma_{decomp.}$ and an
arbitrary $\gamma \in \Gamma$ it is easily verified that $\gamma \cdot \gamma_d \cdot \gamma^{-1}$ is the identity
on $E$, $F/E$ and $F^{\perp}/E$.

In order to determine whole $\Gamma$ (as set) it remains to investigate the possible choices 
for an adapted basis $\beta_{\wt{D}}$ w.r.t. $\wt{D} \in \cal{D}$. 
Fix $D=(E,U,W,E^*) \in \cal{D}$ and ${\beta}_{D}=\{ {\beta}_{{E}},{\beta}_{{W}},{\beta}_{{U}},{\beta}_{{E^*} }\} \in \cal{B}_{D}$.
For another $(V,B,{D})$-adapted basis
$\wt{\beta}_{{D}}=\{\wt{\beta}_{{E}},\wt{\beta}_{{W}},\wt{\beta}_{{U}},\wt{\beta}_{{E}^*} \}$ obviously  exists $P_E \in Gl({E})$,
$P_W \in O({W},B|_{W \times W})$ and  $P_U \in O({U},B|_{U \times U})$ such that
\[
\wt{\beta}_{{E}} = P_E({\beta}_{{E}}),\quad
\wt{\beta}_{{W}} = P_W({\beta}_{{W}}),\quad
\wt{\beta}_{{U}} = P_U({\beta}_{{U}}),\quad
\wt{\beta}_{{E^*}} = P_{E^*}({\beta}_{{E^*}})
\] 
where $P_{E^*} \in Gl(E^*)$ is the dual of $P_E$ i.e.
$B(P_{E^*}e^*,e)=B(e^*,{P_E}^{-1}e)$. 
\label{DualBasis} 
This suggests
\begin{align*}
\text{\mb{$\Gamma_{basis}^{{D}}$}} &:= \{ \gamma \in \Gamma \;\big| \; \gamma {D}={D}\} \\
&= \Big\{ P_E \oplus P_W \oplus P_U \oplus ({P_E}^{-1})^* \in
 Gl({E}) \times O({W}B|_{W \times W}) \times O({U},B|_{U \times U}) \times Gl({E^*}) \Big\}
\end{align*}

Now $\Gamma_{decomp.} \cap \Gamma_{basis}^{{D}} = \{id\}$ and
$\big\langle \Gamma_{decomp.}, \Gamma_{basis}^{{D}} \big\rangle = \Gamma $. Since
$\Gamma_{decomp.}$ is normal, $\Gamma$ is isomorphic to the
semi-direct product $\Gamma_{decomp.} \rtimes_C \Gamma_{basis}^{{D}} $,
where $C(\gamma)=l_{\gamma}\circ r_{{\gamma}^{-1}}$ is the conjugation by $\gamma$ (see \cite[III.3.14.]{HilgertN91}).
Thus we proved the following

\begin{prop}\label{prop:AdapBasis}
With the notation from above the following groups are isomorphic:
\[
\Gamma_{decomp.} \rtimes_C \Gamma_{basis}^{{D}} \simeq \Gamma
\]
Moreover, if ${D}=(E,{W},{U},{E^*})$ is a fixed adapted decomposition with respect to $F$ and $B$,
then there exists a bijection between the following sets:
\[
\begin{array}{cc}
Hom({W},E) \times Hom({U},E) \times Skew({E^*},E,B)\times \\
 \times Gl({E}) \times O({W},B|_{W \times W}) \times O({U},B|_{U \times U})
\end{array}
\quad \longleftrightarrow \quad \{\mathcal{B}_{\wt{D}}\}_{\wt{D} \in \mathcal{D}} \quad \underset{\ref{re:AdapBasisII}(i)}{ \longleftrightarrow } \quad \Gamma
\]
\end{prop}

The first assertion is essentially \cite[p.70,prop.4]{Boubel00} where much attention
is put on the group $\Gamma$ itself (e.g. one knows, that $\Gamma_{decomp.}$ is nilpotent and $\Gamma_{basis}^{D}$
is reductive). The determination of the concrete transformations of one adapted decomposition into another by means of
graphs of linear maps has been elaborated since we will need them in the following.

\subsection{Iterated adapted decompositions of solvable symmetric triples}\label{ss:IteratedAdapDecomp}

Given a solvable symmetric triple $\tau=(\g,\sigma,B)=(\h \oplus \m,B)$. Then $[\g,\g]=\big( [\h,\h]+[\m,\m]\big) \oplus [\h,\m]=
\h \oplus [\h,\m] \subsetneq \g$.
If in addition $\tau$ is indecomposable, $[\h,\m]$ is non-zero (unless $\m$ is one dimensional) 
and its orthogonal complement $[\h,\m]^{\perp} \cap \m=\z(\g)$ is totally
isotropic by lemma \ref{lem:Center}. In particular $[\h,\m] \subsetneq \m$ is a proper subspace. 

Let $V:=\m, \; F:=[\h,\m]$ then in the notation of the preceeding section $E = F^{\perp} \cap F = 
[\h,\m]^{\perp} \cap [\h,\m] = \z \cap \z^{\perp} = \z$. Thus we obtain 

\begin{lem}\label{lem:1.Decomp}
Any solvable indecomposable symmetric triple $\tau=(\h \oplus \m,B)$  admits an $([\h,\m],\m,B|_{\m \times \m})$-adapted
decomposition of the form 
\[
\m= \ub{14}{-2}{0}{= [\h,\m]} W  \;\overset{\perp}{\oplus}\;  \ob{12}{1}{0}{\sst sing. paired}  \big(\,\z  \; \oplus \; \z^* \big) 
 \quad \text{ i.e. } \quad \ad(\h)|_{\m}= 
\left( \begin{array}{c | c | c} \nosp{0}{4}{$\z$} 0  & \nosp{0}{4}{$W$} \,*\, &  \nosp{0}{4}{$\z^*$} * \\ 
\hline 0 & * & * \\ \hline 0 & 0 & 0 \\ 
\end{array} \right)
\]
For two such decompositions $\m = \z \oplus W_i \oplus \z^*_i,\,(i=1,2)$ there exists an $S \in Hom(W_1,W_2)$ 
which is an isometry w.r.t. $B_i=B|_{W_i \times W_i }$
and an equivalence w.r.t. the representations  $\rho_i= pr^{\z \oplus W_i}_2 \circ \ad(\cdot)|_{W_i}$.
\end{lem}

\begin{proof}
As mentioned above we apply the generalized Witt decomposition \ref{lem:WittDecomp} to the subspace $F=[\h,\m]$. 
Since $E=\z$ is totally isotropic the complement $U$ does not occur (cf. remark \ref{re:AdapBasis}(2)). 
Furthermore $\z$ is the center of $B|_{[\h,\m] \times [\h,\m]}$ and $\ad(\h)$ and thus they induce a non-degenerate bilinear form $\wh{B}$
and an $\h$-representation $\wh{\ad}$ on the quotient $[\h,\m]/ \z$. On the other hand 
the complement $W_i$ is equipped with the restrictions $B_i$ and $\rho_i$.
Then the canonical identification $W_i \simeq [\h,\m]/\z$ yields the isometry and $\h$-equivalence 
$\big(W_i,B_i,\rho_i \big) \simeq \big([\h,\m]/\z,\wh{B},\wh{\ad} \big)$. In this sense $(W_i,B_i,\rho_i)$ is independent
of the chosen complement $W_i$. 
\end{proof}

\begin{samepage}
This procedure can be iterated:
\begin{itemize}
\item An initial decomposition is provided by lemma \ref{lem:1.Decomp}: Set $E_0:=\z$, $E^*_0:=\z^*$, $W_0:=W$, $U_0:=\{0\}$.
\item $\m= E_0 \oplus \ldots \oplus E_i \; \oplus \;W_i \;\oplus\; U_0\oplus \ldots \oplus U_i \;\oplus \;E^*_i \oplus \ldots \oplus E^*_0$
is inductively defined by 
\bi
\item $F_{i+1} := pr_{W_i}[\h,W_i] = [\h,W_i]_{W_i} \subset W_i$.
\item $E_{i+1}$, $W_{i+1}$, $U_{i+1}$ and $E_{i+1}^*$ are obtained from an $(F_{i+1},W_i,B|_{W_i \times W_i})$-adapted decomposition:
\[
W_i=\underbrace{\overbrace{ \big(F_{i+1} \cap F_{i+1}^{\perp} \big)}^{=E_{i+1}} \oplus W_{i+1}}_{=F_{i+1}} \oplus U_{i+1} \oplus E_{i+1}^*
\]
\ei
\item $\rho_i(\cdot):= pr_{W_i}\circ \ad(\cdot)|_{W_i} : \h \rightarrow \gl(W_i)$.
\end{itemize} 
\end{samepage}

\begin{prop}\label{prop:AdapDecomp} Consider a decomposition gained by induction as explained above:
\[
\m=\underbrace{\z \oplus E_1 \oplus \ldots \oplus E_n}_{\ds =:E} \oplus W_n \oplus
\underbrace{U_1  \oplus \ldots \oplus U_n}_{\ds =:U} \oplus
\underbrace{E_n^*  \oplus \ldots \oplus E_1^* \oplus \z^*}_{\ds =:E^*}
\]
\begin{enumerate}
\item[(i)]For $0 \le i \le n$ the kernel of the representation $\rho_i$ is given by
\[
\ker(\rho_i)= \menge{ X \in W_i}{ [\h,X]_{W_i}=0 } = E_{i+1} \oplus U_{i+1}
\]
and $\rho_i$ is a nilpotent representation of $\h$.
\item[(ii)] $F_{n+1}=[\h,W_n]_{W_n}$ is non-singular iff $F_{n+1}=\{0\}$ (i.e. $\rho_{n}=0$).
In particular this is the case, if $B|_{W_n^2}$ is definite. 
Thus every inductive decomposition ends with $W_n=U_{n+1}$ (we will keep the notation $W_n$) and 
is said to be \mb{$complete$} (in the sense, that it cannot be decomposed further).

Note, that $E_i \ne \{0\}$ for $0 \le i \le n$, but in contrast some of the subspaces $U_i$ or $W_n$ might be zero.
\item[(iii)] The integers $\dim(E_i)$, $\dim(W_i)$ and $\dim(U_i)$ $(0 \le i \le n)$ are invariants of the symmetric triple.
\item[(iv)] There exists a basis
$\beta=\{\beta_{E_0}, \ldots \beta_{E_n},\beta_{W_n},\beta_{U_1}, \ldots \beta_{U_n},\beta_{E_n}^*, \ldots \beta_{E_0}^*\}$
 adapted to the above decomposition of $\m$ (i.e. $span\{\beta_{E_i}\}=E_i, \ldots$), such that
\[
Mat_{\beta}(B)=
\begin{array}{c}
\begin{array}{r r r r} E & W_n & U & E^*\end{array} \\
\left( \begin{array}{r | rr |r} &&& {I_1}^t \\ \hline  & I_2 &&\\ && I_3 &\\ \hline I_1 &&&
\end{array} \right)
\end{array}
\]
where 
\[
I_1=
\left(
\begin{smallmatrix}
0&& I_{\dim(E_n)}\\
&\text{\rotatebox[origin=c]{80}{$\ddots$}}&\\
I_{\dim(E_0)} &&0
\end{smallmatrix}
\right),\qquad
I_2=I_{\Sig(W_n)},\qquad
I_3=
\left(
\begin{smallmatrix}
I_{\Sig(U_1)}&&0\\
&\ddots &\\
0&&I_{\Sig(U_n)}
\end{smallmatrix}
\right)
\]
\item[(v)] With respect to every complete decomposition and an adapted basis as in (iv) the endomorphisms 
$\ad(\h) \subset \gl(\m)$ are upper triangular matrices of the following form:
{ \renewcommand{\arraystretch}{1.5}
\[
Mat_{\beta}(\ad(h)|_{\m})=
\left( \begin{array}{c | c} \begin{array}{ccc} M_E^h & M_W^h & M_U^h \end{array}  & M_{E^*}^h  \\ \hline
\text{\Large $0$} & \begin{array}{r} -I_2 \,(M_W^h)^t \,{I_1}^t \\
  -I_3 \, (M_U^h)^t \, {I_1}^t \\ - I_1 \, (M_E^h)^t \, {I_1}^t
\end{array}
\end{array} \right)
\]
}
where 
\[
I_1 M_{E^*}^h = -(I_1 M_{E^*}^h)^t, \qquad
M_E^h=
\left( \rca{1.2}{\ba{cccc} \nosp{0}{5}{$\z$}0 & \nosp{0}{5}{$E_1$}\;* & \nosp{0}{5}{$\cdots$}\cdots &\nosp{0}{5}{$E_n$}\;*  \\
0 & 0 & \ddots & \vdots \\ 
\vdots & \vdots & \ddots & * \\
0 & 0 & \cdots & 0 \ea } \right),\qquad
M_U^h=
\left(
\rca{1.2}{\ba{ccc}
\ob{15}{0}{0}{U_1 \oplus \cdots \oplus U_n}*      & \hdots & *      \\
\vdots &        & \vdots \\
0      & \cdots & *      \\
\vdots & \ddots & \vdots \\
0      & \cdots & 0      \\
\ea}
\right)
\]
\item[(vi)] $\quad \h = \big[ E_0^*,W_0 \oplus E_0^* \big] + \ldots + \big[ E_n^*,W_n \oplus E_n^* \big]$

\vspace{.3cm}
In particular $E \oplus W_n \oplus U$ is an abelian subalgebra. Recall that by construction

$W_i = E_{i+1} \oplus W_{i+1} \oplus U_{i+1} \oplus E_{i+1}^* =\bigoplus_{k = i+1}^n (E_k \oplus U_k \oplus E_k^*) \oplus W_n$.
\end{enumerate}
\end{prop}

\begin{proof}(i) First we determine the kernel of $\rho_i$:
\begin{align*}
[\h,X]_{W_i}=0 \quad &\Leftrightarrow \quad B \big( [\h,X],W_i \big)=\big( X,[\h,W_i] \big)=0, \quad \text{($W_i$ non-singular)}\\
&\Leftrightarrow \quad X \in [\h,W_i]^{\perp} \cap W_i = F_{i+1}^{\perp}=E_{i+1} \oplus U_{i+1}
\end{align*}
According to proposition \ref{prop:Levi} $\ad(h)|_{\m}$ $(h \in \h)$ is nilpotent.
$\rho_0(h)=pr_{W_0} \circ \ad(h)|_{W_0}$ inherits this property due to $(id_{\m}-pr_{W_0}) \ad(\h)W_0 \subset \z$. 
Similarly $(id_{W_{i-1}} - pr_{W_i}) \rho_{i-1}(\h)W_i \subset E_i = \ker(\rho_{i-1})$ and
\begin{equation}\label{eq:Rho}
\rho_i(\cdot) =pr_{W_i} \circ \rho_{i-1}(\cdot)|_{W_i}
\end{equation}
shows by induction the nilpotency of all representations $\rho_i$.

(ii) Assume $F_{n+1}=[\h,W_n]_{W_n}$ is non-singular w.r.t. $B$. Then
$W_n = F_{n+1} \oplus (F_{n+1}^{\perp} \cap W_n) = W_{n+1} \oplus U_{n+1}$
since $F_{n+1} \cap F_{n+1}^{\perp} = E_{n+1} = \{0\}$ and we get
\[
W_{n+1}=F_{n+1}=[\h,W_n]_{W_n}=[\h,W_{n+1} \oplus U_{n+1}]_{W_n} \underset{(i)}{=}[\h,W_{n+1}]_{W_n} =
\rho_{n}(\h)\, W_{n+1} \underset{(\ref{eq:Rho})}{=} \rho_{n+1}(\h)\, W_{n+1}
\]
$\rho_{n+1}$ being nilpotent implies $W_{n+1} =\{0\}$ due to the lemma of Engel. 

(iii) Given two complete adapted decompositions of $(\h \oplus \m,B)$, we consider the induced flag of the subspaces
$\m \supset W_0 \supset \ldots \supset W_{n}$ resp. $\m \supset \wt{W}_0 \supset \ldots \supset \wt{W}_{m}$.
Each of these spaces is equipped with the $\h$-representation $\rho_i$ (resp. $\wt{\rho}_i$) and the non-degenerate bilinear form
$B_i=B|_{W_i \times W_i}$ (resp. $\wt{B}_i=B|_{\wt{W}_i \times \wt{W}_i}$). We will show by
induction, that there is a family of linear maps $S_i : \big( W_i,\rho_i,B_i \big) \rightarrow  \big( \wt{W}_i,\wt{\rho}_i,\wt{B}_i \big)$
such that $S_i$ is an isometry and an $\h$-equivalence. The initial case $i=0$ is covered by lemma \ref{lem:1.Decomp}.
Now assume there exists a map $S_i,\;(i\geq 0)$ with the required properties. Since $E_{i+1}$ is in the kernel of $B|_{F_{i+1} \times F_{i+1}}$ and
$\rho_i$ they factor through $E_{i+1}$. Analogous to \ref{lem:1.Decomp} this yields an isometry and \mbox{$\h$-equi}\-valence
$\big( W_{i+1},\rho_{i+1},B_{i+1} \big) \simeq \big( F_{i+1}/E_{i+1} , \wh{\rho_i},\wh{B_i}\big)$. The same holds for the other
decomposition and finally $S_i$ induces an isometry and $\h$-equivalence on the quotients $F_{i+1}/E_{i+1} \simeq \wt{F}_{i+1}/\wt{E}_{i+1}$. 
Obviously the obtained isometries and $\h$-equivalences can be composed to a suitable map $S_{i+1}$.
Moreover the dimensions of the subspaces $U_i$ and $E^*$ are determined by
the dimensions of $W_i$ and $E_i$.  Since the latter are invariant under the
choice of an adapted decomposition, all of them are.

(iv): We fix a complete inductive decomposition and denote in each step the chosen $(\rho_i(W_i),W_i,B_i)$-adapted
decomposition of $W_i$ by $D_i$, $(0 \le i \le n)$ 
(The first decomposition $\m=\z \oplus W_0 \oplus \z^*$ will be denoted simply by $D$).
Then the basis $\beta$ is defined successively starting with an 
$(\m,B,D)$-adapted basis 
$\{\beta_{E_0},\beta_{W_0},\beta_{E_0}^*\}$ (see \ref{dfn:AdapBasis}). For $0 \le i < n$ we replace $\beta_{W_i}$
by an $(W_i,B_i,D_i)$-adapted basis $\{\beta_{E_{i+1}},\beta_{W_{i+1}},\beta_{U_{i+1}},\beta_{E_{i+1}}^*\}$ of $W_i$.
Finally we choose a pseudo-orthonormal basis $\beta_{W_n}$ of $W_n$ and sort the basis as cited. 

(v): The uniform triangulation of the holonomy representation by a basis adapted to the metric
was the purpose the inductive decomposition has been developed for. The essential observation is the following:
\begin{equation}\label{eq:Kernel}
V \subset \ker(\rho_k) \subset W_k, \quad \Rightarrow \quad [\h,V] \subset \bigoplus_{i=0}^k E_i
\end{equation}
$V \subset \ker(\rho_k)$ implies $\rho_k(\h)V \underset{(\ref{eq:Rho})}{=} pr^{E_{k} \oplus W_{k}}_2 \circ \rho_{k-1}(\h)V =\{0\}$,
hence $\rho_{k-1}(\h)V \subset E_k$. Iterating this argument for $\rho_{k_2}, \ldots \rho_0$ finally yields (\ref{eq:Kernel}).
In particular (\ref{eq:Kernel}) can be applied to $\ker(\rho_{k}) \underset{(i)}{=} E_{k+1} \oplus U_{k+1}$:
\begin{equation}\label{eq:[h,E+U]}
[\h,E_{k+1} \oplus U_{k+1}]=\ad(\h)(E_{k+1} \oplus U_{k+1}) \subset \bigoplus_{i=0}^k E_i
\end{equation}
Thus the matrices $M^h_E,\;h \in \h$ are strictly upper triangular bloc-matrices and the
lower left triangle of the matrices $M^h_U,\;h \in \h$ is zero.  
Further $W_n = \ker(\rho_n)$ implies
\[
[\h,W_n]=\ad(\h)(W_n) \subset \bigoplus_{i=0}^n E_i=E
\]

Thus the shape of the columns of $Mat_{\beta}(\ad(h))$ corresponding to the subspaces $E$, $W$ and $U$ is determined.
The remaining part of $Mat_{\beta}(\ad(h))$ is given by $\ad(h) \in \so(\m,B)$. In terms of matrices:
\[
Mat_{\beta}\big( \ad(h) \big)^t \cdot Mat_{\beta}(B) = - Mat_{\beta} \big( \ad(h) \big) \cdot Mat_{\beta} ( B)
\]
Using the identities
${I_1}^{-1} = {I_1}^t,\; {I_j}^{-1} = I_j$ ($j=2,3$) it is easy to verify the cited form of $Mat_{\beta}(\ad(h))$.
Since  $- I_1 \, (M_E^h)^t \, {I_1}^t$ is a strictly upper triangular bloc, the whole matrix is.

(vi) $B$ being non-degenerate and $\ad(\g)$-invariant allows the following conclusions: 
\begin{align*}
&B\Big( [ E^{\perp} , E^{\perp} ] \,,\,\h \Big)=B\big( E^{\perp} ,\underbrace{[E \oplus W_n \oplus U ,\h]}_{\subset E} \big)=\{0\}\\
&\quad \Rightarrow \quad  \big[ E^{\perp},E^{\perp} \big]=\{0\} \; \text{ and } \; \h =[E^*,\m]
\end{align*}
More precisely:
\begin{align*}
&B\Big( [ E_i^* , E_k \oplus U_k ] \,,\,\h \Big) = B \big( E_i^* , [E_k \oplus U_k ,\h] \big) \underset{ (\ref{eq:[h,E+U]}) }{\subset} 
B \Big( E_i^* , \bigoplus_{j=0}^{k-1} E_j \Big) \\
&\quad \Rightarrow \quad  \big[ E_i^*,E_k \oplus U_k \big]=\{0\} \; \text{ for } k \le i
\end{align*}
\end{proof}

The depth of the iteration depends on the signature of the symmetric triple:
Setting $W_{-1}:=\m$ and for brevity $\Sig(V)=\Sig(B_{V \times V})$ $(V \subset \m)$,
we get the following estimates $(0 \le i \le n)$:
\[
\Sig(W_i) = \Sig(W_{i-1})-\Sig(U_i)-\Sig(E_i \oplus E_i^*) \le \Sig(W_{i-1}) - (1,1)
\]
This implies for the \emph{number of iterations in a complete adapted decomposition} $n \le \min(\Sig(\m))$.
The first two cases i.e. Lorentzian signature and signature $(2,n-2)$ will be discussed and classified in the subsequent sections.

\subsection{Solvable symmetric triples with maximal center}\label{subsec:MaxCenter}

An indecomposable symmetric triple $\tau = (\g,\sigma,B)$ of signature $(p,q)$ $(p \le q)$ is said to be of 
\mb{maximal center} if $\dim(\z(\g))=p$, i.e. $\z$ is a maximal isotropic subspace. If we ignore the Riemannian case (i.e. $p=0$)
the center being maximal implies the solvability of $\g$ due to corollary \ref{cor:SeperationProperty}(iv).

For a symmetric triple with maximal center the iterated adapted decomposition reduces to an ordinary Witt-decomposition. 
Indeed, consider the first decomposition $\m = \z \oplus W_0 \oplus \z^*$
as in lemma \ref{lem:1.Decomp}. Then $B|_{W_0 \times W_0}$ is positive definite as remarked in \ref{re:AdapBasis}(2). 
Hence $W_1 \subset [\h,W_0]_{W_0}=\{0\}$ by proposition \ref{prop:AdapDecomp}(ii). Lets simplify the
notation by setting $W=W_0$.  

\begin{re}
The case of split signature $(p,p)$ i.e. $\m = \z \oplus \z^*$ can be treated simultaneously 
minding that some statements and conditions might be empty, since $W=\{0\}$. 
\end{re}

Now fix an adapted decomposition of $\m$ together with an adapted basis:
\[
\beta_{\z} =\{Z_1, \ldots Z_p\},\quad \beta_{\z}^* =\{Z_1^*, \ldots Z_p^*\}, \quad 
\beta_{W} = \{W_1, \ldots W_{q-p}\}
\]
In the following we will use the convention to denote indeces out of $\{1,\ldots,p\}$ by Latin letters $i,j,k,\ldots$ and 
indeces out of $ \{1,\ldots,q-p\}$ by Greek letters $\alpha,\beta,\gamma,\ldots$.
According to \ref{prop:AdapDecomp}(vi) $\h$ is generated by
\[
\rca{1.3}{\begin{array}{ll} Y_{ij}&:=[Z_i^*,Z_j^*] \\ X_{i\alpha}&:=[Z_i^*,W_{\alpha}] \end{array}\quad
\begin{array}{ll} Y =span\{Y_{ij}\}_{i<j} \\ X =span\{X_{i\alpha}\} \end{array} }\quad 
 \h = Y + X
\]
The structure coefficients of the Lie algebra $\g$ concerning $[\m,\m]$ are thus encoded in the linear dependency of the 
$X_{i\alpha},Y_{ij}$, those of $[\h,\h]$ will follow by the Jacobi identity. For the coefficients
of $[\h,\m]$ we use our knowledge about the action of $\h$ on $\m$ due to \ref{prop:AdapDecomp}(v):

\begin{equation}\label{eq:hAction}
Mat_{\beta}(\ad(h)|_{\m})=
\left( \ba{ccc} \nosp{1}{5}{$\z$}0 &  \nosp{2}{5}{$W$} M_W^h & \nosp{2}{5}{$\z^*$} M_{\z^*}^h \\
0 & 0 & -{M_W^h}^t \\ 0&0&0 \ea \right),\qquad
 M_{\z^*}^h=-(M_{\z^*}^h)^t
\end{equation}
Set:
\begin{align*}
\left(M_W^{Y_{ij}}\right)_{k\alpha} &= -a_{ijk\alpha} & \left(M_{\z^*}^{Y_{ij}}\right)_{kl} &= b_{ijkl} \\
\left(M_W^{X_{i\alpha}}\right)_{j\beta} &= -f_{ij\alpha\beta} & \left(M_{\z^*}^{X_{i\alpha}}\right)_{kl} &= c_{kli\alpha} 
\end{align*}
(This notation has been taken from \cite{CahenP70} to enable a direct comparison). 
Our plan is to formulate the conditions for $(\h \oplus \m,B)$ being an indecomposable symmetric triple in terms of the 
structure coefficients $a,b,c,f$. Lets recall what has to be checked:

\begin{samepage}
\begin{enumerate}\label{Checklist}
\item[(A)]: $(\h \oplus \m,[\cdot,\cdot])$ being a Lie algebra (skew-symmetry and Jacoby-identity)
\item[(B)]: $\ad(\h) \subset \so(\m,B|_{\m \times \m})$
\item[(C)]: $\ad(\cdot)|_{\m}:\h \rightarrow \gl(\m)$ being faithful
($\Leftrightarrow$ $B$ being non-degenerate on $\g$ due to
\ref{lem:InvBi}(iii)) 
\item[(D)]: indecomposability of the symmetric triple
\end{enumerate}
\end{samepage}

Note that (B) is ensured by (\ref{eq:hAction}). Before discussing the implications of (A) systematically, it is
advantageous to consider the basic symmetries derived from $\RR=B([[\cdot,\cdot],\cdot],\cdot) \in \mathcal{T}^4(\m)$ 
being the curvature tensor of the symmetric space associated to the symmetric triple.
Algebraically this is a consequence of (A) and (B) (see the proof of
\ref{lem:InvBi}(ii)). Moreover  $a,b,c,f$ are exactly the
nontrivial entries in the curvature tensor (up to its symmetries):
\begin{equation}\label{eq:StrucCoef} 
\rca{1.3}{\begin{array}{rclcl} 
a_{ijk\alpha} &=& -B([Y_{ij},W_{\alpha}],Z_k^*) &=&-\RR(Z_i^*,Z_j^*,W_{\alpha},Z_k^*)\\ 
b_{ijkl}      &=&-B([Y_{ij},Z_k^*],Z_l^*)      &=&  \phantom{-} \RR(Z_i^*,Z_j^*,Z_k^*,Z_l^*)\\
f_{ij\alpha\beta} &=& -B([X_{i\alpha},W_{\beta}],Z_j^*) &=&-\RR(Z_i^*,W_{\alpha},W_{\beta},Z_j^*)\\ 
c_{kli\alpha} &=& \phantom{-}B([X_{i\alpha},Z_k^*],Z_l^*) &=&  \phantom{-}\RR(Z_i^*,W_{\alpha},Z_k^*,Z_l^*)\\ 
\end{array}} 
\end{equation}  
\begin{itemize} \item
{\bf pair symmetry :} \begin{align*} a_{ijk\alpha} =c_{ijk\alpha},\qquad
b_{ijkl} =b_{klij},\qquad f_{ij\alpha\beta} =f_{ji\beta\alpha} \end{align*}
Thus $c$ can be replaced by $a$. \item {\bf skew symmetry in the 1. and 2.
argument:} \begin{align*} a_{ijk\alpha}=-a_{jik\alpha},\qquad b_{ijkl}
=-b_{jikl} \end{align*}
\item {\bf skew symmetry in the 3. and 4. argument:} This is guaranteed by (\ref{eq:hAction}). In particular:
\[
b_{ijkl} =-b_{ijlk}
\]
\item {\bf 1. Bianchi identity :}
\[
\sigma_{ijk}\big(a_{ijk\alpha}\big)=0,\qquad \sigma_{ijk}\big(b_{ijkl}\big) =0, \qquad
f_{ij\alpha\beta} \underset{W \text{ abelian}}{=}f_{ji\alpha\beta} \underset{pair sym.}{=}f_{ij\beta\alpha}
\]
\end{itemize}
Next we will uncover the structure of $\h$:

\begin{lem}
$X \subset \h$ is central in $\h$ i.e. $[X,\h]=\{0\}$ and $[Y,Y] \subset X$.
Thus the holonomy algebra $\h$ is either abelian or two step nilpotent. 

\end{lem}
\begin{proof} 
\begin{align*}
&[X,\h]=\big[[\z^*,W],\h\big]\underset{Jac.-Id.}{=}
\big[\z^*,\underbrace{[W,\h]}_{\subset \z}\big]+\big[W,\underbrace{[\z^*,\h]}_{\subset W \oplus \z}\big]=\{0\}\\
&[Y,\h]\underset{Jac.-Id.}{=}\big[\z^*,\underbrace{[\z^*,\h]}_{\subset \z \oplus W}\big] \subset X
\end{align*} 
\end{proof}  
Concrete in terms of the structure coefficients we get:
\begin{multline}\label{eq:LAh}
[Y_{ij},Y_{kl}] =\big[Y_{ij},[Z_k^*,Z_l^*]\big]\underset{Jac.-id.}{=}-\big[Z_l^*,[Y_{ij},Z_k^*]\big]+\big[Z_k^*,[Y_{ij},Z_l^*]\big] = \\
\underset{(\ref{eq:hAction})}{=}-\big[Z_l^*,\sum_{\alpha}a_{ijk\alpha}W_{\alpha}\big]+\big[Z_k^*,\sum_{\alpha}a_{ijl\alpha}W_{\alpha}\big]
 = -\sum_{\alpha}a_{ijk\alpha}X_{l\alpha}+\sum_{\alpha}a_{ijl\alpha}X_{k\alpha}
\end{multline}
By the $\ad(\g)$-invariance of $B \in S^2(\g)$ the curvature can be considered as non-degenerate bilinear form on $\h$ by 
$B([X_1,Y_1],[X_2,Y_2])=\RR(X_1,Y_1,X_2,Y_2)$. Thus $a,b,f$ may be interpreted as ``metric'' coefficients which will be helpful
for the discussion of (C):
\begin{equation}\label{eq:hMetric}
a_{ijk\alpha}=B(Y_{ij},X_{k\alpha}), \qquad
b_{ijkl}=B(Y_{ij},Y_{kl}), \qquad
f_{ij\alpha\beta}=B(X_{i\alpha},X_{j\beta})
\end{equation} 
Now we return to (A) in order to determine the remaining Jacoby-identities:
\begin{itemize}
\item $\sigma_{1,2,3}\big[[m_1,m_2],m_3\big]=0,\;m_i \in \m$ is fulfilled due to the 1. Bianchi-identity.
\item $\sigma_{1,2,3}\big[[h_1,h_2],h_3\big]=0,\;h_i \in \h$ is clear, since $\h$ is abelian or to 2 step nilpotent.
\item $\big[[m_1,m_2],h\big]=\big[m_1,[m_2,h]\big]-\big[m_1,[m_2,h]\big]$ yields the Lie algebra structure of $\h$ (cf. (\ref{eq:LAh})).
\item $\big[[h_1,h_2],m\big]+\big[[h_2,m],h_1\big]+\big[[m,h_1],h_2\big]=0$
leads to new relations:\\ 
\begin{enumerate}
\item[(1)] If $h_1 \in X \subset \z(\h)$ the Jacobi-identity asserts $\big[\ad(h_1)|_{\m},\ad(h_2)|_{\m}\big]=0$. Hence by (\ref{eq:hAction}):
\begin{align*}
&M_W^{h_1} \circ (M_W^{h_2})^t = M_W^{h_2} \circ (M_W^{h_1})^t = \big( M_W^{h_1} \circ (M_W^{h_2})^t \big)^t\\
\Leftrightarrow \quad & 
\begin{cases}
\quad \sum_\gamma f_{ik\alpha\gamma}\cdot f_{jl\beta\gamma}=\sum_\gamma f_{il\alpha\gamma}\cdot f_{jk\beta\gamma}&
\text{ for } h_1=X_{i\alpha},h_2=X_{j\beta}\\
\quad \sum_\gamma f_{im\alpha\gamma}\cdot a_{kln\gamma}=\sum_\gamma f_{in\alpha\gamma}\cdot a_{klm\gamma}&
\text{ for }h_1=X_{i\alpha},h_2=Y_{kl}
\end{cases}
\end{align*}
\item[(2)] If $h_1,h_2 \in Y$ the Jacobi-identity implies $\ad \Big([h_1,h_2] \Big)|_{\m}= \Big[\ad(h_1)|_{\m},\ad(h_2)|_{\m} \Big]$.
Setting $h_1=Y_{ij},h_2=Y_{kl}$ the evaluation on $W$ yields the same relation between the coefficients $a$ and $f$ 
as obtained in (1) whereas on $\z^*$ we get:
\end{enumerate} 
\begin{align*}
\Big( \big[\ad(Y_{ij}),\ad(Y_{kl}) \big]Z_n^* \Big)_{Z_m} 
&=\big(M_W^{Y_{ij}}\circ (M_W^{Y_{kl}})^t -M_W^{Y_{kl}}\circ (M_W^{Y_{ij}})^t \big)_{mn}\\
&= \sum_{\gamma} \big( a_{ijm\gamma}a_{kln\gamma}- a_{klm\gamma}a_{ijn\gamma}\big)\\
\Big( \ad \big( [Y_{ij},Y_{kl}] \big) Z_n^* \Big)_{Z_m}
&=\sum_{\gamma}\big[ a_{ijk\gamma} X_{l\gamma}- a_{ijl\gamma} X_{k\gamma}, Z^*_n \big]_{Z_m}\\ 
&\underset{(\ref{eq:hAction})}{=}\sum_{\gamma}\big(
a_{ijk\gamma}c_{mnl\gamma} - a_{ijl\gamma}c_{mnk\gamma} \big)=\sum_{\gamma}\big( a_{ijk\gamma}a_{mnl\gamma} - a_{ijl\gamma}a_{mnk\gamma}\big) 
\end{align*} 
\end{itemize}
Thus we proved the following

\begin{samepage}
\begin{prop}\label{prop:1.NFMaxCenter} \cite[page 343]{CahenP70}
Let $(\h \oplus \m,B)$ be an indecomposable symmetric triple with maximal center $\z \ne \{0\}$. If $\m=\z \oplus  W \oplus \z^*$
is an $([\h,\m],\m,B)$-adapted decomposition $D$ and $\beta=(\{Z_i\},\{W_{\alpha}\},\{Z_i^*\})$ is an $(\m,B,D)$-adapted basis,
then the only brackets possibly non-zero are\footnote{The following summation convention is used: 
if in one term an index occurs twice it is a summation index.}:
\begin{align}\label{eq:StructureCoefficients}
[Z_i^*,Z_j^*] &= \phantom{-}Y_{ij} \notag\\
[Z_i^*,W_{\alpha}] &= \phantom{-}X_{i\alpha} \notag\\
[Y_{ij},Z_k^*] &= \phantom{-} a_{ijk\gamma}W_{\gamma} +  b_{ijkl}Z_l \notag\\
[Y_{ij},W_{\alpha}] &= -  a_{ijl\alpha} Z_l\\
[X_{i\alpha},Z_k^*] &= \phantom{-} f_{ik\alpha\gamma}W_{\gamma} +  a_{kli\alpha}Z_l \notag\\
[X_{i\alpha},W_{\beta}] &= - f_{il\alpha\beta} Z_l \notag\\
[Y_{ij},Y_{kl}] &= \phantom{-}a_{ijl\gamma}X_{k\gamma} - a_{ijk\gamma}X_{l\gamma} \notag  
\end{align}
where $\h=span\{Y_{ij},X_{k\alpha}\}$. The structure coefficients obey the following symmetries and relations:
\begin{align*}
a_{ijk\alpha}&=-a_{jik\alpha}, & \sigma_{ijk}(a_{ijk\alpha})&=0 \\
b_{ijkl}&=b_{klij}=-b_{jikl}, &  \sigma_{ijk}(b_{ijkl})&=0 \\
f_{ij\alpha\beta}&= f_{ji\alpha\beta}=f_{ij\beta\alpha}
\end{align*}
\begin{gather}
 f_{ik\alpha\gamma}\cdot f_{jl\beta\gamma}  =  f_{il\alpha\gamma}\cdot f_{jk\beta\gamma} \label{eq:Kommutativ}\\
 f_{im\alpha\gamma}\cdot a_{kln\gamma}  =  f_{in\alpha\gamma}\cdot a_{klm\gamma} \label{eq:af}\\
 a_{ijm\gamma} \cdot a_{kln\gamma}- a_{klm\gamma} \cdot a_{ijn\gamma} =
 a_{ijk\gamma} \cdot a_{mnl\gamma} - a_{ijl\gamma} \cdot a_{mnk\gamma} \label{eq:a}
\end{gather}
\end{prop}
\end{samepage}

\begin{cor}\label{cor:(p,p)}
Let $(\h \oplus \m,B)$ be an indecomposable symmetric triple of signature $(p,p)$ with maximal center. 
If $\m=\z \oplus \z^*$ is a Witt-decomposition and $\beta=(\{Z_i\},\{Z_i^*\})$ a Witt-basis then
the only brackets possibly non-zero are
\begin{align*}
[Z_i^*,Z_j^*] &= Y_{ij}\\
[Y_{ij},Z_k^*] &=   b_{ijkl}Z_l
\end{align*}
where $\h =span\{Y_{ij}\}$. The structure coefficients $b$ satisfy 
\[
b_{ijkl}=b_{klij}=-b_{jikl} \qquad \sigma_{ijk}(b_{ijkl})=0
\]
\end{cor}

\begin{re}\label{re:LAConstruction}
If conversely parameters $a,b,f$ satisfying the relations of \ref{prop:1.NFMaxCenter} are predicted then
(\ref{eq:StructureCoefficients}) defines a Lie algebra structure on the vector space $\wt{\h} \oplus \m$, where 
$\wt{\h}=\Big(\bigoplus_{i,\alpha} \R \cdot X_{i\alpha} \Big) \oplus \Big(\bigoplus_{l<k} \R \cdot Y_{lk} \Big)$ and 
$\{Z_i,W_{\alpha},Z^*_i\}$ is a basis of $\m$ as before.
Moreover $\m$ carries an $\wt{\h}$-invariant bilinear-form $B$ by construction. Thus the conditions
(A) and (B) (see page \pageref{Checklist}) are satisfied. In order to achieve the validity of (C)
we set $\h=\wt{\h}/ker(\rho)$ where $\rho(\cdot)=\ad(\cdot)|_{\m}:\wt{\h} \rightarrow \gl(\m)$ 
(note that $\ker(\rho)=\z(\wt{\h}\oplus \m) \cap \wt{\h}$ is an ideal). Thus the adjoint representation of $\h$ on $\m$
becomes faithful which yields a unique $\ad(\g)$-invariant and non-degenerate extension of $B$ to whole $\g=\h \oplus \m$
s.t. $\h \perp_B \m$ (\ref{lem:InvBi}(iii)). 
In other words $(\h \oplus \m,B)$ is a symmetric triple and by the preceeding proposition every indecomposable symmetric triple
with maximal center arises from such a construction. \hfill \qedsymbol
\end{re}


\subsubsection[The action of $\Gamma=stab({[} \h,\m {]}) {\cap} O(\m,B)$ on the structure coefficients]
{\mb{The action of $\Gamma=stab([\h,\m]) \cap O(\m,B)$ on the structure coefficients}}

The next aim is to investigate how the coefficients $a,b,f$ of a solvable triple behave under the action of $\Gamma=stab([\h,\m]) \cap O(\m,B)$.
On the one hand there will appear new invariants which can be used to separate several classes of symmetric triples. On the other hand 
those parameters which are not invariant are used to fix specific adapted decompositions i.e. a normal form.
For both we need the results of section \ref{ss:AdapDecomp}.

We denote the $([\h,\m],\m,B)$-adapted decomposition by $\mathcal{D}$ and the $(\m,B,D)$-adapted basis by $\mathcal{B}_{D}$.
If we fix $D=(\z,W,\z^*) \in \mathcal{D}$ and $\beta=\{\beta_{\z},\beta_W,\beta_{\z}^*\}\in \mathcal{B}_{D}$, then by
proposition \ref{prop:AdapBasis} we have the following one-to-one correspondences:
\begin{align}\label{eq:one-to-one}
 Hom(W,\z) \times Skew(\z^*,\z,B)\quad &\leftrightarrow \quad \mathcal{D} \quad \simeq \quad \Gamma_{decomp.}\\
Gl(\z) \times O(W,B|_{W \times W}) \quad &\leftrightarrow \quad \mathcal{B}_D \quad \simeq \quad \Gamma_{basis}^D \notag
\end{align}

\begin{description}
\item[\mb{$\boxed{Hom(W,\z) \leftrightarrow \rho_{\wt{W}}}$}]\label{eq:HomWZ} 
To $L=L_{\wt{W}} \in Hom(W,\z)$ we associate $\rho_{\wt{W}} \in \Gamma_{decomp.}$ (cf. page \pageref{dfn:Rho}): 
\[
\wt{W}_{\alpha}:=\rho_{\wt{W}}(W_{\alpha})=(id_W + L)(W_{\alpha})=W_{\alpha} + \sum_k (L)_{k \alpha}Z_k
\]
In order to determine $\wt{Z^*}:=\rho_{\wt{W}}(Z^*_i)$ set $\wt{Z}_i^* = Z_i^* + \sum_{\gamma} g_{\gamma i}W_{\gamma} + \sum_k h_{ki}Z_k$.
Then $\wt{\z^*} \perp \wt{W}$ implies $g_{\alpha i} = - (L)_{i \alpha }$ and $\wt{\z^*} \perp \wt{\z^*}$ requires 
$h_{ij}+h_{ji}=-\sum_{\gamma}(L)_{\gamma i}(L)_{\gamma j}$. Of course $h_{ij}$ is not unique for $p > 1$ due to
the part $Skew(\z^*,\z,B)$ discussed next. We set  $h_{ij}=\frac{1}{2}\sum_{\gamma}(L)_{\gamma i}(L)_{\gamma j}$ and obtain
\begin{equation}\label{eq:MatRhoW}
Mat_{\beta}(\rho_{\wt{W}})=
\left(
\begin{matrix}
I_p & (L)_{i \alpha} & \frac{1}{2}\left(-\sum_{\gamma}(L)_{\gamma i}(L)_{\gamma j}\right)\\
0   &  I_{q-p}       & -(L)_{i \alpha}^t\\
 0  & 0 & I_p
\end{matrix}
\right)
\end{equation}
\item[\mb{ $\boxed{ Skew(\z^*,\z,B) \leftrightarrow \rho_{\wt{\z^*}} }$  }]\label{eq:SkewZZB} 
To  $S=S_{\wt{\z^*}} \in Skew(\z^*,\z,B)$ we associate $\rho_{\wt{\z^*}} \in \Gamma_{decomp.}$ (cf. page \pageref{dfn:Rho}): 
\[
\wt{Z}_i^*:= \rho_{ \wt{\z^*}}(Z_i^*)=(id_{\z^*}+S)(Z_i^*) = Z_i^* + \sum_k (S)_{ki}Z_k
\]
With respect to the Witt-basis $\{\beta_{\z},\beta_{\z^*}\}$ 
in $\z \oplus \z^*$ the matrix of $Mat_{\{\beta_{\z},\beta_{\z^*}\}}(S)_{ij}=(S)_{ij}$ 
is skew-symmetric. Since $\rho_{\wt{\z^*}}$ acts identically on $\z \oplus W$ we get
\begin{equation}\label{eq:MatRhoZ}
Mat_{\beta}(\rho_{\wt{\z^*}})=
\left(
\begin{matrix}
I_p & 0 & (S)_{ij}\\
 0  &  I_{q-p}       &0\\
 0  & 0 & I_p
\end{matrix}
\right), \qquad (S)_{ij}=-(S)_{ji}
\end{equation}
\end{description}
Recall that $\rho_{\wt{W},\wt{\z^*}}=\rho_{\wt{W}} \circ \rho_{\big( (\rho_{\wt{W}})^{-1}\wt{\z^*} \big)} \in \Gamma_{decomp.}$ 
satisfies $\rho_{\wt{W},\wt{\z^*}}(\z,W,\z^*)=(\wt{\z},\wt{W},\wt{\z^*})$.

\begin{lem}\label{lem:SkewInIso}
Let $(\h \oplus \m,B)$ be an indecomposable solvable symmetric triple.
\begin{align*}
&(i) \quad Skew(\z,\z^*,B) \; \subset \; Aut\big(\m,[[\cdot,\cdot],\cdot],B\big) \; \subset \; \Gamma=stab([\h,\m]) \cap O(\m,B)\\
&(ii) \quad exp\big(\ad(\h)|_{\m}\big) \; \subset \; \Gamma_{decomp.}  
\end{align*}
\end{lem}

\begin{proof}
For the first inclusion of (i) note that $\ad(X)=\ad(\rho_{\wt{\z^*}}X)$ for all $X \in \m$ and $\rho_{\wt{\z^*}}$ associated
to $S_{\wt{\z^*}} \in Skew(\z,\z^*,B)$ and thus for $m_i \in \m$:
\[
B\big([[\rho_{\wt{\z^*}} m_1,\rho_{\wt{\z^*}} m_2],\rho_{\wt{\z^*}} m_3],\rho_{\wt{\z^*}} m_4\big)-B\big([[ m_1, m_2], m_3],m_4\big) = 
B\big(\underbrace{ [[ m_1, m_2], m_3] }_{\in \z \oplus W},\underbrace{ (\rho_{\wt{\z^*}}-id_{\m}) m_4 }_{\in \z} \big)=0
\]
The second inclusion is due to prop. \ref{prop:Aut}(ii) since any $\rho \in Aut(\m,[[\cdot,\cdot],\cdot]B)$ 
extends to a Lie algebra automorphism $\wt{\rho}$
of $\h \oplus \m$ s.t. $\wt{\rho}|_{\m}=\rho$. In particular $\rho(\z)=\z$ hence 
$\rho([\h,\m])=\rho(\z^{\perp})=\rho(\z)^{\perp}=\z^{\perp}=[\h,\m]$.
(ii) follows from the fact that $exp\big(\ad(\h)|_{\m}\big)$ acts identically on $[\h,\m]/\z$.
\end{proof}

Next we define a family of operators in $End(W)$ by restriction of the curvature operators $\RR(Z_i^*, \cdot ) Z_j^*$
onto $W$. This will clarify the structure of the coefficients $f$ and reveal an invariant which plays an essential role 
in the classification.
\begin{lem}\label{lem:Fij}
Consider $F_{ij}:=pr_W \circ \ad(Z_i^*) \circ \ad(Z_j^*)|_W\in End(W),\; Z_i^*,Z_j^* \in \beta_{\z}^*$
for a symmetric triple as in proposition \ref{prop:1.NFMaxCenter}. 
\begin{enumerate}
\item[(i)]All $F_{ij}$ are selfadjoint w.r.t. the positive definite scalar-product $B|_{W \times W}$.
\item[(ii)] $\big[ F_{ij} , F_{kl} \big]= 0 \quad \forall i,j,k,l \in \{1,\ldots,p\} $
\item[(iii)] There exists an orthonormal basis $\beta_{W}=\{W_{\alpha}\}$, 
$\lambda^{\alpha}=(\lambda_1^{\alpha},\dots,\lambda_p^{\alpha}) \in \R^p$ and 
$\epsilon \in \{\pm 1\}^{q-p}$ s.t.
\begin{equation}\label{eq:Fij}
F_{ij}(W_{\alpha}) = \epsilon_{\alpha} \lambda_i^{\alpha}\lambda_j^{\alpha} \cdot W_{\alpha}
\quad \overset{ \ref{prop:1.NFMaxCenter}}{\Longleftrightarrow} \quad
f_{ij\alpha\beta}=\delta_{\alpha\beta} \, \epsilon_{\alpha} \lambda_i^{\alpha} \lambda_j^{\alpha}
\end{equation}
If $\{\wt{W_{\alpha}}\}$ is another orthonormal basis of $W$ 
s.t. $\{\wt{W}_{\alpha}\},\,\wt{\lambda}^{\alpha},\,\wt{\epsilon}_{\alpha}$ satisfies also (\ref{eq:Fij})
then there exists a permutation $\Pi \in S_{q-p}$ s.t. 
\[
\epsilon_{\alpha}=\wt{\epsilon}_{ \Pi(\alpha) }, \quad \lambda^{\alpha} \in \big\{\pm \wt{\lambda}^{ \Pi(\alpha) } \big\}
\]
\item[(iv)]  The coefficients $f_{ij\alpha\beta}$ are invariant under $\Gamma_{decomp.}$.
More precisely, any $\rho \in \Gamma_{decomp.}$ satisfies 
\[
\big[ [\rho Z_j^*,\rho Z_j^* ] , \rho W_\alpha  \big] \equiv \rho \big[ [Z_j^*,Z_j^*] , W_\alpha \big] \; \mod \z
\]
%
\item[(v)] $\dim \bigcap_{i,j} \ker(F_{ij})$ is an invariant of the symmetric space.
\end{enumerate}
\end{lem}

\begin{proof}
(i) By \ref{prop:1.NFMaxCenter} 
$F_{ij}(W_{\alpha})=[Z_i^*,[Z_j^*,W_{\alpha}]]_{W}=[Z_i^*,X_{j\alpha}]={\sum}_{\gamma}f_{ij\alpha\gamma}W_{\gamma}$ for any
orthonormal basis $\beta_{W}=\{W_{\alpha}\}$, 
hence $(F_{ij})_{\alpha\beta}:=\big( Mat_{\beta_{W}}(F_{ij}) \big)_{\alpha\beta} = f_{ij\alpha\beta}$.
Thus $F_{ij}$ is selfadjoint, since $f_{ij\alpha\beta}=f_{ij\beta\alpha}$.

(ii) follows from  (\ref{eq:Kommutativ}) which is itself a consequence of the fact that $X$ is abelian:
{\footnotesize
\begin{multline*}
B \Big( (F_{ij} \circ F_{kl})W_{\alpha},W_{\beta} \Big) \underset{(i)}{=} B \Big( F_{kl}W_{\alpha},F_{ij}W_{\beta} \Big)
 = B \Big([Z_k^*,X_{l\alpha}],[Z_i^*,X_{j\beta}] \Big)
 = B \Big([X_{j\beta},[Z_k^*,X_{l\alpha}]],Z_i^* \Big)= \\
 \underset{Jac.-Id.}{=} B \Big([Z_k^*,\underbrace{[X_{j\beta},X_{l\alpha}]}_{\subset [X,X]=\{0\}}],Z_i^*\Big) +  
   B \Big([X_{l\alpha},[Z_k^*,X_{j\beta}]],Z_i^* \Big)
  = B \Big( (F_{il} \circ F_{kj})W_{\alpha},W_{\beta} \Big) 
\end{multline*}}
Finally $F_{ij}=F_{ji}$ ($f_{ij\alpha\beta}=f_{ji\alpha\beta}$), hence the assertion.

(iii) Since $B|_{W \times W}$ is positive definite the abelian Lie algebra $span\{ F_{ij} \}$ of selfadjoint operators 
can be diagonalized simultaneously i.e. there exists an orthonormal basis $\{W_{\alpha}\}$ s.t. 
$(F_{ij})_{\alpha\beta}=\delta_{\alpha\beta}f_{ij\alpha\beta}$. Defining $f_{ij}^{\alpha}:=f_{ij\alpha\alpha}$
condition (\ref{eq:Kommutativ}) can be reformulated as 
\begin{equation}\label{eq:ikjl}
F_{ik} \circ  F_{jl} = F_{il} \circ F_{jk} \quad \text{i.e} \quad 
f_{ik}^{\alpha} \cdot f_{jl}^{\alpha} = f_{il}^{\alpha} \cdot f_{jk}^{\alpha}
\end{equation}
As eigenvalues of the operators $F_{ij}$ the coefficients $f_{ij}^{\alpha}$ are uniquely determined up to a permutation.
Now fix an $\alpha$:
{\begin{itemize}
\item Neglecting the signs $f_{ii}^{\alpha} \cdot f_{kk}^{\alpha} \overset{(\ref{eq:ikjl})}{=} (f_{ik}^{\alpha})^2$ yields
$|f_{ij}|=\wt{\lambda_i^{\alpha}} \cdot \wt{\lambda_j^{\alpha}} $ where $\wt{\lambda_i^{\alpha}}:= |f_{ii}^{\alpha}|$. 
Furthermore $\epsilon_{\alpha}:=sign(f_{ii}^{\alpha})$ is independent of the index $i$ since $f_{ii}^{\alpha} \cdot f_{jj}^{\alpha} \geq 0$.
\item Set $\epsilon^{\alpha}_{ik}:= \epsilon_{\alpha} \cdot sign(f_{ik}^{\alpha})$. Then
$f_{ii}^{\alpha} \cdot f_{lk}^{\alpha}  \overset{(\ref{eq:ikjl})}{=} f_{il}^{\alpha} \cdot f_{ik}^{\alpha}$
evaluated on the non-vanishing indeces $I^{\alpha}_{\ne 0}:=\menge{1 \le j \le p}{\wt{\lambda}_j^{\alpha}\ne0}$
yields $\epsilon^{\alpha}_{lk}=\epsilon^{\alpha}_{il}\cdot \epsilon^{\alpha}_{ik}$. 
Fix $i_0 \in I^{\alpha}_{\ne 0}$ and define $\epsilon_k^{\alpha} := \epsilon_{i_0 k}^{\alpha}$. Thus 
$\epsilon_{kl}^{\alpha}=\epsilon_{i_0 k}^{\alpha} \cdot \epsilon_{i_0 l}^{\alpha} = \epsilon_k^{\alpha} \cdot \epsilon_l^{\alpha}\;
\forall l,k \in I^{\alpha}_{\ne0}$. 
Finally set 
\[
\lambda_i^{\alpha}= 
\left\{ \rca{1.0}{\ba{ll} \wt{\lambda_i}^{\alpha} \cdot \epsilon_i^{\alpha} & i \in I^{\alpha}_{\ne0} \\ 
0 &  \text{else} \ea } \right. 
\]
Then $f_{ij\alpha\beta}=\delta_{\alpha\beta}\cdot f_{ij}^{\alpha}= 
\delta_{\alpha\beta}\cdot \wt{\lambda_i^{\alpha}} \cdot \wt{\lambda_k^{\alpha}} \cdot 
\epsilon^{\alpha}_{ik} \cdot \epsilon_{\alpha} = 
\delta_{\alpha\beta}\cdot\epsilon_{\alpha} \cdot \lambda_i^{\alpha} \cdot \lambda_j^{\alpha}$ as stated.
\end{itemize}}
Finally we remark that for fixed $\alpha$ the only ambiguity is the sign of the vector $\lambda^{\alpha}$.
Together with the uniqueness of the $f_{ij}^{\alpha}$ up to permutations this proves the assertion. 
However note that choosing a different basis $\{\wt{W}_{\alpha}\}$ which satisfies (\ref{eq:Fij}) might effect the $a_{ijk\alpha}$.

(iv): The adjoint representation induces a map $\wh{\ad}:\m / [\h,\m] \rightarrow \gl \lp [\g,\g]/ \z \rp$ since 
$\big[[\h,\m],[\g,\g]\big]=[W \oplus \z, \h \oplus W \oplus \z] \subset \z$.
Then $\big( \m \cap [\g,\g] \big)/ \z = [\h,m] / \z$ is invariant under the composition $\wh{\ad}(Z_i^*) \circ \wh{\ad}(Z_j^*)$
and we obtain the symmetric bilinear map
\begin{align*}
F:\; \m / [\h,\m] \times \m / [\h,\m]  \quad & \rightarrow \quad \gl \lp [\h,m] / \z \rp \\
(Z_i^*,Z_j^*) \quad  & \mapsto \quad \wh{\ad}(Z_i^*) \circ \wh{\ad}(Z_j^*)
\end{align*}
Under the identification 
$\m / [\h,\m] =\m / \z^{\perp} \simeq\z^*$ and $[\h,\m]/\z = \z^{\perp} / \z \simeq W$ 
clearly $F(Z_i^*,Z_j^*)$ corresponds to $F_{ij}$. This proves the equivariance of the $F_{ij}$ 
and thus the invariance of the $f_{ij}^{\alpha}$ under $\Gamma_{decomp.}$ since $\Gamma_{decomp.}$ acts trivially on the quotients 
$\m / [\h,\m]$ and $[\h,\m]/\z$.

(v) is an immediate consequence of (iv).
\end{proof}

\begin{dfn}[least nilpotent] Let \mb{$W_{nil}:=\bigcap_{i,j}\ker(F_{ij})$} and \mb{$W_{reg}:=W_{nil}^{\perp}$} 
An indecomposable solvable symmetric triple $(\h \oplus \m,B)$ with maximal center is said to be
\emph{least nilpotent} if $W=W_{reg}$ i.e. $\lambda^{\alpha} \neq 0,\;\forall \alpha$.
\end{dfn}
This notion - introduced in \cite[page 343]{CahenP70} - is justified due to the following
\begin{lem}\label{lem:NilpotentMaxCenter}
An indecomposable solvable symmetric triple $(\g,\sigma,B)$ of maximal center is nilpotent
if and only if  $W_{nil}=W$ (w.r.t. an arbitrary adapted decomposition) i.e. $f_{ij\alpha,\beta} = 0$.
In particular every solvable triple of signature $(p,p)$ is nilpotent.
\end{lem}

\begin{proof}
``$\Rightarrow$'': Assume $\lambda_i^{\alpha} \ne 0$ and consider $\wh{\g}=\g/\z$.
Then $(\wh{\ad}\wh{Z}_i^* )^2\wh{W}_{\alpha}=-\epsilon_{\alpha} (\lambda_i^{\alpha})^2\wh{W}_{\alpha} \ne 0$ i.e.
$\wh{\ad}( \wh{Z}_i^*) $ is not nilpotent. Hence $\wh{\g}$ is not nilpotent which contradicts
$\g$ being nilpotent.

``$\Leftarrow$'': If $f = 0$ prop. \ref{prop:1.NFMaxCenter} shows that $\ad(\g)$ is contained
in an upper triangular bloc matrix with respect to the decomposition $\g = \z \oplus X \oplus W \oplus Y \oplus \z^*$
Thus $\g$ is nilpotent due to the theorem of Engel.
\end{proof}

In order to distinguish certain adapted decompositions it remains to investigate how the coefficients
$a$ and $b$ transform under $\Gamma_{decomp.}$. Lets concentrate first on $a$. According to (\ref{eq:MatRhoW}) we have
\begin{equation}\label{eq:Xtilde}
\rca{2}{ \ba{l} \wt{W}_{\alpha}=W_{\alpha} + \sum_k (L)_{k \alpha }Z_k  \\
\wt{Z}_i^* = Z_i^* - \sum_{\gamma} (L)_{i \gamma} W_{\gamma} + \sum_k h_{ki}Z_k \ea}
\quad \Rightarrow \quad
\rca{2}{ \ba{l} \wt{Y}_{ij}=[\wt{Z}_i^*,\wt{Z}_j^*]=Y_{ij}-\sum_{\gamma} (L)_{j\gamma }X_{i \gamma} +\sum_{\gamma} (L)_{i \gamma }X_{j \gamma}\\
  \wt{X}_{k\alpha}=[\wt{Z}_k^*,\wt{W}_{\alpha}]=[{Z}_k^*,{W}_{\alpha}]=X_{k\alpha} \ea}
\end{equation}
hence by (\ref{eq:hMetric})
\begin{multline*}
\wt{a}_{ijk\alpha} = B(\wt{Y}_{ij},\wt{X}_{k\alpha})
=B(Y_{ij},X_{k\alpha}) - \sum_{\gamma} (L)_{j \gamma }B(X_{i \gamma},X_{k \alpha}) 
+\sum_{\gamma} (L)_{i \gamma }B(X_{j \gamma},X_{k \alpha}) = \\
= a_{ijk\alpha} - \sum_{\gamma} (L)_{j \gamma }f_{ik\gamma\alpha}
+\sum_{\gamma} (L)_{\gamma i}f_{jk\gamma\alpha}
= a_{ijk\alpha} -(L)_{j\alpha } f_{ik}^{\alpha} +(L)_{i \alpha } f_{jk}^{\alpha}
\end{multline*}
and thus 
\begin{equation}\label{eq:aTilde}
\wt{a}_{ijk\alpha}= a_{ijk\alpha} + 
  \epsilon_{\alpha} \lambda_k^{\alpha} \cdot \big( -(L)_{j \alpha} \lambda_i^{\alpha} + (L)_{i \alpha} \lambda_j^{\alpha} \big)
\end{equation}

\pagebreak[2]

According to the decomposition $W=W_{nil} \oplus W_{reg}$ two cases has to be distinguished:
\begin{description}
\item[ \mb{$I_{nil}$}$:=\menge{ \alpha }{ W_{\alpha} \in W_{nil} $}] 
Those $a_{ijk\alpha}$ are invariant under $\Gamma_{decomp.}$ due to (\ref{eq:aTilde}):
\begin{equation}
\wt{a}_{ijk\alpha}=a_{ijk\alpha}, \quad \forall \alpha \in I_{nil}
\end{equation}
\item[ \mb{$I_{reg}$}$:=\menge{ \alpha }{ W_{\alpha} \in W_{reg} }$] 
In this case there exists  $\lambda_{n_{\alpha}}^{\alpha} \ne 0$. From the relation (\ref{eq:af}) we see   
$f_{m l}^{\alpha} \cdot a_{ijk\alpha}  =  f_{m k}^{\alpha} \cdot a_{ijl\alpha}$ i.e. 
$\lambda_{m}^{\alpha} \lambda_{l}^{\alpha} a_{ijk\alpha} = \lambda_{m}^{\alpha} \lambda_{k}^{\alpha} a_{ijl\alpha} $.
In particular for $m=l=n_{\alpha}$ 
\[
a_{ijk\alpha}= \lambda_{k}^{\alpha} \cdot \mu_{ij}^{\alpha} \quad \text{with} \quad 
\mu_{ij}^{\alpha}:= \frac{ a_{ijn_{\alpha}\alpha} } { \lambda_{n_{\alpha}}^{\alpha} }
\]
The skew-symmetry $a_{ijk\alpha}=-a_{jik\alpha}$ and the 1. Bianchi-identity $\sigma_{ijk}(a_{ijk\alpha})=0$ then
can be rewritten into $\mu_{ij}^{\alpha}=-\mu_{ji}^{\alpha}$ and $\mu_{ij}^{\alpha} \cdot \lambda_{k}^{\alpha} = 
\mu_{ki}^{\alpha}\lambda_{j}^{\alpha}-\mu_{kj}^{\alpha}\lambda_{i}^{\alpha}$. Hence: 
\[
\mu_{ij}^{\alpha}=\lambda_{i}^{\alpha}\cdot \mu_j^{\alpha}-\lambda_{j}^{\alpha}\cdot \mu_i^{\alpha} \quad \text{with} \quad 
\mu_{l}^{\alpha}:= \frac { \mu_{n_{\alpha}l}^{\alpha} } { \lambda_{n_{\alpha}}^{\alpha} }
\]
Summarizing we get
\begin{equation}\label{eq:aForm}
a_{ijk\alpha}=\lambda_{k}^{\alpha} \cdot \big( \lambda_{i}^{\alpha} \mu_j^{\alpha}-\lambda_{j}^{\alpha} \mu_i^{\alpha}\big),
\quad \forall \alpha \in I_{reg}
\end{equation}
This structure of the $a_{ijk\alpha},\; \alpha \in I_{reg}$ ensures an adapted decomposition s.t. they vanish:
\end{description}

\begin{lem}[Splitting property]\label{lem:afSeperation}
Let $(\g,\sigma,B)$ be an indecomposable solvable symmetric triple with maximal center. 
\begin{enumerate}
\item[(i)] There exists an adapted decomposition $\m = \z \oplus W_{nil} \oplus W_{reg} \oplus \z^*$ 
such that $a_{ijk\alpha}=0 \; \forall  \alpha \in I_{reg}$. Thus the relation (\ref{eq:af}) is trivially satisfied
which can be seen as splitting of $a$ and $f$.
\item[(ii)] For every $\alpha \in I_{nil}$ there exists an $a_{ijk\alpha} \ne 0$.
Furthermore  the coefficients $a_{ijk\alpha},\; \alpha \in I_{nil}$ are invariant under $\Gamma_{decomp.}$.
%
\item[(iii)] If the triple is least nilpotent the set of adapted decompositions 
$\cal{D}_a := \menge{D \in \cal{D}}{a_{ijk\alpha}=0}$ is non-empty due to (i). Then 
\[
\Gamma_{decomp.}^{a}:=\menge{\rho \in \Gamma_{decomp.}}{\rho \big(\cal{D}_a\big) \subset \cal{D}_a} \subset Aut\big(\m,[[\cdot,\cdot],\cdot],B\big)
\]
i.e the adapted decomposition is determined up to an isomorphism of the symmetric triple.
\item[(iv)] The holonomy algebra of a least nilpotent symmetric space is abelian.
\end{enumerate}
\end{lem}

\begin{proof}
(i) Inserting (\ref{eq:aForm}) in (\ref{eq:aTilde}) shows that $\wt{a}_{ijk\alpha} =0,\;\alpha \in I_{reg}$ can
be achieved by the transformation
\[
(L)_{i \alpha }=
\begin{cases}
\epsilon_{\alpha} \cdot \mu_i & \text{if } \alpha \in I_{reg}\\
0 & \text{else}
\end{cases}
\]
(ii) Assume there exists an $\alpha \in I_{nil}$ s.t. $a_{ijk\alpha}=0\; \forall ijk$. Then $\R \cdot W_{\alpha}$
would be a proper non-degenerate ($B|_{W \times W}$ is positive definite) an $\ad(\h)$-invariant subspace of $\m$
which contradicts the  indecomposability of the symmetric triple.
The invariance of the $a_{ijk\alpha},\;\alpha \in I_{nil}$ is clear by (\ref{eq:aTilde}).

(iii) According to lemma \ref{lem:SkewInIso} $Skew(\z,\z^*,B) \subset Aut(\m,[[\cdot,\cdot],\cdot]) \subset \Gamma^a_{decomp.}$ hence
it suffices to determine $\Gamma_{decomp.}^a \cap Hom(W,\z)$.  Let $(\z,W,\z^*) \in \cal{D}_a$.
Since $\lambda^{\alpha} \neq 0$ the transformation rule (\ref{eq:aTilde}) shows
\[
\wt{a}_{ijk\alpha}=0\; \forall i,j,k 
\quad \Leftrightarrow \quad
-(L)_{j \alpha } \lambda_i^{\alpha} + (L)_{i \alpha } \lambda_j^{\alpha} =0
\quad \Leftrightarrow \quad
\exists c_{\alpha} \in \R:\; (L)_{i \alpha } =  c_{\alpha} \cdot \lambda_i^{\alpha} \quad \forall i
\]
hence
\begin{equation}\label{eq:Lalphai}
\Gamma_{decomp.}^a \quad \leftrightarrow \quad \Menge{ L \in Hom(W,\z)}{L(W_{\alpha})=c_{\alpha} \cdot {\sum}_i \lambda_i^{\alpha} Z_i } 
\; \simeq \; \R^{q-p}
\end{equation}
Now $\Gamma_{decomp.}^a  \subset Aut\big(\m,[[\cdot,\cdot],\cdot],B\big)$ can be proved by showing the invariance of the $b_{ijkl}$:
\begin{equation}\label{eq:bTilde}
\wt{b}_{ijkl} \underset{ (\ref{eq:hMetric}) }{=} B( \wt{Y}_{ij} , \wt{Y}_{kl} )
 \underset{(\ref{eq:Xtilde})}{=}   b_{ijkl} +
\sum_{\gamma} \Big( (L)_{j \gamma } (L)_{l \gamma } f_{ik}^{\gamma} +
                    (L)_{i \gamma } (L)_{k \gamma } f_{jl}^{\gamma} -
                    (L)_{j \gamma } (L)_{k \gamma } f_{il}^{\gamma}  -
                    (L)_{i \gamma } (L)_{l \gamma } f_{jk}^{\gamma} \Big) 
\end{equation}
since $a = 0$ implies $Y \perp_B X$.  
Applying (\ref{eq:Lalphai}) and \ref{lem:Fij}(iii) shows that each summand in the large bracket of (\ref{eq:bTilde}) is of the form 
$(c_{\gamma})^2 \epsilon_{\gamma} \lambda_{i}^{\gamma}\lambda_{j}^{\gamma}\lambda_{k}^{\gamma}\lambda_{l}^{\gamma}$,
thus they cancel each other due to their signs. It remains $ \wt{b}_{ijkl}=b_{ijkl}$. 

(iv) This is only the observation that $a = 0$ in this case.
\end{proof}

{\bf The case \mb{$W_{nil} \ne \{0\}$} :}\label{case:f=0}
In order to find a parametrisation for those symmetric triples, it seems to be necessary to determine all
$(a_{ijk\alpha}) \in \R^p \wedge \R^p \times \R^p \times \R^{\#(I_{nil})}$ which satisfy $\sigma_{ijk}(a_{ijk\alpha})=0$ and (\ref{eq:a}). 
\begin{itemize}
\item For $q-p=1$ it is known, that all such families of coefficients $(a_{ijk\alpha})$ are of the form 
(\ref{eq:aForm}).
Then (\ref{eq:a}) is trivially fulfilled.
Moreover there exists a construction of coefficients $(a_{ijk\alpha})$ satisfying (\ref{eq:aForm}) in arbitrary signature
modeled on a Riemannian symmetric space (cf. \cite{CahenP70}).
\item If $p \le 2$ the condition (\ref{eq:a}) is again empty. 
\item According to lemma \ref{lem:afSeperation} $a = 0$ contradicts the indecomposability. 
\item The signature restricts the dimension of $W_{nil}$. More precisely: 
Let  $Y=span\{Y_1,\ldots,Y_m\}$. 
$Y \ne \{0\}$ since otherwise $W_{nil}$ would be non-singular and $\ad(\h)$-invariant due to $[X,W_{nil}]=\{0\}$). 
Hence  $1 \le \dim(Y)=m \le p(p-1)/2$.
Each kernel $K_{i}=\ker(\ad(Y_{i})|_{W_{nil}})$ is a subspace of dimension greater than $\dim(W_{nil})-\dim(\z)$ (dimension formula for linear maps)
and thus $\dim \big( \bigcap_{i=1}^m K_i \big) \ge \dim(W_{nil}) - m \cdot p$. 
Since every $w \in \bigcap_{i=1}^m K_i$ would span a non-singular $\ad(\h)$-invariant subspace of $\m$, 
indecomposability implies $\bigcap_{i=1}^m K_i = \{0\}$, hence $\dim(W_{nil}) \le m \cdot p$. 
\begin{equation}\label{eq:Dimf=0}
\dim(W_{nil}) \le p^2(p-1)/2
\end{equation}
Again we see that the complexity increases with the signature.
\end{itemize}
The last three points enable a classification in signature $(2,n-2)$, whereas the case of arbitrary signature is open. In particular
it would be interesting to know if there are examples with non-abelian holonomy.

\subsubsection[Classification of least nilpotent symmetric triples]{Classification of least nilpotent symmetric triples}

In contrast to the case $W_{nil} \neq \{0\}$ the structure of a least nilpotent symmetric triple is encoded in simple algebraic objects.
First we describe the construction of a normal form for predicted data:

\pagebreak
\begin{dfn}[\mb{$\tau_{p,q}(\RR_b,\epsilon,\Lambda)$}]\label{dfn:TripleLeastNilpotent}
Let $\cal{N}_{p,q}$ $(1 \le p \le q)$ be the class of the following symmetric triples 
$\tau_{p,q}(\RR_b,\epsilon,\Lambda)=(\h \oplus \m,B)$: 
\begin{itemize}
\item $(W,B_W)$ a $(q-p)$-dim. Euclidian vector space with a fixed orthonormal basis $\{W_{\alpha}\}$.
\item $(\z \oplus \z^*,B_{\z})$ with $\z$ an $p$-dim. vector space, $\z^*$ its dual space s.t. $\z,\z^*$ are singular paired w.r.t $B_{\z}$, i.e. 
\[
B_{\z}|_{\z \times \z}=B_{\z}|_{\z^* \times \z^*}=0 \qquad B_{\z}(Z^*,Z)=Z^*(Z) \quad \forall Z \in \z,Z^* \in \z^*
\]
\item $\RR_b \in \Re(\z^*) \subset S^2\big(\Lambda^2(\z^*)\big)$ an algebraic curvature tensor on $z^*$. 
Consider the canonical projection $\pi: \Lambda^2(\z^*) \rightarrow {}^{\ds{\Lambda^2(\z^*)}}/_{\ds{\ker(\RR_b)}} = :Y$.
Then $\RR_b$ factors through $\pi$ and yields the pseudo-Euclidian space $(Y,B_Y)$ with $B_Y:=\pi_*(\RR_b)$.
Further consider the dualized map $b:Y \times \z^* \rightarrow \z$ defined by
\[
B_{\z}\big( b( \pi(Z_1^* \wedge Z_2^*),Z_3^*),Z_4^*\big)=\RR_b(Z_1^* \wedge Z_2^*,Z_3^*)\wedge Z_4^*)
\]
\item $(W^*,B_{W^*})$ is a pseudo-Euclidian vector space with $\dim(W^*)=\dim(W)$ and a fixed pseudo-orthonormal basis
$\{W_{\alpha}^*\}$, i.e. $B_{W^*}(W_{\alpha}^*,W_{\alpha}^*)=\epsilon_{\alpha}$, $\epsilon =(\epsilon_{\alpha}) \in \{\pm 1\}^{q-p}$.
\item $\Lambda_1,\ldots,\Lambda_{q-p} \in \z \setminus \{0\}$. Set
$\m:=\z \oplus W \oplus \z^*$ and $\h:=Y \oplus W^*$. Then $\h \oplus \m$ carries the non-generate bilinear form 
$B:= B_Y \oplus B_{W^*} \oplus B_{\z} \oplus B_{W}$ and the 
following Lie algebra structure (only the brackets which are possibly non-zero are cited):
\begin{align*}
[\m,\m]: \quad &\left\{ \hspace{.15cm}
\rca{1.3}{\ba{rcl}
{[Z_1^*,Z_2^*]} &=& \pi(Z_1^* \wedge Z_2^*)\\
{[Z^*,W_{\alpha}]} &=& B(Z^*,\Lambda_{\alpha}) \cdot W_{\alpha}^*
\end{array}}
\right. \\
[\h,\m]:\quad &\left\{
\rca{1.3}{\ba{rcl}
{[y,Z^*]}&=&b(y, Z^*)\\
{[W_{\alpha}^*,Z^*]}&=&\epsilon_{\alpha} \cdot B(\Lambda_{\alpha},Z^*) \cdot W_{\alpha}\\
{[W_{\alpha}^*,W_{\beta}]}&=&-\epsilon_{\alpha}\delta_{\alpha \beta} \cdot \Lambda_{\alpha}
\end{array}}
 \right.
\end{align*}
\item We require $[Y,\z^*]+[W^*,W]=b(Y,\z^*)+span\{\Lambda_{\alpha}\}=\z$ which is equivalent to $[\h,\m]=\z \oplus W$.
\end{itemize}
\end{dfn}

\begin{thm}\label{thm:Main}
Let $\mathcal{T}_{p,q}$ $(1 \le p \le q)$ be the class of indecomposable 
solvable symmetric triples of signature $(p,q)$ with maximal non-zero center which are least nilpotent.
\begin{enumerate}
\item[(i)] Every $\tau \in \mathcal{T}_{p,q}$ is isomorphic to a $\tau_{p,q}(\RR_b,\epsilon,\Lambda) \in \mathcal{N}_{p,q}$.
\item[(ii)] Every $\tau_{p,q}(\RR_b,\epsilon,\Lambda)\in \mathcal{N}_{p,q}$ is isomorphic to a sum $\tau_1 \oplus \ldots \oplus \tau_s$,
where $\tau_i \in \mathcal{T}_{p_i,q_i}$ $(1 \leq p_i \leq q_i)$.
A sufficient condition for $\tau_{p,q}(\RR_b,\epsilon,\Lambda)$ being indecomposable is $\dim(Y)>p(p-2)/2$.
\item[(iii)] Two indecomposable triples $\tau=\tau_{p,q}(\RR_b,\epsilon,\Lambda), 
\wt{\tau}=\tau_{p,q}(\wt{\RR}_b,\wt{\epsilon},\wt{\Lambda}) \in \mathcal{N}_{p,q}$
are isomorphic if and only if there exists a permutation $\Pi \in \mathcal{S}_{q-p}$ and $P \in Gl(\z)$ such that:
\begin{equation}\label{eq:IsoCond}
P^* \cdot \RR_b  = \wt{\RR_b}, \qquad 
\epsilon_{\alpha} =\wt{\epsilon}_{\Pi(\alpha)}, \qquad 
P (\Lambda_{\alpha}) =  \pm \wt{\Lambda}_{\Pi(\alpha)} 
\end{equation}
\item[(iv)] For $\tau=\tau_{p,q}(\RR_b,\epsilon,\Lambda)$ its automorphism group 
$Aut(\tau):=Aut\big(\m,[[\cdot,\cdot],\cdot],B\big) \subset \Gamma_{decomp.} \rtimes \Gamma_{basis}^D$ 
corresponds under the identification (\ref{eq:one-to-one}) to 
\begin{align*}
Aut(\tau) \cap \Gamma_{decomp.} \quad & \leftrightarrow \quad Skew(\z,\z^*,B)  \times  
\Menge{L \in Hom(W,\z)}{L(W_{\alpha})=c_{\alpha}\cdot \Lambda^{\alpha},\; c \in \R^{q-p}} \\
Aut(\tau) \cap \Gamma_{basis}^D \quad & \leftrightarrow \quad 
\ba{l} \Menge{P \in GL(\z)}{P(\Lambda^{\alpha}) = \pm \Lambda^{\alpha},\; P^* \cdot \RR_b=\RR_b} \times \ldots \\
\hspace{2cm} \ldots \times O(W_1,B|_{W_1 \times W_1}) \times \ldots \times  O(W_s,B|_{W_s \times W_s}) \ea
\end{align*}
where $W_i$ are the simultaneous eigenspaces of the operators $pr_W \circ \ad(Z^*)^2|_W$. 
Moreover
\begin{align*}
exp(\ad(Y)|_{\m}) \quad  &\subset \quad  Skew(\z,\z^*,B) \quad \text{with equality iff $\dim(Y)=p(p-1)/2$} \\ 
exp(\ad(W^*)|_{\m})  \quad  & \simeq  \quad  \R^{q-p} \subset Hom(W,\z) \quad \text{as above}  
\end{align*}

\end{enumerate}
\end{thm}

Assertion (i) of this theorem is essentially \cite[prop. 1]{CahenP70}. Although in a strict sense
the theorem cannot be regarded as a classification of $\mathcal{T}_{p,q}$
due to the lack of a necessary and sufficient criterion for indecomposability of elements
in $\mathcal{N}_{p,q}$, the condition $[\h,\m]=\z \oplus W$ guarantees at least
that $\mathcal{N}_{p,q}$ is closed under decomposition which can be seen as 
an "inductive classification".

\begin{proof}
(i) Let $(\m \oplus \h,B) \in \cal{T}_{p,q}$. Then by the splitting property \ref{lem:afSeperation}(iii) 
there is up to an isomorphism a unique adapted decomposition $D$ such that for a suitable adapted basis
as in lemma \ref{lem:Fij} the only non-trivial coefficients are $(\lambda_i^{\alpha})$ and $(b_{ijkl})$. 
The $(\lambda_i^{\alpha})$ can be replaced by 
\[
\Lambda^{\alpha}:=\sum_k \lambda_k^{\alpha} \cdot Z_k \in \z
\]
since $\lambda_i^{\alpha}=B(\Lambda^{\alpha},Z_i^*)$.
It remains to investigate the structure of $\h$.

\begin{itemize}
\item for $X \subset \h$ the linear dependency among its generators $X_{k \alpha}$
can be explicitely determined. Since $I_{reg}=\{1,\ldots,q-p\}$ we choose for every $\alpha$ an index $k_{\alpha}$ 
such that $\lambda_{k_{\alpha}}^{\alpha} \ne 0$ and set
\[
W_{\alpha}^* := \frac{1}{\lambda_{k_{\alpha}}^{\alpha}} X_{k_{\alpha} \alpha}
\]
\begin{enumerate}\label{eq:Walpha*}
\item[a)] 
$B(W_{\alpha}^*,W_{\beta}^*)= 
(\lambda_{k_{\alpha}}^{\alpha} \cdot \lambda_{k_{\beta}}^{\beta})^{-1} \cdot B(X_{k_{\alpha} \alpha},X_{k_{\beta} \beta})
\underset{(\ref{eq:hMetric}), \ref{lem:Fij}(iii)}{=} \delta_{\alpha \beta}\cdot \epsilon_{\alpha}$
\item[b)]
From 
$B \big( X_{i \alpha}\cdot \lambda_j^{\alpha} - X_{j \alpha}\cdot \lambda_i^{\alpha} \;,\; Y_{kl} + X_{n \gamma}\big)  
 \underset{(\ref{eq:hMetric})}{=}    
\delta_{\alpha \gamma} \epsilon_{\alpha}\lambda_n^{\alpha} 
\cdot ( \lambda_i^{\alpha} \lambda_j^{\alpha} - \lambda_j^{\alpha} \lambda_i^{\alpha})=0 $ follows 
$X_{i \alpha}\cdot \lambda_j^{\alpha} =X_{j \alpha}\cdot \lambda_i^{\alpha}$ since  $B|_{\h \times \h}$ is non-degenerate.
Inserting the definition of $W_\alpha^*$ into this relation yields 
\[
X_{i \alpha}= \lambda_{i}^{\alpha} \cdot W_{\alpha}^*
\]
\item[c)] The $W_{\alpha}^*$ are well defined i.e. independent of the chosen index $k_{\alpha}$ with 
$\lambda_{k_{\alpha}}^{\alpha} \ne 0$.
This follows immediately from b).
\end{enumerate}
Hence $X = \bigoplus_{\alpha}^{\perp} \R \cdot W_{\alpha}^*$ is a non-singular subspace of 
$\h$ (w.r.t. $B|_{\h \times \h}$) spanned by the pseudo-orthonormal basis $\{W_{\alpha}^*\}$. Moreover 
$X \cap Y = \{0\}$ due to $X \perp Y$ and $X$ being non-singular. 
In \ref{dfn:TripleLeastNilpotent} we renamed $W^*:=X$.
 
\item 
For $Y \subset \h$ we will be content with a slight reformulation. Since
$b_{ijkl} = \RR(Z_i^*,Z_j^*,Z_k^*,Z_l^*)$
the coefficients $(b_{ijkl})$ correspond in a unique way to an
algebraic curvature tensor 
$\RR_{b} \in \Re (\z^*) \subset S^2 \big( \Lambda^2(\z^*) \big)$.
Then $[Z^*_1,Z^*_2]=0 \; \Leftrightarrow \; \ad([Z^*_1,Z^*_2])=0 
\; \Leftrightarrow \; Z^*_1 \wedge Z^*_2 \in ker (\RR_b)$. Thus 
the pseudo-Euclidian space $(Y,B|_{Y \times Y})$ is obtained from $\big( \Lambda^2(\z^*) ,\RR_b \big) $
by taking the quotient over $\ker(\RR_b)$.
\end{itemize}

(ii) By construction every $\tau_{p,q}(\RR_b,\epsilon,\Lambda)=(\h \oplus \m,B)$ is a solvable symmetric triple (cf. \ref{re:LAConstruction}). 
Thus we can decompose it according to \ref{prop:Decomposability}:
$\tau_{p,q}(\RR_b,\epsilon,\Lambda)=\tau_0 \oplus \ldots \oplus \tau_s$ with
$\tau_i=(\h_i \oplus \m_i,B_i)$, $\Sig(B_i|_{\m_i^2})=(p_i,q_i)$. By assumption $[\h,\m]=\z \oplus W$, hence 
$\z(\h \oplus \m)=[\h,\m]^{\perp} \cap \m = \z$ according to  \ref{lem:Center}(iii). Further note that
\begin{equation}\label{eq:CenterSum}
\z=\z(\h \oplus \m)\underset{(\ref{lem:Center})}{=} \m \cap [\h,\m]^{\perp}= \bigoplus_{i=0}^s \m_i \cap [\h_i,\m_i]^{\perp} = 
\bigoplus_{i=0}^s \underbrace{\z(\h_i \oplus \m_i)}_{=:\z_i} 
\end{equation}
Every $\tau_i$ is solvable since $\h_i \oplus \m_i \lhd \h \oplus \m$ and $\dim(\m_i) > 1$ since 
otherwise $\m_i \subset \z(\h \oplus \m)=\z$ would be a non-singular subspace
of the totally isotropic space $\z$.  Thus by lemma \ref{lem:Center}(iii) $\z_i$ is totally isotropic and non-zero, hence 
$1 \le \dim(\z_i) \le \min\{p_i,q_i\}$. On the other hand
\[
\sum_{i=1}^s p_i = p=\dim(\z) \underset{(\ref{eq:CenterSum})}{=}\sum_{i=1}^s \dim(\z_i)  \le \sum_{i=1}^s  \min\{p_i,q_i\} \le
\sum_{i=1}^s p_i  \quad \Rightarrow \quad  \dim(\z_i) = p_i \leq q_i
\]
Thus every $\tau_i$ has maximal center. Finally we observe that every $\tau_i$ is least nilpotent:
We may assume that all $\tau_i$ are given in the form guaranteed by \ref{prop:1.NFMaxCenter} and denote the corresponding adapted basis
by $\beta_i=\{\beta_{\z_i},\beta_{W_i},{\beta_{\z_i}^*}\}$. Then $\tau$ is in a canonical way given in the form 
\ref{prop:1.NFMaxCenter} by choosing the basis $\beta=\big\{ \{\beta_{\z_i}\}_i,\{\beta_{W_i}\}_i,\{\beta_{\z_i}^*\}_i\big\}$.
Assume one of the triples, say $\tau_i$ is not least nilpotent, i.e. there exists an element $\wh{W} \in W_i$ s.t.
$[\z_i^*,[\z_i^*,\wh{W}]] \subset \z_i$. Then $[\g_i,\g_j]=\{0\}\; (i \ne j)$ implies 
$[\z^*,[\z^*,\wh{W}]] \subset \z$ which contradicts $\tau$ being itself least nilpotent. 

To proof the sufficient criterion concerning indecomposability assume $\tau \simeq \tau_1 \oplus \ldots \oplus \tau_s$.
Then it is easily verified that $\sum_i \dim(Y_i)=\dim(Y)$. Since $Y_i$ is a factor space of $\z_i^* \wedge \z_i^*$,
its dimension is bounded from above by $p_i(p_i-1)/2$. Using $\sum_i p_i=p$ yields 
\[
\sum_{i=1}^s \dim(Y_i) \le \sum_{i+1}^s p_i(p_i-1)/2 
= \frac{1}{2} \Big( p^2 - p - \sum_{i \ne j} p_i p_j \Big) 
\; \underset{\stru{2}(\text{if } s \ge 2)}{\le} \; p(p-2)/2  
\]
Thus $\dim(Y) > p(p-2)/2$ implies $s=1$ i.e. $\tau$ is indecomposable in this case.

(iii) If (\ref{eq:IsoCond}) is satisfied a straight forward computation shows that 
$P \oplus Q_{\Pi} \oplus {P^*}^{-1} \in GL(\z \oplus W \oplus \z^*)$
is the required isomorphism between $\tau$ and $\wt{\tau}$ where $Q_{\Pi}(W_{\alpha})=W_{\Pi(\alpha)}$. 

Conversely assume there exists an isomorphism $\Phi: \tau \rightarrow \wt{\tau}$, i.e. $\Phi \in O(\z,W,\z^*)$
and $\Phi \circ \ad( \cdot ) = \wt{\ad}\big(\Phi (\cdot)\big) \circ \Phi$ (denoting the adjoint representation of $\wt{\tau}$ by  $\wt{\ad}$). 
Since $(\z,W,\z^*),\Phi^{-1}(\z,W,\z^*) \in \cal{D}_a(\tau)$ are both adapted decompositions which satisfy $a=0$ in $\tau$
we may modify $\Phi$ by composition with an isomorphism $\rho \in \Gamma^a_{decomp.}(\tau)$ 
s.t. it leaves the decomposition $(\z,W,\z^*)$ invariant (lemma \ref{lem:afSeperation}(ii)).
Thus $\Phi$ is of the form $\Phi = P \oplus Q \oplus {P^*}^{-1} \in GL(\z) \times O(W,B|_{W\times W}) \times Gl(\z^*)$.

Furthermore let $F(Z^*) := pr_W \circ \ad(Z^*)^2|_W$ and consider the decomposition $W=W_1 \oplus \ldots \oplus W_s$
into simultaneous eigenspaces of the operators $F(Z^*)$ (analog notation for $\wt{\tau}$). Without loss of generality
we may assume that $W_i=span\{W_{\alpha}\}_{\alpha_i \leq \alpha \leq \alpha_{i+1}}$. 
Since $\Phi(W_i) = \wt{W_j}$ and since $W_{\alpha}$ is an eigenvector of both $F(Z^*)$ and $\wt{F}(Z^*)$ their exists a permutation
$\Pi \in S_{q-p}$ s.t. $\Phi(W_i) = span\{W_{\Pi(\alpha)}\}_{\alpha_i \leq \alpha \leq \alpha_{i+1}}$. 
Since $O(W_1,B_{W_1 \times W_1}) \times \ldots \times O(W_1,B_{W_1 \times W_1}) \subset Aut(\tau)$ another modification of $\Phi$ by
an isomorphism of $\tau$ yields $\Phi(W_{\alpha})=W_{\Pi(\alpha)}$. Now the relations (\ref{eq:IsoCond}) can be obtaind immediately
by using the fact that $\Phi$ is an Lie algebra isomorphism.

(iv) The assertion on $Aut(\tau) \cap \Gamma_{decom.}$ is equivalent to lemma \ref{lem:afSeperation}(iii).
As far as $Aut(\tau) \cap \Gamma_{basis}^D$ is concerned we are in the situation of (iii) with two identical triples.
Moreover the holonomy action can be explicitely computed. For $y \in Y$ we have
$\ad(y)|_{\z^*}(\cdot)=b(y,\cdot) \in Skew(\z^*,\z,B)$ and  $\ad(y)|_{\z \oplus W}=0$ hence
\[
Exp(\ad(y)|_{\m})(Z+W+Z^*)= Z+W+Z^* + b(y,Z^*) 
\]
Writing this map down in a basis shows that it is of the form (\ref{eq:MatRhoZ}) with $(S)_{ij}=b_{klij}$ for $y=Y_{kl}$.
For $W_{\alpha}^* \in W^*$ we obtain
\[
\rca{1.5}{\ba{l}
Exp(\ad(W_{\alpha}^*))(Z+W+Z^*)=\\
\qquad =Z+W+Z^* -\epsilon_{\alpha}B(W_{\alpha},W)\Lambda_{\alpha} + \epsilon_{\alpha}B(\Lambda_{\alpha},Z^*)W_{\alpha} -
\frac{1}{2}B(\Lambda_{\alpha},Z^*)\Lambda_{\alpha} 
\ea}
\]
Writing this map down in a basis shows that it is of the form (\ref{eq:MatRhoW}) with
$(L)_{i\beta}= \delta_{\alpha \beta}\epsilon_{\alpha}\lambda_i^{\alpha}$.
\end{proof}

\pagebreak

\subsubsection{Classification of Lorentzian solvable symmetric triples}

Indecomposable Lorentzian symmetric spaces split into three classes (cf. page \pageref{fig:tree}). 
By the separation property \ref{cor:SeperationProperty}(ii) the transvection group is either solvable or semi-simple.
Semi-simple symmetric spaces can be divided in those with irreducible or non-irreducible holonomy.
Since the latter case implies split signature the Lorentzian symmetric space has dimension 2
and thus constant sectional curvature. If the holonomy acts irreducibly the same result holds.
This was first observed by inspecting Berger's list \cite{Berger57} of semi-simple symmetric spaces and computing their
signatures (cf. \cite[Thm. 2]{CahenL90}). A conceptual proof was given in \cite[Thm. 1.5]{DiScala01}.

The indecomposable solvable Lorentzian  symmetric spaces has been classified by Cahen and Wallach 
in \cite{CahenW70},\cite{Cahen72},\cite{Cahen98}. For dimensions $\geq 2$ their center is
non-zero and totally isotropic due to \ref{cor:SeperationProperty}(iv) and \ref{lem:Center}(iii). 
Thus every solvable Lorentzian symmetric triple $\tau=(\h \oplus \m,B)$
has maximal center and we can apply
\ref{prop:1.NFMaxCenter}. In particular $Y=\{0\}$ since $\z^*$ is $1$-dimensional. Thus $a =  0, \; b =  0$
and according to lemma \ref{lem:afSeperation}(ii)
$W_{nil}=\{0\}$, i.e. $\tau$ is least nilpotent. Thus we are lead to theorem \ref{thm:Main}.

\begin{dfn}[\mb{$\tau_{n}(f)$}]\label{dfn:tau(f)}
Let $\mathcal{N}_{Lor}^{n}$ be the class of the following symmetric triples $\tau_{n}(f)=(\h \oplus \m,B)$ depending on 
parameters $f =(f^1,\ldots,f^{n-2}) \in (\R \setminus \{0\})^{n-2}$:
\bi
\item $\m=\z \oplus W \oplus \z^*$, $\z=\R\cdot Z$, $\z^*=\R\cdot Z^*$ and $W = span\{W_1,\ldots,W_{n-2}\}$
\item $\h=span\{X_1,\ldots,X_{n-2}\}$
\item $\z \oplus \z^*,\; W, \h$ are mutual orthogonal w.r.t. $B$ and 
\begin{align*}
 &B(Z,Z) = B(Z^*,Z^*)=0 & &B(W_{\alpha},W_{\beta}) = \delta_{\alpha\beta}\\
 &B(Z^*,Z) = 1          & &B(X_{\alpha},X_{\beta}) = \delta_{\alpha\beta} f^{\alpha}
\end{align*}
\item The non-zero commutators of the Lie algebra structure are 
\[
\rca{1.3}
{\ba{rcl}
{[Z^*,W_{\alpha}]} &=& \phantom{-}X_{\alpha}\\
{[X_{\alpha},Z^*]} &=& \phantom{-}f^{\alpha} \cdot W_{\alpha}\\
{[X_{\alpha},W_{\beta}]}&=& -\delta_{\alpha\beta} f^{\alpha} \cdot Z
\ea}
\]
\ei
\end{dfn}

\begin{prop}\label{triple:tau(f)}
Let $\mathcal{T}_{Lor}^{n}$ $(n \geq 3)$ be the class of $n$-dimensional indecomposable solvable Lorentzian symmetric triples.
\begin{enumerate}
\item[(i)] $\mathcal{N}_{Lor}^{n} \subset \mathcal{T}_{Lor}^{n} $ and for every $\tau \in \mathcal{T}_{Lor}^{n}$
exists an $f \in (\R \setminus \{0\})^{n-2} $ s.t.  $\tau$ is isomorphic to $\tau_{n}(f)$.

\item[(ii)] Two triples $\tau_{n}(f),\tau_{n}(\wt{f}) \in \cal{N}_{Lor}^n$ are isomorphic if there exists a permutation 
$\Pi \in \mathcal{S}_{n-2}$ and
$c \in \R^+$, such that $\forall \alpha,\; f^{\alpha}=c \cdot \wt{f}^{\phantom{,}\Pi(\alpha)}$.
\item[(iii)] If $E_i$ $(1 \le i \le s)$ denotes the eigenspaces of $pr_W \circ \ad(Z^*)^2|_{W}$ for  a triple $\tau=\tau_n(f)$ then
\[
Aut(\tau) \;  \simeq  \; Exp(\ad(\h)|_{\m}) \; \rtimes \;   
\Big( \{\pm id_{\z \oplus \z^*}\} \times O(E_1,B|_{E_1 \times E_1}) \times \ldots \times  O(E_s,B|_{E_s \times E_s}) \Big)
\]
Moreover $\mathfrak{aut}(\tau)$ admits a proper Levi decomposition if and only if $\max(\dim(E_i)) > 2$.
The radical $\rad \big( \mathfrak{aut}(\tau)\big)$ is abelian if $\forall i,\;\dim(E_i)\ne 2$, otherwise still nilpotent.
\end{enumerate}
\end{prop}
\nopagebreak

\begin{re}\label{re:LorentzDim2}\be
\item[(i)] A solvable Lorentzian symmetric triple of dimension two is flat, hence decomposable. This is because
$\z \ne \{0\}$ and thus $\h =\{0\}$.  
\item[(ii)] The classification of indecomposable solvable Lorentzian symmetric space derives as well 
from a more general class of symmetric triples considered in \ref{section:(III)}.
\ee
\end{re}

\begin{proof}
(i): As mentioned above every triple $\tau \in \mathcal{T}_{Lor}^n$ has a maximal center and is least nilpotent.
Hence, with the notation of theorem \ref{thm:Main} $\mathcal{T}_{Lor}^n=\mathcal{T}_{1,n-1}$.
Now we will show that for any triple in $\mathcal{N}_{Lor}^n$ there is an isomorphic triple in  $\mathcal{N}_{1,n-1}$ and vice versa.
Moreover every triple in $\mathcal{N}_{1,n-11}$ is indecomposable due to  \ref{thm:Main}(ii). Together with  \ref{thm:Main}(i) this yields the
assertion.

Consider $\tau_{1,n-1}(\R_b,\epsilon,\Lambda) \in \mathcal{N}_{1,n-1}$ with $\Lambda_1,\ldots,\Lambda_{n-2} \in \z \setminus \{0\}$ and 
$\epsilon \in \{\pm 1\}^{n-2}$. Since $\z^*$ is 1-dimensional, $\RR_b = 0$. Hence the Lie algebra structure of 
$\tau_{1,n-1}(0,\epsilon,\Lambda)$ is given by
\begin{align*}
[Z^*,W_{\alpha}]&=\phantom{-}B(Z^*,\Lambda_{\alpha}) \cdot {W}_{\alpha}^*\\
[{W}_{\alpha}^*,Z^*]&=\phantom{-}\epsilon_{\alpha}B(Z^*,\Lambda_{\alpha})\cdot W_{\alpha} \\
[{W}_{\alpha}^*,W_{\beta}] &= - \epsilon_{\alpha}\delta_{\alpha\beta} \cdot \Lambda_{\alpha}
\end{align*}
Now we replace $W_{\alpha}^*,\Lambda_{\alpha}$ by the original data $X_{\alpha},f^{\alpha}$ as gained in \ref{prop:1.NFMaxCenter}. 
Since $p=1$ the coefficients $f_{ij\alpha\beta}$ reduces to 
$f_{11\alpha\beta}=\delta_{\alpha\beta}f_{11}^{\alpha}=\delta_{\alpha\beta} \,\epsilon_{\alpha}
(\lambda_1^{\alpha})^2$. Thus it is not advantageous to split $f_{11}^{\alpha}$ into a sign and
a square root which causes even ambiguity. But as far as the other assertions are concerned it is quite convenient 
to have an ``embedding'' of the Lorentzian case into the more general one of theorem \ref{thm:Main}. Hence set
\[ 
X_{\alpha}=B(Z^*,\Lambda_{\alpha})W_{\alpha}^* \qquad 
f_{\alpha}=\epsilon_{\alpha}B(Z^*,\Lambda_{\alpha})^2 \qquad 
X=span\{X_{\alpha}\}
\]
It is easy to check that with this substitutions the above triple transfers into the triple $\tau_n(f)$
as cited in the proposition. Hence for every triple in $\mathcal{N}_{1,n-1}$ can be found an isomorphic one in
$\mathcal{N}_{Lor}^n$ and conversely. 

(ii): We know from theorem \ref{thm:Main}(iii), that $\tau_{1,n-1}(0,\epsilon,\Lambda) \simeq \tau_{1,n-1}(0,\wt{\epsilon},\wt{\Lambda})$
if and  only if $\exists \;\Pi \in \mathcal{S}_{n-2},\; c \in \R \setminus \{0\}$ such that
$\epsilon_{\Pi(\alpha)}=\wt{\epsilon_{\alpha}}$ and $c \cdot \Lambda_{\Pi(\alpha)} \in \{\pm \wt{\Lambda}_{\alpha}\}$.
Translating this into the coefficients $f$ yields:
\[
\tau_n(f) \simeq \tau_n(\wt{f})
\qquad  \Leftrightarrow \qquad   \exists \;\Pi \in \mathcal{S}_{n-2}, \;c \in \R \setminus \{0\} \text{ such that }
c^2 \cdot f^{\Pi(\alpha)} =  \wt{f}^{\alpha}
\]
since $\wt{f}^{\alpha}=\wt{\epsilon}_{\alpha}B(Z^*,\wt{\Lambda}_{\alpha})^2 =
\epsilon_{\Pi(\alpha)} \, c^2B(Z^*,\Lambda_{\Pi(\alpha)})^2 = c^2 f^{\alpha}$.

(iii): We apply \ref{thm:Main}(iv) to the current situation:
\bi
\item  $Aut\big( \tau_n(f) \big) \cap \Gamma_{decomp.}=exp \big( \ad(\h)|_{\m} \big)$ since  $\dim(\z)=1$, hence $Skew(\z^*,\z,B)=\{0\}$.
\item  $Aut\big( \tau_n(f) \big) \cap \Gamma_{basis}^D=\menge{c \in \R}{c \cdot \Lambda_{\alpha} \in \{\pm \Lambda_{\alpha}\}} = \{\pm 1\}$.
The induced maps are $\left(\pm id_{\z \oplus \z^*}\right) \oplus id_{W}$.
\end{itemize}
Descending to the Lie algebra level yields:
\[
\mathfrak{aut}(\m,B,\RR) \simeq \h \rtimes
\left( \bigoplus_{i=1}^s \mathfrak{o}(E_i,B|_{W_i \times W_i}) \right)
\]
Further $\mathfrak{o}(E_i,B|_{E_i \times E_i})$ is semi-simple if and only if $\dim(E_i) > 2$.
In particular if $\dim(E_i) \ne 2 \; \forall i$, $\rad(\mathfrak{aut}(\m,B,\RR))=\h$ is abelian. Otherwise it
is easy to see, that for $\dim(E_i)=2$, $\mathfrak{o}(E_i,B|_{E_i \times E_i}) \simeq \R$ does not commute with $\h$,
but is contained in the radical.
\end{proof}

\pagebreak

\section[Classification of Solvable Symmetric Triples  in Signature \mb{$(2,n-2)$}]
{\begin{minipage}[t]{14.5cm} 
Classification and Normal forms of Solvable Symmetric Triples in Signature \mb{$(2,n-2)$} \end{minipage} }\label{ss:Classification(2,n-2)}

An indecomposable solvable symmetric triple of signature $(2,n-2)$ has either a one- or two-dimensional hence maximal center.
The treatment of triples with maximal center is based on the results of section \ref{subsec:MaxCenter}, whereas 
a triple with 1-dimensional center requires in the case $[\h,\z^{\perp}] \not \subset \z$ already a 2-step iterated
adapted decomposition as discussed in section \ref{ss:IteratedAdapDecomp}. 
Finally the case $[\h,\z^{\perp}] \subset \z$ is a slight generalization of the Lorentzian case. 

The aim of this chapter is to develop for each case sketched above a normal form depending on certain parameters.
Next we determine the isomorphism 
classes among the normal forms by equivalence relations on the parameter spaces and discuss indecomposability.
Thus a classification is obtained.
Among the triples with maximal center the nilpotent ones are characterized by the property $W=W_{nil}$ (lemma \ref{lem:NilpotentMaxCenter}).
On the other hand it has been shown in \cite[lemma 2]{CahenP70} that a nilpotent symmetric triple of dimension greater
than two has at least a two dimensional center. This is a general fact on nilpotent metric Lie algebras (prop. \ref{prop:Bordemann}).
Thus for signature $(2,n-2)$ nilpotent triples occur only among those with maximal center. 
It will turn out that their are exactly 4 such triples up to isomorphism.
Their normal forms are extracted on page \pageref{NorNil}. 

\begin{prop}\cite[prop. 2.2]{Bordemann97}\label{prop:Bordemann}
Let $\g$ be a nilpotent Lie algebra of dimension $\geq 2$ which admits an $\ad$-invariant
non-degenerate symmetric bilinear form. Then $\dim \big( \z(\g) \big) \ge 2$.
\end{prop}

\subsection{Triples with maximal center}

By lemma \ref{lem:afSeperation} one has to distinguish the case of 
non-trivial $W_{nil}$ from the least nilpotent one.

\vspace{2em}
{\large\mb{(I) \hspace{.3cm} The least nilpotent case $W_{nil}=\{0\}$ }}
\addcontentsline{toc}{\protect subsubsection}{(I) $\;W_{nil}=\{0\}$  (least nilpotent case)}
\vspace{1em}

This case is essentially the declination of theorem \ref{thm:Main} for signature $(2,n-2)$.
Every such triple is isomorphic to a triple $\tau_{2,n-2}(\RR_b,\epsilon,\Lambda)$ (cf. \ref{dfn:TripleLeastNilpotent}).

\enlargethispage{2em}
\vspace{1em}
{\mb{(Ia) \hspace{.3cm} $W_{nil}=\{0\}$, $\RR_b \neq 0$}} \nopagebreak \\
\fbox{
\parbox{15cm}{
\begin{empt}\label{triple:tau(epsilonY)}
{\mb{Normal form $\tau_n(\epsilon_{Y},r,\lambda)$}}
Consider parameters $\epsilon_{Y} \in \{\pm 1\}$, $r \in \{0,1,\ldots,n-4\}$ and 
$\lambda_i \in \R^{n-4}$, $(i=1,2)$ with $(\lambda_1^{\alpha},\lambda_2^{\alpha}) \ne 0, \; \forall \alpha$.
Then with respect to a basis
\[ 
\beta_{\h} =\{y,W_1^*,\ldots,W_{n-4}^*\} \qquad \beta_{\m} =\{Z_1,Z_2,W_1,\ldots,W_{n-4},Z_1^*,Z_2^*\}
\]
the metric and Lie algebra structure of $\tau_n(\epsilon_{Y},r,\lambda)=(\h \oplus \m,B)$ are as follows:

\[
\ba{rcl}
B|_{\h \times \h}&=&
\left(
\begin{smallmatrix}
       \epsilon_Y &     \\
                  &  I_{r,n-r-4} 
\end{smallmatrix}
\right) \\{}\\
B|_{\m \times \m}&=&
\left(
\begin{smallmatrix}
           &           & I_2   \\
           &I_{n-4}    &       \\
 I_2       &           &       
\end{smallmatrix}
\right) \\{}\\
\epsilon_{\alpha}&=&
\left\{ \begin{smallmatrix}
-1& 1 \le \alpha \le r\\
+1& \quad r < \alpha \le n-4
\end{smallmatrix} \right.
\ea 
\hspace{2cm}
\ba{rcl}
{[Z_1^*,Z_2^*]} &=& \phantom{-}y\\
{[Z_i^*,W_{\alpha}]} &=& \phantom{-}\lambda_i^{\alpha} \cdot W_{\alpha}^*\\
{[y,Z_1^*]} &=& \phantom{-}\epsilon_Y \cdot Z_2\\
{[y,Z_2^*]} &=& -\epsilon_Y \cdot Z_1\\
{[W_{\alpha}^*,Z_i^*]} &=& \phantom{-}\epsilon_{\alpha}\lambda_i^{\alpha} \cdot W_{\alpha}\\
{[W_{\alpha}^*,W_{\beta}]} &=&
   -\delta_{\alpha\beta}\, \epsilon_{\alpha}(\lambda_1^{\alpha} \cdot Z_1 +\lambda_2^{\alpha} \cdot Z_2)
\ea 
\]
\end{empt} 
}
}

\begin{proof}
Consider a triple $\tau_{2,n-2}(\RR_b,\epsilon,\Lambda)$ together with a fixed basis $\{Z^*_1,Z^*_2\}$ of $\z^*$
and the dual one $\{Z_1,Z_2\}$ in $\z$.
Since  $\Re(\z^*)$ is 1-dimensional $\RR_b \neq 0$ implies $\ker(\RR_b)=\{0\}$, 
hence $Y = \Lambda^2(\z^*)= span\{y\}$ and $b=b_{1212}=B(y,y)$ where $y:=Y_{12}=[Z_1^*,Z_2^*]$. 
A scaling of $Z^*_i$ by $|b|^{-1/4}$ provides $|b|=1$
whereas the signature $\epsilon_Y:=b$ of $B$ restricted to $Y=[\z^*,\z^*]$ is an invariance
since the adapted decomposition $\z \oplus W \oplus \z^*$ is unique up to an isometry.
This factor is missing in \cite{CahenP70}. Analog $(r,n-r-4):=\Sig(W^*)$ is an invariance and 
by a permutation of the basis in $W$
the $\epsilon_{\alpha}$ can be arranged by their sign.
Finally set $\lambda_i^{\alpha}:=B(\Lambda_{\alpha},Z_i^*)$.
\end{proof}

The remaining freedom of the normal form is due to theorem \ref{thm:Main}(iii) first a permutation of the basis $\{W_{\alpha}\}$
respecting the partition $\{1,\ldots,r\}\cup\{r+1,\ldots,n-r-4\}$, second multiplication of $(\lambda_1^{\alpha},\lambda_2^{\alpha})$
by a sign  and third a change of the basis in $\z^*$ respecting $|b|=1$, i.e.
a transformation $S \in Gl(\z^*) \cap Aut(\z^*,\RR_b)$. 
The equivalence relation on the parameters $\lambda_i^{\alpha}$ caused by the first two points
is described by orbits of the following group acting on $\R^{n-4}$:

\mb{$\cal{S}_r \times \cal{S}_{s-r} \times \{\ \pm 1\}^s$ :} Consider $\R^s$ with its canonical basis $\{e_i\}_{i=1,\ldots,s}$. 
The subgroup $\cal{S}_r \times \cal{S}_{s-r} \subset S_s$ of the symmetric group respecting the partition 
$\{1,\ldots,r\}\cup\{r+1,\ldots,n-r-4\}$ acts on $\R^s$ by
\[
\cal{S}_r \times \cal{S}_{s-r} \;\rightarrow \; GL(s,\R), \qquad \pi \; \mapsto \;\{e_i \mapsto e_{\pi(i)}\} 
\]
and the group $\{\pm 1\}^{s}$ acts by reflections on the coordinate hypersurfaces ${e_i}^{\perp}$, that is
\[
\{\pm 1\}^{s} \;\rightarrow \; GL(s,\R), \qquad  (\sigma_1,\ldots,\sigma_s)  \; \mapsto \;\{e_i \mapsto \sigma_i \cdot e_i \} 
\]

For the third point we observe 
$\wt{b}=B \big( \, [S Z_1^*,S Z_2^*],[S Z_1^*,S Z_2^*] \, \big)=det(S)^2 \cdot b$, hence
\[
Gl(\z^*) \cap Aut(\z^*,\RR_b)=Sl(\z^*)\times \{\pm id_{\z^*}\}
\]
If we set $S =\left(\begin{smallmatrix}s_{11}&s_{12}\\s_{21}&s_{22}\end{smallmatrix}\right)$ w.r.t. the basis $\{Z_1^*,Z^*_2\}$ 
the transformation of the $\lambda_i$ is given by
\begin{equation}\label{eq:Lambda}
\rca{1.5}{\ba{c} \wt{\lambda}_i^{\alpha} = B \big( \Lambda_{\alpha},S(Z_i^*) \big)= 
B(\Lambda_{\alpha},s_{1i}Z_1^* + s_{2i}Z_2^*)=s_{1i}\lambda_1^{\alpha} + s_{2i}\lambda_2^{\alpha} \\
 \Leftrightarrow \quad (\wt{\lambda}_1^t,\wt{\lambda}_2^t)=({\lambda}_1^t,{\lambda}_2^t)(S)
\ea}
\end{equation}
\vspace{1em}

{\bf \mb{ $\lambda_1,\lambda_2$} linear independent:}
Every orbit of $(\lambda_1,\lambda_2)$ under (\ref{eq:Lambda}) can be uniquely represented by its span as element in the 
Grassmanian $\mathbb{G}_{2,n-4}$ and the volume of its span in $\R^+$. Thus we get the following bijection
\begin{equation}\label{eq:Vlambda}
(\lambda_1,\lambda_2) \Big( SL(2,\R) \times \{\pm1 \}\Big) \quad \mapsto \quad  (P_{\lambda},V_{\lambda}) :=
\Big( span\{\lambda_1,\lambda_2\},\sum_{\alpha,\beta}(\lambda_1^{\alpha} \lambda_2^{\beta}-\lambda_2^{\alpha} \lambda_1^{\beta})^2 \Big)
\end{equation}
Finally $(\lambda_1^{\alpha},\lambda_2^{\alpha}) \ne 0$ is a necessary and sufficient condition for indecomposability by \ref{thm:Main}(ii),
since in this case $[\h,\m]=\z \oplus W$ and $1=\dim(Y) > p(p-2)/2=0$. 
Under the correspondence (\ref{eq:Vlambda}) $(\lambda_1^{\alpha},\lambda_2^{\alpha}) \ne 0$ 
translates into $P_{\lambda} \not\subset \bigcup e_{\alpha}^{\perp}$. Thus we proved:

\begin{prop}\label{prop:IaRank2}
The isomorphism classes of indecomposable symmetric triples $\tau_n(\epsilon_{Y},r,\lambda)$,
$(\lambda_i)_{i=1,2}$ of rank $2$, are in one-to-one
correspondence with 
{\small
\[
\Big(\epsilon_Y,r,V_{\lambda},\big[P_{\lambda}\big] \Big) \; \in 
\; \{\pm 1\}\; \times \;\{0,\ldots,n-4\} \;\times \R^+ \times \; 
{}^{\ds \Menge{P \in \mathbb{G}_{2,n-4}}{ P \not\subset \bigcup e_{\alpha}^\perp}}
\Big/
{}_{ \ds \mathcal{S}_r \times \mathcal{S}_{n-r-4} \times \{\pm 1\}^{n-4} }
\]
}
where $\lambda_1,\lambda_2$ are determined
by $V_{\lambda}$ and $P_{\lambda}$ according to (\ref{eq:Vlambda}). In particular $n>5$.
\end{prop}

\vspace{1em}
{\bf \mb{ $\lambda_1,\lambda_2$} linear dependent:}
Let $d_1 \lambda_1 = d_2 \lambda_2,\; (d_1,d_2)\ne 0$. In this case one can achieve $\left\|\lambda_1 \right\|=1$ and $\lambda_2=0$
by application of $S=\frac{1}{(d_1)^2+(d_2)^2} 
\left(\begin{smallmatrix} d_1&-d_2\\d_2&d_1\end{smallmatrix}\right) \in Sl(\z^*) \times \{\pm id_{\z^*}\}$ 
followed by  a transformation of the form $S=\left(\begin{smallmatrix} \delta & 0 \\ 0 & \delta^{-1}\end{smallmatrix}\right)$. 
Analog to the preceeding case the triple is indecomposable if and only if $\lambda_1 \not\subset \bigcup e_{\alpha}^{\perp}$.
Thus we obtain:

\begin{prop}\label{prop:IaRank1} 
The isomorphism classes of indecomposable symmetric triples $\tau_n(\epsilon_{Y},r,\lambda)$,
$(\lambda_i)_{i=1,2}$ of rank $1$, are in one-to-one
correspondence with 
{\small
\[
\Big(\epsilon_Y,r,[\lambda_1] \Big) \; 
\in \; \{\pm 1\}\; \times \;\{0,\ldots,n-4\}  \times \; 
{}^ {\ds \Menge{\lambda \in S^{n-5}}{\lambda \not\subset \bigcup e_{\alpha}^{\perp}}} 
\Big/
{}_{\ds \mathcal{S}_r \times \mathcal{S}_{n-r-4} \times \{\pm 1\}^{n-4} }
\]
}
with $\lambda_2=0$. In particular $n>4$.
\end{prop}

\vspace{1em}
{\mb{(Ib) \hspace{.3cm} $W_{nil}=\{0\}$, $\RR_b = 0$}}\\
\fbox{
\parbox{15cm}{
\begin{empt}\label{triple:tau(epsilon)}
\mb{Normal form $\tau_n(r,\lambda)$}
Consider parameters $r \in \{0,1,\ldots,n-4\}$ and
$\lambda_i \in \R^{n-4}$, $(i=1,2)$ with $(\lambda_1^{\alpha},\lambda_2^{\alpha}) \ne 0, \; \forall \alpha$.
Then with respect to a basis
\[ 
\beta_{\h} =\{W_1^*,\ldots,W_{n-4}^*\} \qquad \beta_{\m} =\{Z_1,Z_2,W_1,\ldots,W_{n-4},Z_1^*,Z_2^*\}
\]
the metric and Lie algebra structure of $\tau_n(r,\lambda)=(\h \oplus \m,B)$ is as follows ($\epsilon_{\alpha}$ as in (Ia)):
\[
\ba{rcl}
B|_{\h \times \h}&=& I_{r,n-r-4} \\ {} \\
B|_{\m \times \m}&=&
\left(
\begin{smallmatrix}
           &           & I_2   \\
           &I_{n-4}    &       \\
 I_2       &           &       
\end{smallmatrix}
\right)
\ea
\hspace{2cm}
\ba{rcl}
{[Z_i^*,W_{\alpha}]} &=& \phantom{-}\lambda_i^{\alpha} \cdot W_{\alpha}^*\\
{[W_{\alpha}^*,Z_i^*]} &=& \phantom{-}\epsilon_{\alpha}\lambda_i^{\alpha} \cdot W_{\alpha}\\
{[W_{\alpha}^*,W_{\beta}]} &=&
   -\delta_{\alpha\beta} \, \epsilon_{\alpha}(\lambda_1^{\alpha} \cdot Z_1 +\lambda_2^{\alpha} \cdot Z_2)
\ea
\]
\end{empt} 
}
}

In comparison with (Ia) the condition $[\h,\m]=\z \oplus W$ is satisfied if and only the rank of $(\lambda_i^{\alpha})$ is $2$,
hence $n \ge 6$. Due to theorem \ref{thm:Main}(ii) $\tau_n(r,\lambda)$ still may decompose 
into two solvable Lorentzian triples of dimension $\geq 3$.
Indeed, let $\tau_{n_i}(f_i) \in \mathcal{N}_{Lor}^{n_i}$, $n_i \ge 3$ ($i=1,2$). Then it is easily verified, that
\[
\tau_{n_1}(f_1) \oplus \tau_{n_2}(f_2) \simeq \tau_{n_1+n_2}(r,\lambda)
\]
where $r:=\# \{ f_i^{\alpha} < 0 \}$ and $(\lambda_1^t,\lambda_2^t) := \left( \begin{smallmatrix}
\sqrt{|f_1|^t}&0\\
0&\sqrt{|f_2|^t}
\end{smallmatrix} 
\right)$ (up to a permutation of the rows).   
Since there is no restriction on a transformation in $Gl(\z^*)$ we get:

\begin{prop}\label{prop:Ib} 
The isomorphism classes of indecomposable symmetric triples $\tau_n(r,\lambda)$ are in one-to-one
correspondence with 
{\small
\[
\Big(r,\big[P_{\lambda}\big] \Big) \; \in \;
\{0,\ldots,n-4\}  \times \; {}^ {\ds 
\Menge{ P \in \mathbb{G}_{2,n-4} }{ \ba{ll}  
\text{if}& P \cap E_{\alpha_1, \ldots \alpha_s}\neq \{0\} \\
\Rightarrow & P \cap {E_{\alpha_1, \ldots \alpha_s}}^{\perp} = \{0\} \ea}
} \Big/
{}_{\ds \mathcal{S}_r \times 
\mathcal{S}_{n-r-4} \times \{\pm 1\}^{n-4} }
\]
}
where $E_{\alpha_1, \ldots \alpha_s}:= span\{e_{\alpha_1},\ldots,e_{\alpha_s}\}$ $(0<s<n-4)$ and
$\lambda_1,\lambda_2$ are determined by $span\{\lambda_1,\lambda_2\} = P_{\lambda}$.
In particular $n > 6$.
\end{prop}

\pagebreak


\subsubsection*{(II) \mb{ The case $W_{nil} \neq \{0\}$}}
\addcontentsline{toc}{\protect subsubsection}{(II) $W_{nil} \neq \{0\}$ }

As pointed out before there is not much known about these triples in general signature (cf. p.\pageref{case:f=0}). 
The essential fact which enables a classification in signature $(2,n-2)$ is that the dimension of $W_{nil}$
is restricted due to (\ref{eq:Dimf=0}) by $\dim(W_{nil}) \leq 2$. The nilpotent triples among them, i.e.
those with $W=W_{nil}$ are cited in \cite{CahenP70}.
Due to the lack of a general theory (comparable to theorem \ref{thm:Main} for least nilpotent
triples) we will study in each case the action of $\Gamma$ on the structure coefficients by hand.
Indecomposability is easily verified in this case:

\begin{lem} 
A solvable symmetric triple $\tau$ with a $1$-dimensional center is indecomposable.
\end{lem}

\begin{proof}
Assume $\tau$ decomposes into indecomposable triples $\tau=\tau_1 \oplus \ldots, \oplus \tau_s$. 
Then each $\tau _i$ is again solvable i.e. $\dim(\z(\tau_i)) \geq 1$. Hence $s=1$ due to (\ref{eq:CenterSum}).
\end{proof}

\vspace{1em}
\mb{(IIa) \hspace{.3cm} $\dim(W_{nil})=1$ }\\
\fbox{
\parbox{15cm}{
\begin{empt}{\mb{ Normal form $\tau_n \big(r,\lambda,1\big)$}} \label{triple:tau(If=1)}
Consider parameters $r \in \{1,\ldots,n-5\}$ and $\lambda_i \in \R^{n-5}$, $(i=1,2)$ with 
$(\lambda_1^{\alpha},\lambda_2^{\alpha}) \ne 0,\; \forall \alpha$. 
Then with respect to a basis
\[
\beta_{\h} =\{y,W_1^*\ldots,W_{n-5}^*, w^* \} \quad \beta_{\m} =\{Z_1,Z_2,W_1,\ldots,W_{n-5},w,Z_1^*,Z_2^*\}
\]
the metric and Lie algebra structure are given by
\[
\!\!\!\!\ba{rcl}
B|_{\h \times \h} &=&
\left(
\begin{smallmatrix}
           &                & 1 \\
           &  I_{r,n-r-5}   &   \\
1          
\end{smallmatrix}
\right) \\{} \\ 
B|_{\m \times \m}&=&
\left(
\begin{smallmatrix}
           &           & I_2   \\
           &I_{n-4}    &       \\
 I_2       &           &       
\end{smallmatrix}
\right)\\{}\\
\epsilon_{\alpha}&=&
\left\{ \begin{smallmatrix}
-1& 1 \le \alpha \le r\\
+1& \quad r < \alpha \le n-5
\end{smallmatrix} \right.\ea
\ba{rcl}
{[Z_1^*,Z_2^*]} &=& y\\
{[Z_1^*,w]} &=& w^*\\
{[Z_i^*,W_{\alpha}]} &=& \lambda_i^{\alpha} \cdot W_{\alpha}^* \\
\ea
\ba{rcl}
{[y,Z_1^*]} &=& \phantom{-}w\\
{[y,w]} &=& -Z_1\\
{[w^*,Z_1^*]} &=& \phantom{-}Z_2\\
{[w^*,Z_2^*]} &=& -Z_1\\ 
{[W_{\alpha}^*,Z_i^*]} &=& \phantom{-}\epsilon_{\alpha}\lambda_i^{\alpha} \cdot W_{\alpha} \\
{[W_{\alpha}^*,W_{\beta}]} &=&
   -\delta_{\alpha\beta} \epsilon_{\alpha}(\lambda_1^{\alpha} \cdot Z_1 +\lambda_2^{\alpha} \cdot Z_2) 
\ea
\]
\end{empt} 
}
}

\begin{proof}
Any triple of this type can be represented in the form of prop. \ref{prop:1.NFMaxCenter}. Furthermore we may choose an adapted decomposition
$\z \oplus W \oplus \z^*$ and  a basis of $W$ with $W_{reg}=span\{W_1,\ldots,W_{n-5}\}$ and $W_{nil}=\R \cdot W_{n-4}$ s.t. 
$a_{ijk\alpha}=0,\; \forall \alpha \leq n-5$ and $f_{ij\alpha\beta}=\delta_{\alpha\beta} \epsilon_{\alpha} \lambda_i^{\alpha}\lambda_j^{\alpha}$
(lemma \ref{lem:Fij}, \ref{lem:afSeperation}). We simplify the notation to $a_k:=a_{12k\alpha}$ since $a_{12k\alpha}=-a_{21k\alpha}$
and $a_{11k\alpha}=a_{22k\alpha}=0$. The 1.Bianchi identity $\sigma_{ijk\alpha}(a_{ijk\alpha})=0$ and the relation  (\ref{eq:a})
are trivially satisfied. Thus we get the following parametrisation
\[
\begin{array}{rcl}
{[Z_1^*,Z_2^*]} &=& y\\
{[Z_i^*,W_{\alpha}]} &=& X_{i\alpha}\\
\\
{[y,Z_1^*]} &=&  \phantom{-}a_1 \cdot W_{n-4} + b \cdot Z_2\\
{[y,Z_2^*]} &=& \phantom{-} a_2 \cdot W_{n-4} - b \cdot Z_1\\
{[y,W_{n-4}]} &=&  -  a_1 \cdot Z_1 - a_2 \cdot Z_2 
\end{array}
\quad
\begin{array}{rcl}
{[X_{i(n-4)},Z_1^*]} &=& \phantom{-}a_i \cdot Z_2\\
{[X_{i(n-4)},Z_2^*]} &=& -a_i \cdot Z_1\\
\\
(\alpha \leq n-5): &&\\
{[X_{i\alpha},Z_k^*]} &=& \phantom{-}\epsilon_{\alpha}\lambda_i^{\alpha} \lambda_k^{\alpha} \cdot W_{\alpha} \\
{[X_{i\alpha},W_{\beta}]} &=& - \epsilon_{\alpha}\lambda_i^{\alpha}\delta_{\alpha\beta} 
  \big( \lambda_1^{\alpha} \cdot Z_1 +  \lambda_2^{\alpha} \cdot Z_2\big) \\
\end{array}
\]
First lets determine the linear dependency among the generators of $\h$:  
Since $a,\lambda^2,\ldots,\lambda^{n-4} \ne 0$ one can find $W_{\alpha}^* \in \h$ 
s.t. (cf. page \pageref{eq:Walpha*}):
\[
X_{i\alpha}=
\begin{cases}
\lambda_i^{\alpha} \cdot W_{\alpha}^* & \text{ for } \alpha \leq n-5 \\
a_i \cdot W_{n-4}^* &  \text{ for } \alpha = n-4
\end{cases}
\]
Thus $\{y,W_1^*,W_2^*,\ldots,W_{n-4}^*\}$ is a basis of $\h$ (consider its action on $\z^*$).
The coefficients $a_k,\pm \lambda^{\alpha}$ are invariant under $\Gamma_{decomp.}^a$ (\ref{lem:Fij}, \ref{lem:afSeperation})
whereas - in contrast to the least nilpotent case - $b$ may change in the following way as a slight variation of (\ref{eq:bTilde}) shows:
\[
\wt{b}=B(\wt{y},\wt{y}) = 
 b + 2 \big( -(L)_{2(n-4)} a_1 + (L)_{1(n-4)} a_2 \big)
\]
where the only nontrivial entries of $L \in \Gamma_{decomp.}$ modulo $\Gamma_{decomp.}^a$ are $(L_{i(n-4)})_{i=1,2}$.
Hence we can achieve $b=0$ by such a transformation. 
Up to now this triple looks very similar to the case (Ia) interchanging
the role of $y$ with $W_{n-4}^*$.
The main difference occurs in the transformation behavior of $(a_1,a_2)$ under $P \in Gl(\z^*)$ (cf. (\ref{eq:Lambda}) ):
\begin{equation}\label{eq:det(P)}
(\wt{a}_1,\wt{a}_2)=(a_1,a_2)(P)\cdot det(P)
\end{equation}
One can find easily $P \in Gl(\z^*)$ s.t. $(\wt{a}_1,\wt{a}_2)=(1,0)$.
Then the triples receives the form $\tau_n(r,\lambda,1)$ with the substitution $w:=W_{n-4}$ and $w^*:=W^*_{n-4}$.
\end{proof}

In order to determine the isomorphism classes of the triples $\tau_n(r,\lambda,1)$ observe that w.r.t the action of $Gl(\z^*)$ on $(a_1,a_2)$
given by (\ref{eq:det(P)}) we have 
$stab(a)= \big\{ \left(\begin{smallmatrix}
1/\alpha &0 \\ \beta & \alpha^2  \end{smallmatrix} \right) \;|\; \alpha,\beta  \in \R \big\}$ if we require $(a_1,a_2)=(1,0)$. 
Two cases has to be distinguished: 

\mb{$\lambda_1,\lambda_2$ linear dependent:} Then one can change the basis of $\z^*$ s.t. $\lambda_i=0$, $|| \lambda_j || =1$, $\{i,j\}=\{1,2\}$.

\mb{$\lambda_1,\lambda_2$ linear independent:} Then one obtains $\lambda_1 \perp \lambda_2$ and $|| \lambda_1 || = || \lambda_2 ||$ 
by suitable $\beta$ and $\alpha$. The remaining freedom is a sign of $\lambda_1$. Thus every orbit of $(\lambda_1,\lambda_2)$ under $stab(a)$
can be represented uniquely by a unit-tangent vector of $\R \mathbb{P}^{n-6} $ and positive number through
\[
\big( [l_1,l_2],\delta\big) \in S(T\R\mathbb{P}^{n-6}) \times \R^+ \quad \mapsto \quad (\delta \lambda_1,\delta \lambda_2)
\]
Hence:
 
\begin{prop}\label{prop:IIa}
The isomorphism classes of symmetric triples of the type $\tau_n \big( r,\lambda,1 \big)$,
are in one-to-one correspondence with 
\[
\Big(r,k,[l] \Big) 
\; \in \; \{0,\ldots,n-5\} \; \times \; \{1,2\} \; \times \; 
{}^ {\ds \Menge{\lambda \in S^{n-6}}{\lambda \not\subset \bigcup e_{\alpha}^{\perp}}} 
\Big/
{}_{\ds \mathcal{S}_r \times \mathcal{S}_{n-r-5} \times \{\pm 1\}^{n-5} }
\]
with $\lambda_i= \delta_{ik} \cdot l$ and 
\begin{multline*}
\Big(r,\delta,[l_1,l_2] \Big) \; \in \\
 \{0,\ldots,n-5\} \; \times \; \R^+ \; \times \; 
{}^ { \ds \Menge{(l_1,l_2) \in S(T\R\mathbb{P}^{n-6})}{span\{l_1,l_2\} \not \subset \bigcup e_{\alpha^{\perp}}} }  
\Big/_{\ds \mathcal{S}_r \times \mathcal{S}_{n-r-5} \times \{\pm 1\}^{n-5} }  
\end{multline*}
with $(\lambda_1,\lambda_2)=(\delta l_1 , \delta l_2 )$.
The two case correspond to $(\lambda_i^{\alpha})$ being of rank $1$ resp. rank $2$.
\end{prop}

\pagebreak

\mb{(IIb) \hspace{.3cm} $\dim(W_{nil})=2$ }\\
\fbox{
\parbox{15cm}{
\begin{empt}{\mb{ Normal form $\tau_n \big(r,\lambda,2 \big)$}} \label{triple:tau(If=2)}
Consider parameters $r \in \{1,\ldots,n-6\}$ and $\lambda_i \in \R^{n-6}$, $(i=1,2)$ with 
$(\lambda_1^{\alpha},\lambda_2^{\alpha}) \ne 0,\; \forall \alpha$. 
Then with respect to a basis
\[
\beta_{\h} =\{y,W_1^*\ldots,W_{n-5}^*, w^* \} \quad \beta_{\m} =\{Z_1,Z_2,W_1,\ldots,W_{n-6},w_1,w_2,Z_1^*,Z_2^*\}
\]
the metric and Lie algebra structure are given by
\[
\!\!\!\!\ba{rcl}
B|_{\h \times \h} &=&
\left(
\begin{smallmatrix}
           &                & 1 \\
           &  I_{r,n-r-5}   &   \\
1          
\end{smallmatrix}
\right) \\{} \\ 
B|_{\m \times \m}&=&
\left(
\begin{smallmatrix}
           &           & I_2   \\
           &I_{n-4}    &       \\
 I_2       &           &       
\end{smallmatrix}
\right)\\{}\\
\epsilon_{\alpha}&=&
\left\{ \begin{smallmatrix}
-1& 1 \le \alpha \le r\\
+1& \quad r < \alpha \le n-5
\end{smallmatrix} \right.\ea
\ba{rcl}
{[Z_1^*,Z_2^*]} &=& y\\
{[Z_i^*,w_j]} &=& \delta_{ij} w_i \\
{[Z_i^*,W_{\alpha}]} &=& \lambda_i^{\alpha} \cdot W_{\alpha}^* \\
\ea
\ba{rcl}
{[y,Z_i^*]} &=& \phantom{-}w_i\\
{[y,w_i]} &=& -Z_i\\
{[w^*,Z_1^*]} &=& \phantom{-}Z_2\\
{[w^*,Z_2^*]} &=& -Z_1\\ 
{[W_{\alpha}^*,Z_i^*]} &=& \phantom{-}\epsilon_{\alpha}\lambda_i^{\alpha} \cdot W_{\alpha} \\
{[W_{\alpha}^*,W_{\beta}]} &=&
   -\delta_{\alpha\beta} \epsilon_{\alpha}(\lambda_1^{\alpha} \cdot Z_1 +\lambda_2^{\alpha} \cdot Z_2) 
\ea
\]
\end{empt} 
}
}

\begin{proof}
A similar reasoning as in the case (IIa) leads to a parametrisation of these triples by
\[
\begin{array}{rcl}
{[Z_1^*,Z_2^*]} &=& y\\
{[Z_i^*,W_{\alpha}]} &=& X_{i\alpha}\\
\\
{[X_{i(n-6+j)},Z_1^*]} &=& \phantom{-}a_i^j \cdot Z_2\\
{[X_{i(n-6+j)},Z_2^*]} &=& -a_i^j \cdot Z_1\\
\end{array}
\quad
\begin{array}{rcl}
{[y,Z_1^*]} &=&  \phantom{-}a_1^1 \cdot W_{n-5} + a_1^2 \cdot W_{n-4} + b \cdot Z_2\\
{[y,Z_2^*]} &=& \phantom{-} a_2^1 \cdot W_{n-5} + a_2^2 \cdot W_{n-4} - b \cdot Z_1\\
{[y,W_{n-6+i}]} &=&  -  a_1^i \cdot Z_1 -a_2^i \cdot Z_2 
\\
(\alpha \leq n-6): &&\\
{[X_{i\alpha},Z_k^*]} &=& \phantom{-}\epsilon_{\alpha}\lambda_i^{\alpha} \lambda_k^{\alpha} \cdot W_{\alpha} \\
{[X_{i\alpha},W_{\beta}]} &=& - \epsilon_{\alpha}\lambda_i^{\alpha}\delta_{\alpha\beta} 
  \big( \lambda_1^{\alpha} \cdot Z_1 +  \lambda_2^{\alpha} \cdot Z_2\big) \\
\end{array}
\]
Analog to the preceeding case we can find a basis $\{y,W_1^*,\ldots,W_{n-6}^*,w^*\}$ of $\h$ by
\[
X_{i\alpha}=
\begin{cases}
\lambda_i^{\alpha} \cdot W_{\alpha}^* & \text{ for } \alpha \leq n-6 \\
a_i^{\alpha-n-6} \cdot w^* & \text{ for } \alpha =n-4,n-5\\
\end{cases}
\]
The condition $[\h,\m]=\z \oplus W$ is satisfied if 
$
(A)=
\left( \begin{smallmatrix}
a_1^1&a_2^1\\a_1^2&a_2^2 
\end{smallmatrix} \right) $ is invertible.
The same arguments as in (IIa) show that we can achieve $b=0$ by the choice of a suitable decomposition.
Finally (\ref{eq:det(P)}) takes in this situation the particular nice form
\[
(\wt{A})=(A)\cdot (P) \cdot det(P),\qquad P \in Gl(\z^*),
\]
We set $(P)=\alpha \cdot (A)^{-1}$ and compute $(A)\cdot (P) \cdot det(P)= \alpha^3/det(A) \cdot  I_2$. Then
$\alpha^3=\big(det(A)\big)^{1/3}$ yields $(A)=I_2$. Finally set $w_i:=W_{n-6+i}$ $(i=1,2)$.
\end{proof}

\begin{prop}\label{prop:IIb}
The isomorphism classes of symmetric triples of the type
$\tau_n \big( r,\lambda,2 \big)$,
are in one-to-one correspondence with 
\begin{multline*}
\Big(r,[\lambda_1,\lambda_2] \Big) \; \in  \\
\{0,\ldots,n-6\} \; \times \;
{}^{ \ds \Menge{ (l_1,l_2) \in \R^{n-6} \times \R^{n-6})}{span\{l_1,l_2\}  \not \subset \bigcup e_{\alpha}^{\perp}}} 
\Big/_{\ds \mathcal{S}_r \times \mathcal{S}_{n-r-6} \times \{\pm 1\}^{n-6} }
\end{multline*}
\end{prop}

\pagebreak

{\bf The nilpotent triples}\label{NorNil}
\addcontentsline{toc}{\protect subsubsection}{\hspace{5mm} The nilpotent triples}
\vspace{1em}

According to lemma \ref{lem:afSeperation} a nilpotent triple satisfies $W=W_{nil}$.
In the case (I) this implies $W=\{0\}$. Thus the triple is of signature $(2,2)$. Moreover such a triple is indecomposable 
if and only if $[\h,\m] = \z$ i.e. $\RR_b \ne 0$, hence we are lead into the case (Ia).
In the cases (IIa) and (IIb) the signature is $(2,3)$ and $(2,4)$.

{\mb{(Ia') \hspace{.3cm} $W=\{0\}$}:\\
\fbox{
\parbox{15cm}{
\begin{empt}\label{triple:tau(2,2)}
{\mb{Normal form $\tau_{2,2}(\epsilon_Y)$}}
Consider $\epsilon_{Y} \in \{\pm 1\}$. Then with respect to a basis
\[ 
\beta_{\h} =\{y\} \qquad \beta_{\m} =\{Z_1,Z_2,Z_1^*,Z_2^*\}
\]
the metric and Lie algebra structure are given by:

\[
\ba{rcl}
B|_{\h \times \h}&=& \epsilon_Y \\
{}\\
B|_{\m \times \m}&=&
\left(
\begin{smallmatrix}
         & I_2   \\
    I_2         
\end{smallmatrix}
\right)
\ea 
\hspace{2cm}
\ba{rcl}
{[Z_1^*,Z_2^*]} &=& \phantom{-}y\\
{[y,Z_1^*]} &=& \phantom{-}\epsilon_Y \cdot Z_2\\
{[y,Z_2^*]} &=& -\epsilon_Y \cdot Z_1\\
\ea 
\]
\end{empt} 
}
}

\vspace{2em}
\mb{(IIa') \hspace{.3cm} $\dim(W)=1$ } (cf.\cite[5.2]{CahenP70}): \\
\fbox{
\parbox{15cm}{
\begin{empt}\mb{ Normal form $\tau_{2,3}$}\label{triple:tau(2,3)}
With respect to a basis
\[
\beta_{\h} =\{y,w^* \} \quad \beta_{\m} =\{Z_1,Z_2,w,Z_1^*,Z_2^*\}
\]
the metric and Lie algebra structure are given by
\[
\ba{rcl}
B|_{\h \times \h} &=&
\left(
\begin{smallmatrix}
           &    1 \\
1          
\end{smallmatrix}
\right)\\
{}\\
B|_{\m \times \m}&=&
\left(
\begin{smallmatrix}
           &           & I_2   \\
           &1    &       \\
 I_2       &           &       
\end{smallmatrix}
\right)
\ea
\ba{rcl}
{[Z_1^*,Z_2^*]} &=& y\\
{[Z_1^*,w]} &=& w^*\\
\ea
\ba{rcl}
{[y,Z_1^*]} &=& \phantom{-}w\\
{[y,w]} &=& -Z_1\\
{[w^*,Z_1^*]} &=& \phantom{-}Z_2\\
{[w^*,Z_2^*]} &=& -Z_1\\ 
\ea
\]
\end{empt} 
}
}

\vspace{2em}
\mb{(IIb') \hspace{.3cm} $\dim(W)=2$ }(cf.\cite[5.4]{CahenP70}):\\
\fbox{
\parbox{15cm}{
\begin{empt}{\mb{ Normal form $\tau_{2,4}$}} \label{triple:tau(2,4)}
With respect to a basis
\[
\beta_{\h} =\{y,w^* \} \quad \beta_{\m} =\{Z_1,Z_2,w_1,w_2,Z_1^*,Z_2^*\}
\]
the metric and Lie algebra structure are given by
\[
\ba{rcl}
B|_{\h \times \h} &=&
\left(
\begin{smallmatrix}
           &   1 \\
1          
\end{smallmatrix}
\right) \\{} \\ 
B|_{\m \times \m}&=&
\left(
\begin{smallmatrix}
           &           & I_2   \\
           &I_2    &       \\
 I_2       &           &       
\end{smallmatrix}
\right)
\ea
\ba{rcl}
{[Z_1^*,Z_2^*]} &=& y\\
{[Z_i^*,w_j]} &=& \delta_{ij} w_i
\ea
\ba{rcl}
{[y,Z_i^*]} &=& \phantom{-}w_i\\
{[y,w_i]} &=& -Z_i\\
{[w^*,Z_1^*]} &=& \phantom{-}Z_2\\
{[w^*,Z_2^*]} &=& -Z_1\\ 
\ea
\]
\end{empt} 
}
}

\subsection{1-dimensional center}

According to proposition \ref{prop:AdapDecomp}(ii) there are two possibilities for a \emph{complete} adapted decomposition
of $\m$:
\begin{enumerate}
\item[(III)] $\m = \z \oplus W \oplus \z^*$, i.e. $[\h,W] \subset \z$ and $(W,B|_{W \times W})$ 
is a Lorentzian subspace of $\m$.
\item[(IV)] $\m = \z \oplus E \oplus W \oplus U \oplus E^* \oplus \z^*$, i.e. 
$[\h,W_0]_{W_0}=E \oplus W$, $[\h,W] \subset \z \oplus E$ (with $W_0=E \oplus W \oplus U \oplus E^*$),
and $(W,B|_{W \times W})$ is an Euclidian subspace of $\m$.
\end{enumerate}

\subsubsection*{(III)\label{section:(III)} \mb{ Triples with $[\h,\z^{\perp}] \subset \z$  }}
\addcontentsline{toc}{\protect subsubsection}{(III) $[\h,\z^{\perp}] \subset \z$ }

There are no obstructions to discuss these triples in arbitrary signature. 
Restriction to the Lorentzian case yields ones more the classification \ref{triple:tau(f)}.
Lets recapitulate what we know from proposition \ref{prop:AdapDecomp}:

There exists a basis $\beta_{\m}=(Z,W_1,\ldots,W_{n-2},Z^*)$ in $\m$ such that
\[
Mat_{\beta_{\m}}(B)= \left(
\begin{matrix}
&&1\\ &I_{(p,n-p-2)}&\\1&&
\end{matrix} \right)
\qquad
Mat_{\beta_{\m}}(\ad(h)|_{\m})= \left(
\begin{matrix}
0 & M_W^h & 0\\ 
0 & 0     & -I_{(p,n-p-2)} (M_W^h)^t \\
0 & 0     &0
\end{matrix} \right)
\]
Moreover $\z \oplus W$ is abelian and $\h=[\z^*,W]$. Hence $\dim(\h) \le \dim(W)$. 
On the other hand $[\h,\m]=\z \oplus W$, more precisely $[\h,W]=\z$ and $[\h,\z^*]=W$. Hence $\dim(\h) \ge \dim(W)$.
Thus $\dim(\h)=\dim(W)$ and $\beta_{\h}=\{X_1,\ldots,X_{n-2}\}$ with $X_{\alpha}=[Z^*,W_{\alpha}]$
forms a basis of $\h$.
Setting $(M_W^{X_{\alpha}})_{\beta}=-f_{\alpha\beta}$ yields 
$B(X_{\alpha},X_{\beta})=\RR \big( Z^*,W_{\alpha},W_{\beta},Z^* \big) =B \Big( \big[ [Z^*,W_{\alpha}],W_{\beta} \big],Z^* \Big)=-f_{\alpha\beta}$.
In particular $(f_{\alpha\beta})$ is a symmetric non-degenerate matrix.
It remains to discuss the Jacobi-identity:
\begin{itemize}
\item In $\m$ this gives nothing new since $W$ is abelian.
\item $\big[ [\m,\m],\h \big]$ defines the Lie algebra structure of $\h$: 
\[
[X_{\alpha},X_{\beta}]= \big[ [Z^*,W_{\alpha}],X_{\beta} \big]=\big[Z^*,\underbrace{[W_{\alpha},X_{\beta}]}_{\subset \z}\big] - 
\big[ W_{\alpha},\underbrace{[Z^*,X_{\beta}]}_{\subset W} \big]=0
\]
Thus the holonomy is abelian.
In particular $\big[ [\h,\h],\h \big]=\{0\}$.
\item $ \big[ [\m,\h],\h \big]$ is similar to (\ref{eq:Kommutativ}) but trivially satisfied since $\dim(\z)=1$.
\end{itemize}

\vspace{1em}
\begin{dfn}[\mb{$\tau_{p,q}(F,I)$}]\label{dfn:tau(F,I)}
Consider he following class of symmetric triples $\tau_{p,q}(F,I)=(\h \oplus \m,B)$:
\bi
\item $\m=\z \oplus W \oplus \z^*$, $\z=\R\cdot Z$, $\z^*=\R\cdot Z^*$
\item $\h=W^*$ the dual of $W$ with a distinguished isomorphism  $I:W \rightarrow W^*$
\item $B(Z,Z) = B(Z^*,Z^*)=0$, $B(Z^*,Z) = 1$, $\z \oplus \z^* \perp_B  W$ and $\Sig(B|_{W^2})=(p-1,q-1)$
\item $F \in Gl(W)$ a $B|_{W^2}$-selfadjoint, invertible map
\item The non-zero commutators of the Lie algebra structure are 
\[
\begin{array}{rcl}
{[Z^*,w]} &=& \phantom{-}I(w)\\
{[I(w),Z^*]} &=& \phantom{-}F(w)\\
{[I(w),\wt{w}]} &=& -B(F(w),\wt{w}) \cdot Z
\end{array}
\]
\ei
\end{dfn}

\pagebreak
\begin{prop}
\begin{enumerate}
\item[(i)]If $\tau$ is a solvable triple of signature $(p,q)$ with $1$-dimensional center and
$[\h,\z^{\perp}] \subset \z$, then there exists a triple $\tau_{p,q}(F,I)$ which is isomorphic to $\tau$ and conversely. 
\item[(ii)] Two triples $\tau_{p,q}(F,I),\tau_{p,q}(\wt{F},\wt{I})$ 
are isomorphic if and only if there exists a map $P \in O(W,B|_{W^2})$ and $c \in \R \setminus \{0\}$ such that
$\wt{I}=c \cdot (P^{-1})^* \circ I= c \cdot I \circ P^{-1}$ and $F= c^2 \cdot P^{-1} \circ \wt{F} \circ P$.
\end{enumerate}
\end{prop}
\begin{proof}
Let $\tau$ be a triple as assumed in (i). Any $([\h,\m],\m,B)$-adapted decomposition $\m=\z \oplus W \oplus \z^*$ yields
the operators
\[  
I= \ad(Z^*)|_W, \qquad F=pr_W \circ \big( \ad(Z^*) \big)^2 \big|_W
\]
From the preceeding discussion we know, that $I$ is an isomorphism and $F$ is selfadjoint and invertible ($Mat(F)=(f_{\alpha\beta})$). 
On the other hand it is easy to check that each triple $\tau(I,F)$ is a symmetric triple with the required properties.

(ii): As usual in this context we have to study the effect of $\Gamma_{decomp.} \rtimes \Gamma_{basis}^D$ on the triple.
First we show $\Gamma_{decomp.} \subset Aut(\m,B,\RR)$. Indeed, $\Gamma_{decomp.} \simeq Hom(W,\z) \simeq W^*$ since $\dim(\z)=1$, hence
$Skew(\z,\z^*,B)=\{0\}$. Further we observe, that for $L \in Hom(W,\z)$ exists a unique $w_L \in W$ such that
$L(\wt{w})=[I(w_L),\wt{w}]=-B(F(w_L),\wt{w})Z$ since $F$ is invertible. 
Thus the isometry corresponding to $L$ is given by $\rho_L =Exp(\ad(I(w_L))|_{\m})$ and consequently $\Gamma_{decomp.} = Exp(\ad(\h)|_{\m}) \subset Aut(\m,B,\RR)$.
Now we fix a decomposition $D$. Then $\gamma \in \Gamma_{basis}^D$, can be written as
$\gamma=c \cdot id_{\z^*} \oplus P \oplus 1/c \cdot id_{\z} \in Gl(\z^*) \times B(W,B|_{W^2}) \times Gl(\z)$.
\end{proof}

In order to obtain a normal form for such a symmetric triple, it is necessary to determine normal forms of selfadjoint
operators with respect to a non-degenerate, possibly indefinite symmetric bilinear-form on a real vector space. 
This has been done in \cite{Klingenberg54} and reviewed in \cite{Boubel00}. For signature $(2,n-2)$ we get:

\mb{ (III) \hspace{.3cm} $\dim(\z)=1$, $[\h,z^{\perp}]\subset \z$: }\\
\fbox{
\parbox{15cm}{
\begin{empt}\label{triple:tau(Phi)}
\mb{Normal form  $\tau_{n}(g,\Phi,f)$}
Consider parameters $f=(f_3,\ldots,f_{n-2}) \in (\R \setminus \{0\} )^{n-4}$, $\Phi,g \in Gl(2,\R)$. Then with respect to a basis
\[ 
\beta_{\h} =\{W_1^*,\ldots,W_{n-2}^*\} \qquad \beta_{\m} =\{ Z,W_1,\ldots,W_{n-2},Z^* \}
\]
the metric and Lie algebra structure are given by
\[
\ba{rcl}
B|_{\h \times \h}&=& 
\left(
\begin{smallmatrix}
  \Phi_{11}& \Phi_{12}  &           &           &           \\
  \Phi_{21}& \Phi_{22}  &           &           &           \\
           &            &   f_3     &           &           \\
           &            &           & \ddots    &           \\
           &            &           &           & f_{n-4}   
\end{smallmatrix}
\right)\\
{}\\
B|_{\m \times \m}&=&
\left(
\begin{smallmatrix}
     &              &              &        &   1    \\
     &g_{11} &g_{12} &        &        \\
     &g_{21} &g_{22} &        &        \\
     &              &              & I_{n-4}&        \\
 1   &              &              &        &        
\end{smallmatrix}
\right)
\ea 
\hspace{2cm}
\begin{array}{l}
\,\begin{array}{rcl}
{[Z^*,W_{\alpha}]} &=&  \phantom{-}W_{\alpha}^*\\
\end{array}\\ {}\\
\left.
\begin{array}{rcl}
{[W_i^*,Z^*]} &=& \phantom{-} \Phi_{ik}\cdot W_k \\
{[W_i^*,W_j]} &=& - g_{ik} \Phi_{kj} \cdot  Z \\
\end{array} \right\} \quad i,j=1,2 \\ {}\\ 
\left.
\begin{array}{rcl}
{[W_{\alpha}^*,Z^*]} &=& \phantom{-}f_{\alpha} \cdot W_{\alpha} \\
{[W_{\alpha}^*,W_{\beta}]} &=& - \delta_{\alpha \beta} f_{\alpha} \cdot Z \\
\end{array} \right\} \quad \alpha \ge 3
\end{array} 
\]
According to \cite[page 98]{Boubel00} there are four possibilities for the matrices $g,\Phi$:
\begin{align*}
g &= \left( \begin{smallmatrix} -1 & \\ & 1      \end{smallmatrix} \right) & 
\Phi     &= \left( \begin{smallmatrix}  \phi &  \\ & \phi \end{smallmatrix} \right)  \\
g &= \left( \begin{smallmatrix} & \pm 1 \\ \pm 1 & \end{smallmatrix} \right) & 
\Phi     &= \left( \begin{smallmatrix} \phi & 1 \\ & \phi \end{smallmatrix} \right)  \\
g &= \left( \begin{smallmatrix} & 1 \\ 1 & \\                       \end{smallmatrix} \right) & 
\Phi     &= \left( \begin{smallmatrix} \phi_1 & \phi_2 \\ -\phi_2 & \phi_1 \end{smallmatrix} \right)  
\end{align*}
\end{empt} 
}
}

\subsubsection*{(IV) Triples with \mb{ $[\h,\z^{\perp}] \not\subset \z$  }}
\addcontentsline{toc}{\protect subsubsection}{(IV) $[\h,\z^{\perp}] \not\subset \z$  }

This is the last remaining case, which might give a small insight in what can happen in higher signatures as far as
adapted decompositions are concerned. On the other hand this case is quite exceptional from the point of view that the
``interesting'' subspaces $\z^*$ and $E^*$ are one-dimensional. 
This expresses in the fact that finally the Jacobi-identity imposes
no conditions on the parameters of such a symmetric triple. 
First of all proposition  \ref{prop:AdapDecomp} provides a basis 
$\beta_{\m}=\{Z,e,W_1,\ldots,W_r,U_1,\ldots,U_s,e^*,Z^*\}$ in $\m$ such that
\[
Mat_{\beta}(B)=\left(
\begin{matrix}
&&J\\
&I_{n-4}&\\
J&
\end{matrix} \right) , \;
Mat_{\beta}(\ad(h)|_{\m})=\left(
\begin{matrix}
M_E^h&M_W^h&M_U^h&M_{E^*}^h\\
     &     &     &-(M_W^h)^t \cdot J\\
     & 0   &     &-(M_U^h)^t \cdot J\\
     &     &     &-J \cdot (M_E^h)^t \cdot J\\
\end
{matrix} \right)
\]
where
{\small
\begin{gather*}
J=\left( \begin{matrix} 0&1\\1&0 \end{matrix} \right), \quad
M_E^h=\left( \begin{matrix}  0&-B(h,\ol{E})\\0&0\end{matrix} \right), \quad
M_{E^*}^h=\left( \begin{matrix} -B(h,\ol{E^*})&0\\0&B(h,\ol{E^*}) \end{matrix} \right)\\
M_W^h=\left( \begin{matrix} -B(h,\ol{W}_1) & \hdots & -B(h,\ol{W}_r) \\
  -B(h,\wh{W_1}) & \hdots & -B(h,\wh{W}_r)  \end{matrix} \right), \quad
M_U^h=\left( \begin{matrix} -B(h,\ol{U}_1) & \hdots & -B(h,\ol{U}_s) \\
  0  & \hdots & 0  \end{matrix} \right)
\end{gather*}
}

The notation for the entries of the above matrices derives from the following consideration:
\[
\h = [\z^*,E \oplus W \oplus U \oplus E^*] + [E^*,W]
\]
due to  \ref{prop:AdapDecomp}(vi) and the fact that $\dim(\z^*)=\dim(E^*)=1$. Thus
\[
\h = span \Big\{ \ol{X},\wh{w}\; \Big|\;\ol{X}=[Z^*,X],\;X \in \z^{\perp}/\z \; \text{ and } \; \wh{w}=[e^*,w],\; w \in W  \Big\}  
\]
As usual we can consider the structure coefficients of the Lie algebra concerning $[\h,\m]$ (resp. the entries of 
the curvature tensor) as metrical coefficients on $\h$. 
 
Note, that we have two distinguished $B|_{(W \oplus U)^2}$-selfadjoint, hence diagonalizable operators, namely 
\begin{equation}\label{eq:OpF}
F=pr_{W \oplus U} \circ \big( \ad(Z^*) \big)^2 \big|_{W \oplus U} \quad \quad 
pr_{W \oplus U} \circ \big( \ad(e^*) \big)^2 \big|_{W \oplus U}
\end{equation}
A priori we do not know something about the
eigenvalues of $F$ as for the analogous operators $F_{ij}$ in the least nilpotent case,
whereas $pr_{W \oplus U} \circ \big( \ad(e^*) \big)^2 \big|_{W \oplus U}$
is nilpotent, since $\h \oplus E^* \oplus W \oplus U \subset [\g,\g]$ and the first derivative of a solvable
Lie algebra is nilpotent. Thus it is identically zero:
\begin{equation}\label{eq:E^*}
\big[ [W,E^*],E^* \big] \subset \z \quad \Leftrightarrow \quad B ( \wh{W},\wh{W} ) =\{0\}
\end{equation}

\begin{itemize} \item
Jacobi-identity on $\m$: We cite only those relations which give new conditions:
\begin{align}
\big[[E^*,W],E \oplus U \big]=\{0\} \quad &\Leftrightarrow \quad B(\wh{W},\ol{E} \oplus \ol{U})=\{0\} \label{eq:EU}\\
\big[[e^*,W_i],W_j \big]=\big[[e^*,W_j],W_i \big] \quad &\Leftrightarrow \quad B(\wh{W}_i,\ol{W}_j)=B(\wh{W}_j,\ol{W}_i)
\label{eq:WbarWhat}
\end{align}
The other relations $\sigma_{1,2,3} \big[ [m_1,m_2],m_3 \big] = 0$ are either equivalent to the pair-symmetry or trivial.

\vspace{1em}
\item $\sigma_{1,2,3}\big[[m_1,m_2],h_3 \big]=0$: On the one hand this defines the Lie algebra structure on $\h$:
\[
[h_{12},h]= \big[ [m_1,m_2],h \big]= \big[ [h,m_2],m_1 \big]- \big[ [h,m_1],m_2 \big]
\]
On the other hand it imposes conditions on the structure coefficients, 
which essentially leads to a restriction of the dimension of $W$:
\begin{align*}
[\wh{W},\h] &= \big[[E^*,W],\h \big] \subset \big[\underbrace{ [\h,W] }_{ \subset \z \oplus E },E^* \big] - 
  \big[\underbrace{ [\h,E^*] }_{ \subset \z \oplus W },W \big] =\{0\} \\
[\ol{E} \oplus \ol{U},\h] &=\big[[\z^*,E \oplus U],\h \big] 
  \subset  \big[ \underbrace{[\h,E \oplus U] }_{ \subset \z },\z^*\big] - 
\underset{\subset E \oplus W \oplus U \oplus E^* }{\big[\underbrace{[\h,Z^*]},E \big]} =\{0\} 
\end{align*}
Hence
\begin{equation}\label{eq:Central}
\wh{W} + \ol{E} + \ol{U} \quad \text{ is central in } \h
\end{equation}

Setting $m_1=Z^*,\;m_2=W_i \;(resp.\; e^*),\;X \in E \oplus W \oplus U \oplus E^*$ yields:
\begin{align}
[\ol{W}_i,\ol{X}]=\big[[Z^*,W_i],\ol{X} \big] &= \big[[\ol{X},W_i],Z^* \big]- \big[[\ol{X},Z^*],W_i \big]\label{eq:WX} \\
 &= B(\ol{X},\wh{W}_i)\cdot \ol{e} - B(\ol{X},\ol{e})\cdot \wh{W}_i \notag \\
\notag \\
[\ol{e^*},\ol{X}]= \big[[Z^*,e^*],\ol{X} \big] &= \big[[\ol{X},e^*],Z^* \big]- \big[[\ol{X},Z^*],e^* \big] \label{eq:E^*X}\\
 &= -\sum_k B(\ol{X},\wh{W}_k)\cdot \ol{W}_k + \sum_k B(\ol{X},\ol{W}_k)\cdot \wh{W}_k \notag
\end{align}
The other Jacobi-identities concerning $\big[ [\h,\m],\m \big]$ are trivial due to (\ref{eq:E^*}) and (\ref{eq:EU}). 
\vspace{1em}
\item $\sigma_{1,2,3}[[h_1,h_2],m_3]=0$ and $\sigma_{1,2,3}[[h_1,h_2],h_3]=0$ are treated on page \pageref{RemainingJacobi}.
\end{itemize}

Before we proceed in the discussion of (\ref{eq:WX}) and (\ref{eq:E^*X}) it will be necessary to investigate the dimensions
of some subspaces of $\h$:
\begin{align}
&[\h,\m]=\z^{\perp} \cap \m:&
E^* = [\h,Z^*]_{E^*} &\;\Rightarrow\; \{0\} \ne B \big( [\h,Z^*],E \big)=B \big( \h,\ol{E} \big) \notag \\
&&& \;\Rightarrow\; \dim(\ol{E})=\dim(E)=1 \label{eq:dim(olE)}\\
&&U = [\h,Z^*]_{U} &\;\Rightarrow\; \{0\} \ne B \big( [\h,Z^*],U_i \big)=B(\h,\ol{U}_i) \notag \\  
&&&\;\Rightarrow\; \dim(\ol{U})=\dim(U) \label{eq:dim(olU)} \\
&[\h,\z^{\perp}]_W=W:&
 W = [\h,E^*]_W  &\;\Rightarrow\; \{0\} \ne B \big( [\h,E^*],W_i \big) =B(\h,\wh{W}_i) \notag \\
&&& \;\Rightarrow\; \dim(\wh{W})=\dim(W) \label{eq:dim(hatW)}
\end{align} 
Nevertheless the subspaces $\ol{E},\ol{U},\wh{W}$ might have non-zero intersections pairwise. Then:
\[
0 \underset{(\ref{eq:Central})}{=} [\ol{e^*},\ol{e}] \underset{(\ref{eq:E^*X})}{=}
-\sum_k B(\ol{e},\wh{W}_k)\cdot \ol{W}_k + \sum_k B(\ol{e},\ol{W}_k)\cdot \wh{W}_k 
\underset{(\ref{eq:EU})}{=}  \sum_k B(\ol{e},\ol{W}_k)\cdot \wh{W}_k 
\]
The linear independency of the $\wh{W}_k$ due to (\ref{eq:dim(hatW)}) implies 
\begin{equation}\label{eq:EW}
B(\ol{E},\ol{W})=\{0\}
\end{equation}
Further
\[
[\ol{W}_i,\ol{W}_j] \underset{(\ref{eq:WX})}{=}B(\ol{W}_j,\wh{W}_i) \cdot \ol{e} - B(\ol{W}_j,\ol{e}) \cdot \wh{W}_i
\underset{(\ref{eq:EW})}{=} B(\ol{W}_j,\wh{W}_i) \cdot \ol{e} 
\]
The left side is skew-symmetric in $i$ and $j$, whereas the right side is symmetric due to (\ref{eq:WbarWhat}), hence
$\ol{e} \cdot B(\ol{W}_i,\wh{W}_j)=0$. Then (\ref{eq:dim(olE)}) implies:
\begin{equation}\label{eq:WW}
B(\ol{W},\wh{W})=\{0\}
\end{equation}
Collecting the above relations shows $B(\wh{W},\wh{W}+\ol{W}+\ol{E}+\ol{U})=\{0\}$. 
Thus by (\ref{eq:dim(hatW)}) $\dim(W) =  \dim(\wh{W}) \le \dim(\ol{E^*}) \le 1 $, since $B$ is non-degenerate on $\h$.
On the other hand, if the triple at issue exists, necessarily $\dim(W) \ge 1$ (otherwise $[\h,\z^{\perp}] \subset \z$, which
leads into the preceeding case). Summarizing we get
\begin{equation}\label{eq:dimension}
\boxed{ \dim(W)=\dim(\wh{W})=\dim(\ol{E^*})=\dim(\ol{E})=1 }
\end{equation}
Further implications from (\ref{eq:WX}) are:
\begin{equation}\label{eq:E=W}
\begin{array}{ll}
\{0\} \underset{(\ref{eq:Central})}{=}[\ol{W},\ol{E}]= -B(\ol{E},\ol{E}) \cdot \wh{W} &\;{\Rightarrow}\; B(\ol{E},\ol{E})=\{0\}\\
\{0\} \underset{(\ref{eq:Central})}{=}[\ol{W},\ol{U}]= -B(\ol{U},\ol{E}) \cdot \wh{W} &\;{\Rightarrow}\; B(\ol{U},\ol{E})=\{0\} 
\end{array} \quad
\begin{array}{ll} \Longrightarrow & \ol{E}=\wh{W} \\ {\sst (\ref{eq:EU}),(\ref{eq:EW})} & \end{array}
\end{equation}
Thus equation (\ref{eq:WX}) turns into $[\ol{W},\ol{X}]=\{0\}$. Since $[\ol{W},\wh{W}]=\{0\}$ by (\ref{eq:Central}) we get 
\begin{equation}\label{eq:Wcentral}
B(\ol{W},\h)=\{0\}
\end{equation}
In particular the holonomy is abelian.
Finally (\ref{eq:E^*X}) gives:
\begin{equation}\label{eq:W=W}
\begin{array}{ll}
\{0\}\underset{(\ref{eq:Wcentral})}{=}[\ol{E^*},\ol{W}]\underset{(\ref{eq:WW})}=B(\ol{W},\ol{W})\cdot \wh{W} &\;{\Rightarrow}\; B(\ol{W},\ol{W})=\{0\} \\
\{0\}\underset{(\ref{eq:Central})}{=}[\ol{E^*},\ol{U}]=B(\ol{U},\ol{W})\cdot \wh{W} &\;{\Rightarrow}\; B(\ol{U},\ol{W})=\{0\} 
\end{array} \quad
\begin{array}{ll} \Longrightarrow & \ol{W} \subset \wh{W} \\ {\sst (\ref{eq:EW}),(\ref{eq:WW})} & \end{array}
\end{equation}
We return now to the operator $F$ from (\ref{eq:OpF})
\[
-B \big( \, F(W),W \oplus U \, \big) =B \big( \, [Z^*,W],[Z^*,W \oplus U] \, \big)= B(  \ol{W},\ol{W} + \ol{U} )
\underset{(\ref{eq:E=W}),(\ref{eq:W=W})}{\subset} B(\wh{W},\ol{E} + \ol{U}) \underset{(\ref{eq:EU})}{=} \{0\}
\]
which means $W \in \ker(F)$. Since $F$ is selfadjoint it leaves the orthogonal space $ker(F)^{\perp}\cap W \oplus U \subset U$ invariant.
We will fix a basis in $U$ consisting of eigenvectors of $F$, hence $Mat(F)=diag(0,f_1,\ldots,f_s)$ with respect to the 
orthonormal basis $\{w,U_1,\ldots,U_s\}$. Then we get:
{\small
\begin{align*}
\ad(\ol{e}) &= 
\begin{array}{cc}
\begin{array}{ccc}\hphantom{a} \dots & E^* & \z^* \end{array} & \\
\left( \begin{array}{cc|c} \hphantom{aa} & -a &  \\ \hline  &  & a \\  \vphantom{\vdots} &  &  \end{array} \right) & 
\begin{array}{c} \overset{\vphantom{A}}{\z} \\ E \\ \vdots \end{array} 
\end{array}, \quad 
\ad(\ol{U_i})= 
\begin{array}{cc}
\begin{array}{ccccc}  \hphantom{a} \dots & U_i & \dots & E^* & \z^* \end{array} & \\
\left( \begin{array}{cccc|c}
  \hphantom{aa}    & -f_i & \hphantom{aa} & -b_i &  \\ \hline & & &  & b_i \\ & & &\vphantom{\vdots} &  \\ & & & & f_i  \\ & & & & \vphantom{\vdots} \end{array}\right)& 
\begin{array}{c}\overset{\vphantom{A}}{\z} \\ {E} \\ \vdots \\ U_i \\ \vdots   \end{array} 
\end{array}  \\
\ad(\ol{e^*})&=
\begin{array}{cc}
\begin{array}{cccccccc}\; \z & \phantom{-} E & \phantom{-}W& \; U_1 & \, \dots & \, U_s & \phantom{-} E^* & \z^* \end{array} & \\
\left( \begin{array}{ccccccc|c}
 \hphantom{a} &  -a & -\ol{a} & -b_1 & \dots &  -b_s & -c &  \\ \hline & & -\wh{a} & & & & & c \\ & & & & & & \phantom{-}\wh{a} & \ol{a} \\ 
& & & & & & & b_1 \\ & & & & & & & \vdots \\ & & & & & & & b_s \\ & & & & & & & a \\ & & & & & & & \end{array}\right)& 
\begin{array}{c} \z \\  E \\ W \\ U_1 \\ \vdots \\ U_s \\ E^* \\ \z^*  \end{array} 
\end{array} \quad
\begin{array}{l} \ad(\ol{w})= \dfrac{\ol{a}}{a}\cdot \ad(\ol{e}) \\ \\  \ad(\wh{w})= \dfrac{\wh{a}}{a}\cdot \ad(\ol{e}) 
\end{array}
\end{align*}
}
where empty entries denote zeros and
\begin{align}\label{eq:abcf}
a&=B(\ol{e^*},\ol{e}) & b_i&=B(\ol{e^*},\ol{U}_i) \notag \\
\ol{a}&=B(\ol{e^*},\ol{w}) & c&=B(\ol{e^*},\ol{e^*}) \\
\wh{a}&=B(\ol{e^*},\wh{w}) & f_i&=B(\ol{U}_i,\ol{U}_i) \notag
\end{align}

\pagebreak
From the above immediately follows:
\begin{itemize}
\item The Jacobi-identities $\sigma_{1,2,3}[[h_1,h_2],m_3]=0$ and $\sigma_{1,2,3}[[h_1,h_2],h_3]=0$ are
trivially satisfied.\label{RemainingJacobi}
\item $a,\wh{a} \ne 0$ since $\dim(\ol{E^*})=\dim(\ol{E})=\dim(\wh{W})=1$ by (\ref{eq:dimension}).
\item $f_i \ne 0 \; \forall i$: otherwise $U_i-e/a \in \z(\g)$ which contradicts $\z(\g)=\z$.
\item $\ol{a},b_i,c \in \R$ are free parameters.
\end{itemize}
Further $\beta_{\h}=\{\ol{e},\ol{U}_1,\ldots,\ol{U}_r,\ol{e^*}\}$ is a basis in $\h$, since the corresponding endomorphisms
induced by the adjoint action on $\m$ are obviously linear independent. 
Thus a \emph{parametrisation} of the symmetric triples at issue is obtained.

It is worth to summarize the situation in terms of Jacobi-operators and their restrictions
\[
 J(X)   =\big( \ad(X) \big)^2\big|_{\m} \qquad  F_V(X) =pr_V \circ J(X)|_V \quad (\m=V \oplus \wt{V})
\]
\begin{itemize}
\item $F_{E \oplus U \oplus E^*}(Z^*)$ leaves $E$ invariant and is invertible. Conversely every selfadjoint operator with this
properties is predicted by the data $a,b_i,c,f_i$.
\item $F_{E \oplus W \oplus E^*}(Z^*)$ is predicted by the data $a,\ol{a},c$.
\item $F_{\z \oplus W \oplus z^*}(e^*)$ is nilpotent and predicted by the data $\wh{a},c$.
\end{itemize}

Since all spaces the above operators are acting on are Lorentzian the usual spectral theorem for selfadjoint operators
does not work. On the other hand there are two difficulties which prevent an immediate application of the general theory of \cite{Boubel00}:
First $e^*,Z^*$ are not canonical given due to the ambiguity of the adapted decomposition. And second it is not clear
(in contrast to the definite case) what a simultaneous normal form for $J_{e^*},J_{Z^*}$ might be. So we will
continue with the usual procedure of discussing the effect on the structure coefficients 
while changing the adapted decomposition:

The \emph{first} adapted decomposition $\m = \z \oplus W_0 \oplus \z^*$ neglects the finer structure of 
$W_0=E \oplus W \oplus U \oplus E^*$. Transformation into a different decomposition
$\m = \z \oplus \wt{W}_0 \oplus \wt{\z^*}$ can be described by a homomorphism $L \in Hom(W_0,\z) \simeq W_0^*$ since $\dim(\z)=1$, 
hence $Skew(\z^*,\z,B)=\{0\}$
\begin{description}
\item $L_{E^*}: e^* \mapsto e^* + l \cdot Z$ is induced by $\ad(\ol{e})$.
\item $L_{U_i}: U_i \mapsto U_i + l \cdot Z$ is induced by $\ad(\ol{U}_i)$ modulo $\ad(\ol{e})$.
\item $L_{E}: e \mapsto e + l \cdot Z$ is induced by $\ad(\ol{e^*})$ modulo $\ad(\ol{e})$, $\ad(\ol{U}_i)$ and a transformation
in $W_0$ which is not of interest at this point.
\end{description}
Thus the only change which might effect the coefficients is
\begin{align*} 
w &\mapsto w+l \cdot Z \\ Z^* &\mapsto Z^*-l \cdot w - l^2/2 \cdot Z
\end{align*}
The induced change in $\beta_{\h}$ is $\ol{e^*} \mapsto [Z^*-l \cdot w , e^*]=\ol{e^*}+l \cdot \wh{w}$.
The only parameter which is not invariant under this transformation is
\begin{equation}\label{eq:c=0}
\wt{c}=B(\ol{e^*}+l \cdot \wh{w},\ol{e^*}+l \cdot \wh{w}) \underset{(\ref{eq:abcf})}{=}c + 2l \cdot \wh{a}
\end{equation} 
Thus we can achieve $c=0$ without a change of the other parameters ($\wh{a}\ne 0$). Whenever $c$ becomes unequal zero
by a further transformation it can be reset to zero. Thus we assume from now $c=0$ without mentioning the possibly necessary
additional transformations in the following.

Now we play the same game for the \emph{second} decomposition $W_0=E \oplus W \oplus U \oplus E^*$. 
\begin{align*} 
w &\mapsto w+l \cdot e \\ e^* &\mapsto e^*-l \cdot w - l^2/2 \cdot e
\end{align*}
Since $\ol{a} \cdot \ol{e} = a \cdot \ol{w}$ by (\ref{eq:abcf}), the choice  $l=-\ol{a}/a$ yields 
$ \ol{w} \mapsto [Z^*,w - \ol{a}/a \cdot e]=\ol{w} - \ol{a}/a \cdot \ol{e} =0$ (note $a \ne 0$), hence we can achieve
\begin{equation}\label{eq:abar=0}
\ol{a}=0 
\end{equation}
For the transformation
\begin{align*} 
U_i &\mapsto U_i+l \cdot e \\ e^* &\mapsto e^*-l \cdot U_i - l^2/2 \cdot e
\end{align*}
we have to distinguish two cases:
\begin{description}
\item[$f_i \ne a$] In this case $b_i=0$ can be achieved by a suitable transformation 
( note: $ b_i \mapsto l(a -f_i)+b_i$).
\item[$f_i = a$] By a permutation we can assume, that $f_i=a$ for $0 \le i \le s_a \le s$. 
Then $\beta_a=\{e,U_i,e^*\}_{0 \le i \le s_a}$ is a basis in the eigenspace of $J_{Z^*}$ with respect to the eigenvalue $-a$.
If $s_a=0$ there is nothing to do. Otherwise it is easy to see that by an orthogonal transformation in 
$span\{U_1,\ldots,U_{s_a}\}$ we can achieve $(b_1,\ldots,b_{s_a})=(b,0,\ldots,0)$.  
\end{description}

Finally we can scale the Witt-basis in $\z \oplus \z^*$ and $E \oplus E^*$. A unique choice (up to $\pm 1$) can be forced by the
condition
{\small
\[ \left\{
\begin{array}{ll} \wh{a}= |a|=1  & \text{ for } b=0 \text{ or } s_a=0 \\ \wh{a}=b=1 & \text{ for } b\ne 0 \end{array} \right.
\; \begin{array}{c}\text{\footnotesize with the free } \\ \text{\footnotesize parameters } \end{array}\; \left\{
\begin{array}{l}   f_1,\ldots,f_{s} \in \R \setminus \{0\} \\ (a=f_1),f_2,\ldots,f_{s} \in \R \setminus \{0\} \end{array} \right.
\]
}

In terms of $F$ this means
\begin{itemize}
\item $F=F_{E \oplus U \oplus E^*}(Z^*)$ is either diagonalizable with $F_{E \oplus E^*}(Z^*)= \pm id_{E \oplus E^*}$, 
or it contains one Jordan-block of length $3$:
$F_{E \oplus U_1 \oplus E^*}(Z^*)=-\left(\begin{smallmatrix} a&1&\\&a&1\\&&a  \end{smallmatrix}\right)$
\item $F_{E \oplus W \oplus E^*}(Z^*) =-\left(\begin{smallmatrix} a&&\\&0&\\&&a  \end{smallmatrix}\right)$.
\item $F_{\z \oplus W \oplus \z^*}(e^*)=-\left(\begin{smallmatrix}0&1&\\&0&1\\&&0  \end{smallmatrix}\right)$
\end{itemize}
Thus the resulting symmetric triple can be described completely in terms of $F$:
\begin{dfn}[\mb{$\tau_n(F)$}]
Consider the following class of symmetric triples $\tau_n(F)=(\h \oplus \m , B)$:
\bi
\item $\m =\z \oplus E  \oplus W \oplus U \oplus E^* \oplus \z^*= span\{Z,e,w,U_1,\ldots,U_{n-5},e^*,Z^*\}$
\item $\h = \ol{E} \oplus \ol{U} \oplus \ol{E^*}$  where $ X \mapsto \ol{X}$ denotes an isomorphism
\item $B|_{(\z \oplus E)^2} = B|_{(\z^* \oplus E^*)^2} = 0,\quad B(Z^*,Z) = B(e^*,e)=1$,
$B|_{(W \oplus U)^2}$ is positive definite and 
 $\z \oplus \z^*,E \oplus E^* ,W ,U$ are mutual orthogonal 
\item $F \in Gl({E} \oplus {U} \oplus {E^*})$  a selfadjoint invertible map which leaves $E$ invariant  
\item The non-zero commutators of the Lie algebra structure are
\[
\begin{array}{rcl}
{[Z^*,X]} &=&  \phantom{-}\ol{X}\\
{[e^*,w]} &=& -B \big( F(e),e^* \big)^{-1} \cdot \ol{e}\\
\end{array}
\qquad
\begin{array}{rcl}
{[\ol{X},Z^*]} &=& \phantom{-}F(X)\\
{[\ol{X},Y]}={[\ol{Y},X]}&=& -B \big( F(X),Y \big) \cdot Z \\
{[\ol{e^*},e^*]} &=& \phantom{-}w\\
{[\ol{e^*},w]} &=& -e
\end{array}
\]
where $X \in E \oplus U \oplus E^* $ and $Y \in E \oplus U$
\ei
\end{dfn}

\begin{prop}
If $\tau$ is a solvable triple of signature $(p,n-p)$ with $1$-dimensional center and
$[\h,\z^{\perp}] \not\subset \z$ there exists a triple $\tau_n(F)$ which is isomorphic to $\tau$ and conversely.
\end{prop}
for completeness we cite also a normal form for this case:

\mb{ (IV) \hspace{.3cm} $\dim(\z)=1$, $[\h,z^{\perp}] \not\subset \z$: }\\
\fbox{
\parbox{15cm}{
\begin{empt}\label{triple:tau(a,b,f)}
\mb{Normal form  $\tau_{n}(a,b,f)$}
Consider parameters $a,f_1,\ldots,f_{n-5} \in \R \setminus \{0\}$, $b \in \R$. Then with respect to a basis
\[ 
\beta_{\h} =\{\ol{e},\ol{U}_1,\ldots,\ol{U}_{n-5},\ol{e^*}\} \quad
\beta_{\m} =\{Z,e,w,U_1,\ldots,U_{n-5},e^*,Z^*\}
\]
the metric and Lie algebra structure are given by
\[
\ba{rcl}
B|_{\h \times \h}&=& 
\left(
\begin{smallmatrix}
  0     & 0   &\dots  &   0    &a \\
   0    & f_1 &       &        &b \\
\vdots  &     &\ddots &        &\vdots \\   
  0     &     &       &f_{n-5} & 0  \\
 a      & b   &\dots  &    0   &0  
\end{smallmatrix}
\right)\\
{}\\
B|_{\m \times \m}&=&
\left(
\begin{smallmatrix}
    &   &        &   & 1 \\
    &   &        & 1 &   \\
    &   &I_{n-4} &   &   \\
    & 1 &        &   &   \\ 
  1 &   &        &   &   
\end{smallmatrix}
\right)\ea 
\quad
\begin{array}{rcl}
{[Z^*,e]} &=&  \ol{e}\\
{[Z^*,U_i]} &=& \ol{U}_i\\
{[Z^*,e^*]} &=& \ol{e^*}\\
{[e^*,w]} &=&1/a \cdot \ol{e}
\end{array} \quad
\begin{array}{rcl}
{[\ol{e},Z^*]} &=& \phantom{-} a \cdot e\\
{[\ol{U}_1,Z^*]} &=& \phantom{-}b \cdot e + f_1 \cdot U_1\\
{[\ol{U}_{\alpha},Z^*]} &=& \phantom{-}f_{\alpha} \cdot U_{\alpha} ,\quad (\alpha > 2)\\
{[\ol{e^*},Z^*]} &=& \phantom{-}b \cdot U_1 + a \cdot e^*\\
{[\ol{e^*},e^*]} &=& \phantom{-}w\\
{[\ol{e^*},w]} &=& -e\\
{[\ol{U}_1,e^*]} &=& -b \cdot Z\\
{[\ol{e^*},U_1]} &=& -b \cdot Z\\
{[\ol{e},e^*]} &=& -a \cdot Z\\
{[\ol{e^*},e]} &=& -a \cdot Z\\
{[\ol{U}_{\alpha},\ol{U}_{\beta}]} &=& -\delta_{\alpha\beta}f_{\alpha} \cdot Z\\
\end{array}
\]
There are two possibilities for the parameters:
\begin{enumerate}
\item[(i)] $b=0$ and $|a|=1$. $\tau(a,0,f),\tau(\wt{a},0,\wt{f})$ are isomorphic if there exists a permutation $\Pi \in S_{n-5}$
such that $\wt{f}_i=f_{\Pi(i)}$ and $a=\wt{a}$. 
\item[(ii)] $b=1$ and $a=f_1$. $\tau(f_1,1,f),\tau(\wt{f}_1,1,\wt{f})$ are isomorphic if there exists a permutation 
$\Pi \in S_{n-6}$ such that $\wt{f}_i=f_{\Pi(i)},\;i \ge 2$ and $\wt{f}_1=f_1$. 
\end{enumerate}
\end{empt} 
}
}

\begin{re}
This triple corresponds to (5.6) in \cite{CahenP70}. There seems to be a little mistake: 
The holonomy action as cited is not faithful -
more precisely (in the notation of \cite{CahenP70}): 
$n^*,y^*$ are linear dependent. This becomes obvious observing that $y \in W_2$ (The definition of $W_2$ as
complement of  $u$ in $W_1$ is not sufficient). Moreover one can find a ``better'' adapted decomposition. 
\end{re}

\pagebreak
\subsection[Conclusion and table of normal forms]
{Conclusion and table of normal forms}

Using the concept of iterated adapted decompositions we have been able to reproduce the classifications of solvable
symmetric triples with signature $(1,n-1)$ and $(2,n-2)$ obtained in \cite{CahenW70} and \cite{CahenP70}. 
In particular
the group $\Gamma$ acting transitively on the adapted decompositions is of geometric interest since it contains the
linear isotropy representation of the isometry group and on the other hand it is a quite helpful tool in order to
find normal forms for solvable symmetric triples. The classification of the triples in question is enabled mainly by
the simplicity of the Jacobi-identity for low signatures. 
As example for the difficulties which arise in higher signature the non-least-nilpotent case with maximal center has been pointed out. 
The dimensions of the subspaces of $\m$ occurring in the iterated adapted decomposition
are not arbitrary as the case (IV) showed. This indicates that the Jacoby-identity might be quite restrictive for
triples with multiple iterated adapted decompositions. 

The results suggest that a complete classification of indecomposable solvable symmetric triples
in the sense of determining their isomorphism classes is probably impossible (cf.  \cite[page 354]{Wu67}).
On the other hand in a recent work by I.Kath and M.Olbricht \cite{KathO02} more structural results has been uncovered by means of
Lie algebra cohomology - an idea inspired by L.Berard-Bergery. 
A question in this context is to find distinguished adapted decompositions. 
For example a less computational argument for the separation property \ref{lem:afSeperation} would be desirable.

The holonomy algebras of solvable symmetric spaces in signature $(1,n-1)$ and $(2,n-2)$ are all abelian.
In \cite[Theorem 2]{Wu67} a double-extension of $\{0\}$ by a solvable Lie algebra (which yields signature $(p,p)$) has been used to show that
there exists solvable symmetric triples with non-abelian holonomy. It seems to be that
there is nothing known about examples of solvable symmetric triples with non-abelian holonomy which does \emph{not} come from a 
Lie group with a bi-invariant metric. 

{\footnotesize
\begin{center}
\renewcommand{\arraystretch}{2}
\begin{tabular}{| l | l | l | l |  p{2.2cm} | l | p{4.5cm}  | } \hline
Page&  Type    & Signature      &$\dim(\z)$     &   Description      & Triple   &  Ricci-curvature        \\ \hline \hline

\pageref{thm:Main}   &  & $(p,q) $    & $p$  & least nilpotent: \linebreak - general case & $\tau_{p,q}(\RR_b,\epsilon,\Lambda)$ & 
   \hspace{-.3cm} \raisebox{-1mm}{$ \renewcommand{\arraystretch}{1} \begin{array}{l} 
   \RR ic(Z^*,\wt{Z^*})= \\ \phantom{=}\sum_{\alpha} \; \epsilon_{\alpha} B(\Lambda_{\alpha},Z^*)B(\Lambda_{\alpha},\wt{Z^*}) \end{array} $}    
\\ \cline{1-4} \cline{6-7}

\pageref{triple:tau(f)}  &  & $(1,n-1)$  &  $1$   & -  Lorentzian    &   $\tau_n(f)$  &    
    {$\RR ic(Z^*,Z^*)=\sum_{\alpha} \; f_{\alpha}$}
\\ \cline{1-4} \cline{6-7}           

\pageref{triple:tau(epsilonY)}   &   (Ia)     & $(2,n-2)$    & $2$   & - $\dim(Y)=1$  & $\tau_n (\epsilon_Y,\epsilon,\lambda)$ &  
          {$\RR ic(Z_i^*,Z_j^*) = \sum_{\alpha} \; \epsilon_{\alpha} \lambda_i^{\alpha}\lambda_j^{\alpha}$   }
\\ \cline{1-2} \cline{6-6}

\pageref{triple:tau(epsilon)} &    (Ib)     &    &   & - $\dim(Y)=0$  &  $\tau_n (\epsilon,\lambda)$ &   
\\ \cline{1-2} \cline{5-6}
 
\pageref{triple:tau(If=1)}  &   (IIa)    &    &    & $\dim(W_{nil})=1$ & $\tau_n (\epsilon,\lambda,1)$ &  
\\ \cline{1-2} \cline{5-6}

\pageref{triple:tau(If=2)}  &  (IIb)    &  &    &  $\dim(W_{nil})=2$ & $\tau_n (\epsilon,\lambda,2)$ &  
\\ \cline{1-3} \cline{5-7}

\pageref{triple:tau(2,2)}    & (Ia')    & $(2,2)$    &   &  nilpotent  & $\tau_{2,2}(\epsilon_Y)$ &    {$\RR ic \equiv 0$   }
\\ \cline{1-3}  \cline{6-6}

\pageref{triple:tau(2,3)}    & (IIa')    & $(2,3)$    &     & & $\tau_{2,3}$ &    
\\ \cline{1-3}  \cline{6-6}

\pageref{triple:tau(2,4)}    &  (IIb')    & $(2,4)$    &  &  & $\tau_{2,4}$ &   
\\ \hline

\pageref{triple:tau(Phi)}    &  (III)    & $(2,n-2)$    & 1 & $[\h,\z^{\perp}] \subset \z$ & $\tau_n(g,\Phi,f)$ &  
   {\small   $\RR ic(Z^*,Z^*)=tr(\Phi) + \sum f_{\alpha}$ }
\\ \cline{1-2}  \cline{5-7}

\pageref{triple:tau(a,b,f)}    &  (IV)    &    &  & $[\h,\z^{\perp}] \not\subset \z$ & $\tau_n(a,b,f)$ &  
 { $\RR ic(Z^*,Z^*)= 2a + \sum f_{\alpha}$ }
\\ \hline

\end{tabular}

\end{center}
}

\section{On the Topology of Symmetric Spaces}

\subsection{Symmetric Coverings}\label{ss:SymmetricCovering}
Simply connected symmetric spaces are in one-to-one correspondence to symmetric triples. In this section 
coverings of symmetric spaces and the topological structure of non-simply-connected symmetric spaces will be investigated. 

\begin{dfn}
A  pseudo-Riemannian covering map $c: (M^1,g^1) \rightarrow (M^2,g^2)$ 
is said to be a \emph{symmetric covering} if $(M^i,g^i),\; i=1,2$ are both symmetric spaces.
\end{dfn}
In general the discrete quotient of a symmetric space is only \emph{locally} symmetric.
The symmetries are denoted by $s^{(i)}_x$, the transvection groups by $\G^i$ and the conjugation with $s_{o^i}^{(i)}$ by $\Sigma^{(i)}$
for base points $c(o_1)=o_2$.
\begin{prop}\label{prop:SymCov}
Let $c:(M^1,g^1) \rightarrow (M^2,g^2)$ be a symmetric covering. Then
\begin{enumerate}
\item[(i)] $c \circ s^{(1)}_x = s^{(2)}_{c(x)} \circ c$
\item[(ii)] The map $C: \G^1 \rightarrow \G^2$ induced by $C(s^{(1)}_x \circ s^{(1)}_y)=s^{(2)}_{c(x)} \circ s^{(2)}_{c(y)}$ is a 
Lie group homomorphism and a covering map. In particular $\ker(C) \subset Z(\G^1)$ is discrete.
\item[(iii)] $C(g) \circ c = c \circ g, \; \forall g \in \G^1$ and $C \circ \Sigma^{(1)} = \Sigma^{(2)} \circ C$. 
\end{enumerate}
\end{prop}

\begin{proof}
(i) follows immediately from the fact that the two local isometries
$c \circ s^{(1)}_x$ and $s^{(2)}_{c(x)} \circ c$ and their differentials coincide at $x$.

(ii) First of all it should be mentioned that $C$ is well defined since $c$ is onto (cf. \cite[I.3.5]{Bertram00}).
To prove smoothness of $C$ we use the well-known fact that under the submersion $\pi^i:\G^i \rightarrow M^i$, $g \mapsto g(o^i)$
1-parameter subgroups corresponding to $X \in \m^i$ are mapped to geodesics, more precisely 
$\pi^i\big(exp_{\G^i}(tX)\big)=exp_{M^i}\big(d \pi^i (tX)\big)$. Since $c$ is a local isometry $c \circ exp_{M^1}=exp_{M^2} \circ dc$
an thus (cf. proof of \ref{prop:CharTrans}):
\begin{multline*}
(c \circ \pi^1) \big(exp_{\G^1}(X)\big) = exp_{M^2} \big( (dc \circ d\pi^1) X \big) 
    = \pi^2  exp_{\G^2} \Big( \big( (d\pi^2|_{\m^2})^{-1} \circ dc \circ d\pi^1 \big) X \Big)\\
\Rightarrow \qquad  C\big( exp_{\G^1}(X) \big) = C \Big( s^{(1)}_{\pi^1 \big( exp_{\G^1}(X/2) \big)} \circ s^{(1)}_{o^1} \Big)
= exp_{\G^2} \Big( \big( (d\pi^2|_{\m^2})^{-1} \circ dc \circ d\pi^1 \big) X \Big)
\end{multline*}
The proof is completed by the observation that $\G^i$ is generated by $exp_{\G^i}(\m^i)$ 
and that the invertible map $(d\pi^2|_{\m^2})^{-1} \circ (dc) \circ (d\pi^1)|_{\m^1}:\m^1 \rightarrow \m^2$ 
extends to a unique Lie algebra isomorphism $dC:\g^1 \rightarrow \g^2$.

(iii) is verified by direct calculations using (i). 
\end{proof}

Now we specialize to the case where $M^1$ is simply connected. Then $(M^2,g^2)$ is isometric to the quotient 
$(M^1/D,g^1/D)$ by the automorphism group $D$ of the covering $c$
\[
D=Aut(c,M^1,M^2)= \menge{ \delta \in Iso(M^1,g^1)}{ c \circ \delta = c }
\]
where $g^1/D$ denotes the metric on the factor space $M^1/D$ by requiring the covering map $M^1 \rightarrow M^1/D$ to be a local isometry.

\pagebreak

\begin{lem}\label{lem:SymCov}
\be
\item[(i)] Let $c:(M^1,g^1) \rightarrow (M^2,g^2)$ be a symmetric covering. Then $D=Aut(c,M^1,M^2)$ satisfies
$D \subset  Z_{Iso(M^1,g^1)}\G^1$ and $\Sigma^{(1)}(D) \subset D$.
\item[(ii)] If conversely $(M,g)$ is a symmetric space and $D \subset Z_{Iso(M,g)}\G(M,g)$ a discrete subgroup 
which is invariant under $\Sigma=C_{s_o}$ and acts properly
discontinuously on $M$, then $c: M \rightarrow M/ D$ induces a symmetric covering. 
\ee
\end{lem}

\begin{proof}
(i) $c \circ g  \delta  g^{-1} = C(g)  c  \delta  g^{-1} = C(g)  c  g^{-1} = c,\; \forall \delta \in D,\;\forall g\in \G^1$,
hence $ g  \delta  g^{-1} \in  D$. Further the map $\Delta_{\delta}:\G^1 \rightarrow D$, $g \mapsto g  \delta  g^{-1}$ is continuous.
Since $\G^1$ is connected and $D$ is discrete, $\Delta_{\delta}$ is constant. 
Then $\Delta_{\delta}(e)=\delta$ implies  $ g  \delta  g^{-1} = \delta,\; \forall g \in \G^1$.
Secondly note 
$c  s_{o_1}^1  \delta  s_{o_1}^1 = c  s_{o_1}^1   s_{\delta(o_1)}^1  \delta =
 c  \delta =  s_{c(o_1)}^2   s_{(c  \delta)(o_1)}^2   c  \delta = c$, 
hence $s_{o_1}^1  \delta  s_{o_1}^1 = \Sigma^{(1)}(\delta) \subset D$.

(ii) The topological assumption about $D$ provides a pseudo-Riemannian covering $c:M \rightarrow M/D$. 
Thus it suffices to derive from the algebraic conditions on $D$ the symmetry of $M/D$:
An isometry $\wt{\Phi} \in Iso(M,g)$ induces an isometry $\Phi \in Iso(M/D,g/D)$ (i.e. $c \circ \wt{\Phi}=\Phi \circ c$) 
if and only if 
$\wt{\Phi} \in N_{Iso(M,g)}D= \big\{ \Psi \in Iso(M,g) \;|\; \Psi \, D  \, \Psi^{-1} \subset D\}$. In particular 
$\Sigma(D)=D$ signifies $s_o \in N_{Iso(M,g)}(D)$. This provides a symmetry at the base point
$o D \in M/D$. Further 
$D \subset Z_{Iso(M,g)}\G(M,g)$ implies $\G(M,g) \subset Z_{Iso(M,g)}D \subset  N_{Iso(M,g)}D $. 
Thus we have transvections along arbitrary geodesics in $M/ D$. 
Consequently $M/ D$ is homogenous which yields a symmetry at every point of $M/D$. Hence $M/D$ is symmetric.
\end{proof}

\begin{re}
$D \subset Z_{Iso(M,g)}\G(M,g)$ ensures that $D$ acts freely on $M$: Assume $\delta(x)=x$ for some $x \in M$. Then for every $y \in M$ exists
a $g \in \G(M,g)$ with $y=g(x)$, hence $\delta(y)=\delta \big(g(x)\big)=g \big( \delta (x) \big)=y$. 
Since $\G(M,g) \subset Iso(M,g)$ acts effectively on $M$, $\delta$ is the identity in $\G$. \hfill \qedsymbol
\end{re}
Lemma \ref{lem:SymCov} is not completely satisfying for several reasons: First one cannot not expect the isometry group to be computable in general.
The essential idea to avoid this problem is to consider only discrete subgroups in the center of the transvection group 
$Z\big(\G(M,g) \big) \subset Z_{Iso(M,g)}\G(M,g)$. 
We will discuss in the following - at least topologically - what we lose by this restriction. 
The results and technics are due to \cite{Koh65} and \cite{Loos69}. 
Further it would be nice to have a better criterion for the action being properly discontinuous. 
In the Riemannian case this turns out to be an empty condition:

\begin{re}
Let $(M,g)$ be a Riemannian homogenous manifold. Then every discrete subgroup $D \subset Iso(M,g)$ which acts freely on $M$ already acts
properly discontinuously:

We represent $M$ as $Iso(M,g)/K$, where $K$ denotes the isotropy subgroup with respect to a fixed base point. 
The action of $D$ on $M$ corresponds
under this identification to $\delta(gK) \mapsto (\delta g)K,\; \delta \in D$. It is known that $K$ is compact in the Riemannian case 
\cite[Ch.VI, thm. 3.4]{KobayashiN2}. Then the action of $D$ on $Iso(M,g)/K$ and thus on $M$ is properly discontinuous due to 
\cite[Ch.I, prop. 4.4 \& 4.5]{KobayashiN2}. \hfill \qedsymbol
\end{re}

\begin{re}The same has been asserted for symmetric spaces in the Pseudo-Riemannian context \cite[Ch.I, prop. 3.9]{CahenP80}, 
but we have not been able to complete the proof. 
However, assume for every $x \in M$ exists a neighborhood 
$U_x \subset M$ such that $(D \setminus \{e\}) (U_x) \cap U_x = \emptyset$ 
( i.e. $c: M \rightarrow M/D$ becomes a topological covering with the factor topology on $M/D$ ). 
Then the factor space is already Hausdorff, hence a manifold:

Let $x,y \in M$ with $c(x) \ne c(y)$ i.e. $y \not\in D(x)$. Further let $U_x$ be a neighborhood of $x$ as presumed above. 
If $y \in \delta(U_x)$ for some $\delta \in D$  we can choose disjunctive neighborhoods $V_y,\,V_{\delta(x)}$ of $y,\,\delta(x)$ within $\delta(U_x)$.
Then $D(V_y) \cap D(V_{\delta(x)}) = \emptyset$. In other words the points $c(x),c(y)$ can be separated by open sets in $M/D$.
Now assume $y \not\in D(U_x)$. For $z \in M$ consider the 
submersion $\pi_z :\G(M,g) \rightarrow M,\; g \mapsto g(z)$. Then ${\pi_x}^{-1}(U_x) \subset \G$ is an open subset which contains the
identity ($U_x$ as above) and there exists a possible smaller neighborhood $W \subset \G$ of $e$ with $W^{-1}\cdot W \subset {\pi_x}^{-1}(U_x)$. Since a submersion
is an open map, the subsets $\pi_x(W)=W(x)$ and $\pi_y(W)=W(y)$ are again neighborhoods of $x$ resp. $y$ and 
\[
W(y) \cap D \big(W(x) \big) = W \Big(  y \cap W^{-1} D \big( W(x) \big) \Big) 
\underset{D \subset Z_{Iso}\G}{=}W \Big(  y \cap D \big( W^{-1}W(x) \big) \Big) \subset 
  W \big(  y \cap D (U_x) \big) = \emptyset
\]
hence $c(W(x)),c(W(y))$ are disjunctive. \hfill \qedsymbol
\end{re}

Now we return to the situation of a given symmetric covering $c:(M^1,g^1)  \rightarrow (M^2,g^2)$ where $M^1$ is assumed to be simply connected.
We represent these symmetric spaces by their transvection groups $\G^i$ and the isotropy subgroups with respect to base points $o^i$:
\[
{}^{\ds \G^i}/_{\ds \H^i} \rightarrow  M^i, \quad g\H^i \mapsto g(o^i)
\]
Then $\H^1$ is connected due to the exact homotopy sequence of $\G^1 / \H^1$.
From $C \circ \Sigma^{(1)}=\Sigma^{(2)} \circ C$ (\ref{prop:SymCov}(iii)) we get $C\big( (\G^1)^{\Sigma^{(1)}} \big) \subset (\G^2)^{\Sigma^{(2)}}$.
Then $(dC)_{e^1}$ being a Lie algebra isomorphism yields:
\begin{equation}\label{eq:H1H2}
C ( \H^1 )=C\big( (\G^1)^{\Sigma^{(1)}}_0 \big) = C\big( (\G^1)^{\Sigma^{(1)}} \big)_0 = (\G^2)^{\Sigma^{(2)}}_0 = (\H^2)_0 \subset \H^2
\end{equation}
The covering $c$ is then described by (see \ref{prop:SymCov}(iii))
\[
\ol{C}:\G^1/\H^1 \rightarrow \G^2/\H^2,\quad \ol{C}(g\H^1)=C(g)\H^2
\]
which is well defined due to (\ref{eq:H1H2}).
This map naturally splits into two coverings:

\begin{prop}\label{prop:2Cov}
\begin{enumerate} 
\item[(i)] Every symmetric covering {\small $\,c:(M^1,g^1)  \rightarrow (M^2,g^2)$}, $M^1$ simply connected yields 
the following commutative diagram of symmetric coverings:
\[
\begin{diagram}
\node[1]{M^1 \simeq {}^{\ds \G^1}/_{\ds \H^1}} 
\arrow[2]{e,t}{\ol{C}} 
\arrow{se,l}{C_{\G}}
\node[2]{{}^{\ds \G^2}/_{\ds \H^2} \simeq M^2 } 
\\
\node[2]{ {}^{\ds \G^2}/_{\ds \H^2_0}} \arrow{ne,r}{C_{\H}}
\end{diagram}
\]
$C_{\G}(g\H^1)=C(g)\H^2_0$ and $C_{\H}(g\H^2_0)=g\H^2$. Moreover $\G^2$ is the transvection group of $\G^2/\H^2_0$.
\item[(ii)] The covering $C_{\G}$ is induced by the discrete subgroup $\ker(C) \subset Z(\G^1)$.
\item[(iii)] The covering $C_{\H}$ is finite . 
Moreover the group $\H^2_0/\H^2$ is isomorphic to $(\Z_2)^k$ for some $k \in \N$. Hence $C_{\H}$ 
can be comprehended as a sequence of $k$ coverings with fibers of cardinality $2$.
\end{enumerate}
\end{prop}
\begin{proof}
(i): $M^{1,2}:=\G^2/\H^2_0$ becomes a symmetric space 
with the metric $g^{1,2}$ corresponding to the $Ad(\H^2)$- and thus $Ad(\H^2_0)$-invariant scalar product on $\m^2$.
Thus $\G^2$ is a subgroup of $Iso(M^{1,2},g^{1,2})$. Being connected it coincides with the transvection group. 
Further the map $C_{\G}$ is well defined due to (\ref{eq:H1H2}). Then by the very definition 
$C = C_{\H} \circ C_{\G}$ and all three maps are local isometries. 
It is well known that $C_{\H}$ is a pseudo-Riemannian covering map. Then it is clear that $C_{\G}$ is a covering as well, since $\ol{C}$ 
is one by assumption. 

(ii): $C|_{\H^1}$ is injective since $\G^1$ acts effectively on $\G^1/\H^1$ and thus $\ker(C) \cap \H^1 \subset Z(\G^1) \cap \H^1 = \{e\}$.
Thus $C|_{\H^1}:\H^1 \rightarrow \H^2_0$ is bijective and 
${C_{\G}}^{-1}\big(g\H^2_0 \big)=C^{-1}(g)H^1=\ker(C)\big(gH^1\big)$.

(iii): Of course this is the heart of the proposition. It follows directly from a theorem due to S. Koh cited hereafter. 
\end{proof}

\begin{thm}[S.Koh]\label{thm:Koh}
Let $G$ be a connected Lie group, $\Sigma$ an involutive automorphism of $G$ and $\g = \h \oplus \m$ the usual decomposition of $\g=LA(G)$
into the $\pm 1$-eigenspaces of $d\Sigma$. Then there exists a maximal compact $\Sigma$-stable subgroup $K \subset G$. Furthermore there exists
linear subspaces $\h_1,\ldots ,\h_p \subset \h$ and $\m_1, \ldots , \m_q \subset \m$ with $Ad\big(K^{\Sigma}\big)(\m_i) \subset \m_i$ such that
the following map is a diffeomorphism:
\begin{align*}
K \times \m_1 \times \ldots \times \m_q  \times \h_1 \times \ldots \times \h_p \quad &\rightarrow \quad G \\
(k,a_1,\ldots,a_p,b_1,\ldots,b_q) \quad &\mapsto \quad
k \cdot exp(a_1) \cdot \ldots \cdot exp(a_p) \cdot exp(b_1) \cdot \ldots \cdot exp(b_q)
\end{align*}
Shortly $G=K \cdot A \cdot B$, with $A=exp(\m_1)\cdot \ldots \cdot exp(\m_q)$ and $B=exp(\h_1)\cdot \ldots \cdot exp(\h_p)$.
\end{thm}
\begin{proof} See \cite[Main lemma, page 293]{Koh65} and \cite[Ch.IV, thm.3.2]{Loos69}. 
\end{proof}
\begin{cor}\label{cor:Koh} 
Let $G$ be a connected Lie group and $\Sigma$ an involutive automorphism of $G$. Then the fixed-point set $G^{\Sigma}$ has finitely
many connected components and $G^{\Sigma}/G^{\Sigma}_0$ is isomorphic to $(\Z_2)^k$ for some $k$. If $G$ is simply connected, 
then $G^{\Sigma}$ is connected. 
\end{cor}
\begin{proof} See \cite[Cor., page 293]{Koh65} and \cite[Ch.IV, thm.3.4]{Loos69}.
\end{proof}
 
According to proposition \ref{prop:2Cov} there are two extremes of symmetric coverings with simply connected total space:
\begin{enumerate}
\item[1)] $C=C_{\G}$, i.e. $C_{\H}$ is an isometry. These coverings are exactly those whose associated covering of the transvection groups
is non-trivial and the corresponding isotropy groups are connected.
\item[2)] $C=C_{\H}$, i.e. $C_{\G}$ is an isometry. In this case the associated covering of the transvection groups is an isomorphism 
\end{enumerate}
In particular we are in the second case if $Z(\G^1)=\{e\}$. The prototype is the following 
\begin{exa}
Consider the sphere $S^2$ with the standard metric. Its representation as symmetric space by its transvection group is $SO(3)/SO(2)$.
The center of $SO(3)$ is trivial but $Z_{Iso(S^2)}\G(S^2)=Z_{O(3)}SO(3)=\{\pm id\}$. The factor space is the real projective plane
$S^2/\{\pm id\}=\R P^2$.
\end{exa}

The fact that $C_{\H}$ consists of a sequence of double coverings can be motivated geometrically by the following observation:
\begin{re}
If the symmetric spaces in question are Riemannian, the points in a fiber of the covering $C_{\H}$ are lying pairwise on a closed geodesic.

Indeed, let be $x_1,x_2 \in {C_{\H}}^{-1}(y)$. Then $s_{x_1} \circ s_{x_2}$ lies in the kernel of the covering homomorphism of 
transvection groups corresponding to the symmetric covering $C_{\H}$. Due to \ref{prop:2Cov} the first covering is an isomorphism, hence
$s_{x_1} \circ s_{x_2}=id_{\G^2/\H^2_0}$.
Now let $\gamma$ be a geodesic through $x_1,x_2 $ (which exists since the manifold was assumed to be Riemannian). 
Since $s_{x_1} \circ s_{x_2}$ is on the one hand a transvection along $\gamma$ on the other hand 
the identity $\gamma$ has to be a closed geodesic. \hfill \qedsymbol
\end{re}

The aim of proposition \ref{prop:2Cov} was to show that symmetric coverings are  determined essentially 
(i.e. modulo the finite covering $C_{\H}$) by discrete subgroups of the center $Z(\G^1)$. 
Moreover the center of a connected Lie group lies in the image of the exponential map and can be 
described completely in terms of its Lie algebra (cf.  \cite[III.7.11]{HilgertN91}).




\subsection{Application to Solvable Symmetric Spaces}

\begin{prop}\label{prop:SolvSym}
Let $(M,g)$ be a simply connected solvable symmetric space. 
\begin{enumerate}
\item[(i)] $\G(M,g)$ is diffeomorphic to $\R^{\dim\G(M,g)}$.
\item[(ii)] The compact subgroup $K\subset \G$ as in theorem \ref{thm:Koh} is trivial. Thus $\H=B$ and $A$ acts
simply transitively on $M$. 
\end{enumerate} 
\end{prop}

\begin{proof}
(i) Lets denote by $\wt{\G}$ the simply connected covering of $\G$ and by $\wt{\H}$ the connected subgroup generated by $\h$ (identifying
$LA(\wt{\G})$ with $LA(\G)$). Then $\wt{\G}/\wt{\H} \simeq \G/\H$ and $\G \simeq \wt{\G}/N$, where $N$ is the largest normal subgroup of $\wt{\G}$ 
contained in $\wt{\H}$. $N$ being discrete yields $N \subset Z(\wt{\G})$ and thus $Ad(n)=id_{\g}, \; \forall n \in N$.
On the other hand $\wt{\H} \subset \wt{\G}$ is a nilpotent subgroup since $\h =[\m,\m] \subset [\g,\g]$ is a nilpotent subalgebra of $\g$.
Thus $n=exp(h_n)$ for some $h_n \in \h$
since the exponential map of a connected nilpotent Lie group is surjective (see \cite[III.3.24]{HilgertN91}). 
Hence $id_{\g}=Ad(n)=Ad(exp(h_n))=exp(\ad(h_n))$. Then $\ad(h_n)$ being a nilpotent endomorphism 
and $\ad(\cdot): \h \rightarrow \gl(\m)$ being faithful implies $h_n=0$.
This shows $N=\{e\}$ or equivalently $\wt{\G}=\G$. But the simply connected
covering of a solvable Lie group is diffeomorphic to $\R^{\dim\G(M,g)}$  (see
\cite[III.3.30]{HilgertN91}).

(ii)  Let $\G=K \cdot A \cdot B$ be the decomposition of the transvection group guaranteed by \ref{thm:Koh}. 
Then $LA(K)$ is reductive and solvable, hence
abelian. Thus $K$ is diffeomorphic to a torus (see
\cite[III.3.25]{HilgertN91}). On the other hand $\G$ is simply connected due
to (i), hence $K=\{e\}$.
Then the general identity
\begin{equation}\label{eq:H=KB}
G^{\Sigma}=K^{\Sigma} \cdot B
\end{equation}
leads in this case to $\H=\G^{\Sigma}=B$, hence $\G/\H=(A\cdot B)/B \simeq A$.
\end{proof}

\begin{cor}
Every simply connected solvable symmetric space is diffeomorphic to
$\R^{\dim(M)}$ and  admits a global pseudo-orthonormal frame-field, i.e. the
manifold is parallelizable. \end{cor}

Now we proceed to coverings of such spaces. The following theorem is known:
\begin{thm}[{\protect \cite[thm. 1 \& 2]{Koh65}}] 
Let $G$ be a Lie group with an involutive automorphism $\Sigma$ and a Lie subgroup $H$ such that
$G^{\Sigma}_0 \subset H \subset G^{\Sigma}$.
\begin{enumerate}
\item[(i)] With the notation from \ref{thm:Koh}, $G/\H$ is diffeomorphic to
the fiber-bundle $K \times_{K \cap H} A$ associated to the principle-bundle 
$K/(K \cap H)$, where $(K \cap H)\subset K^{\Sigma}$ acts on $A$ by conjugation.
\item[(ii)] If $G$ is
nilpotent, then $G/H$ is diffeomorphic to $T^k \times \R^l$ for some integers
$k,l \in \N_0$, where $T^k = (S^1)^k$ denotes the $k$-dimensional torus. 
\item[(iii)]  If $G$ is solvable, then $G/H$ is finite times covered by a space diffeomorphic to $T^k
\times \R^l$ for some integers $k,l \in \N_0$.  If $H$ is connected this
covering is trivial. 
\end{enumerate} 
\end{thm} %

The various theorems of S.Koh does not make any use of the representation of an (affine) symmetric space by its transvection group. We will investigate
the solvable case under this point of view:

Let $M^2=\G^2/\H^2$ be a solvable symmetric space and $M^1=\G^1/\H^1$ its simply connected covering. 
Assuming $\H^2$ to be connected 
the covering homomorphism $C:\G^1 \rightarrow \G^2$
restricts to an isomorphism $C|_{\H^1}:\H^1 \rightarrow \H^2$ (see proof of \ref{prop:2Cov}(ii)).
  
On the other hand \ref{thm:Koh} yields decompositions $\G^i=K^i \cdot A^i \cdot B^i$ for which $\H^i=(K^i)^{\Sigma^{(i)}}_0 \cdot B^i$ due to (\ref{eq:H=KB})
and the connectedness of $\H^i$. We have seen in \ref{prop:SolvSym}(ii) that
$\H^1$ is diffeomorphic to $B \simeq \R^{dim\H^1}$.  Since $\H^1,\H^2$ are
isomorphic, $\H^2$ is diffeomorphic to $\R^{dim\H^2}$, hence
$(K^i)^{\Sigma^{(i)}}_0=\{e\}$. Thus  
\[ M^2 \quad \simeq \quad {}^{\ds
\G^2}/_{\ds \H^2} \quad \simeq \quad {}^{\ds K^2 \cdot A^2 \cdot B^2}/_{\ds
B^2} \quad \simeq \quad  K^2 \cdot A^2 \] 
The space on the right is
diffeomorphic to $T^{\dim K^2} \times \R^{\dim A^2}$ as pointed out above.
In particular $ K^2 \cdot A^2$ acts again simply transitively on $M^2$. Finally $M^2$ being compact implies $A=\{e\}$. Thus $K$ being $\Sigma$-stable
and acting transitively on $M^2$ shows $\G^2 \subset K^2$ (compare \ref{prop:CharTrans}(ii)). 
Since $K^2$ is abelian as a compact subgroup of a solvable Lie 
group the corresponding symmetric space is flat. Thus we proved:

\begin{thm} For every solvable symmetric space $(M,g)$ their exists a finite symmetric covering $(\wt{M},\wt{g}) \rightarrow (M,g)$, where
$(\wt{M},\wt{g})$ is parallelizable and diffeomorphic to $T^k \times \R^l$ for
some integers $l,k \in \N_0$. If $M$ is compact (i.e. $l=0$) $(M,g)$ is flat.
\end{thm}

In \cite[\S4]{CahenW70} it has been asserted that there exist Lorentzian symmetric tori which are non flat. Later on in \cite{CahenK78} this
has been disproved while determining all possible symmetric coverings of solvable Lorentzian symmetric spaces by explicite calculation
of  their isometry  groups. 
In the subsequent section we will apply the theory developed above
to reproduce these results (up to a 2-fold covering).


\subsection{The realization of Lorentzian Symmetric Spaces}

Let $\g=\h \oplus \m $ be a solvable symmetric triple with its (first) adapted decomposition $\m= \z \oplus W_0 \oplus \z^*$
(\ref{lem:1.Decomp}).
Then $[\g,\g]=\h \oplus \z \oplus W_0 \subset \g$ is a nilpotent ideal and 
\begin{equation}\label{eq:Semi-directSequence}
\g \simeq  \Big( \ldots \big( [\g,\g] \rtimes \R Z^*_1 \big) \rtimes \ldots \Big) \rtimes \R Z^*_k
\end{equation}
for a basis $\{Z^*_k\}$ of $\z^*$ (the derivations of the semi-direct sums are induced by the corresponding adjoint representations). 
The simply connected Lie group $\G=LG(\g)$ of a semi-direct sum is obtained by an semi-direct product (see e.g.\cite[III.3.16]{HilgertN91})

\begin{prop} Let $\alpha :\mathfrak{b} \rightarrow \mathfrak{aut}(\mathfrak{a})$ be a Lie algebra homomorphism. Then
\[
LG(\mathfrak{a} \rtimes_{\alpha} \mathfrak{b}) \simeq LG(\mathfrak{a}) \rtimes_{\gamma} LG(\mathfrak{b})
\]
where $\gamma: LG(\mathfrak{b}) \rightarrow Aut \big(LG(\mathfrak{a}) \big)$ is the unique Lie group homomorphism with
$d\big(\gamma(b)\big)_e = \beta(b)$, $\forall b \in LG(\mathfrak{b})$ for the unique Lie group homomorphism
$\beta: LG(\mathfrak{b}) \rightarrow Aut(\mathfrak{a})$ with $d\beta_e=\alpha$.
\end{prop}

Since $[\g,\g]$ is nilpotent its simply connected Lie group
$\G^1=LG\big([\g,\g]\big)$ is isomorphic to $\big( [\g,\g],* \big)$ via $exp=id_{[\g,\g]}$
and the Campbell-Hausdorff product (cf. \cite[III.3.24]{HilgertN91}):
\[
X * Y =X + Y + \frac{1}{2}[X,Y]+ \ldots   
\]
Now set $\mathfrak{a}=[\g,\g]$ and $\mathfrak{b}=\R \cdot Z^*_1$. Identifying 
$ \mathfrak{b}\simeq LG(\mathfrak{b})$ the Lie algebra homomorphism 
$\alpha(\cdot)=\ad(\cdot)|_{[\g,\g]}:\mathfrak{b}\rightarrow \mathfrak{aut}(\mathfrak{a})$
lifts to the Lie group automorphism $\beta(\cdot)=exp(\ad(\cdot)|_{[\g,\g]})$ and further
$\beta=\gamma$. 
Hence $LG \big([\g,\g] \rtimes \R Z^*_1 \big)$ is explicitely known. Unless the signature of
the triple is $(1,n-1)$ this procedure has to be continued. But the further steps are more
subtle, since the exponential map of 
$LG \big( \big( \ldots ( [\g,\g] \rtimes \R Z^*_1) \rtimes \ldots \big) \rtimes \R Z^*_j \big)$
is in general not as easy to determine.  
The transvection group of simply connected nilpotent symmetric triples could be derived
directly from the Campbell-Hausdorff formula as done for the triple $\tau_{2,2}(\epsilon_Y)$
in \cite{CahenP70}. 
In the following we restrict ourself to the Lorentzian case.


\begin{prop} Let $\tau_n(f)$ be the normal form of a solvable symmetric Lorentz triple (cf. p. \pageref{triple:tau(f)}).
\begin{enumerate}
\item[(i)] The simply connected group $\G_f=LG(\g_f)$ is diffeomorphic to 
$\g_f=\h \oplus \m = X \oplus W \oplus \z \oplus \z^*$ with the product given by
\[
\big(x + w,z,z^* \big)\cdot \big( \ol{x} + \ol{w},\ol{z},\ol{z^*} \big)= 
\Big( e^{-\ad \ol{z^*}}(x + w) + \ol{x} + \ol{w},   z + \ol{z} + 
\frac{1}{2} \big[e^{-\ad \ol{z^*}}(x + w),\ol{x} + \ol{w} \big] , z^* + \ol{z^*}\Big)
\]
\item[(ii)] The simply connected symmetric space $\big( \G_f / \H , g_f \big)$ associated to $\tau_n(f)$
is isometric to $W \oplus \z \oplus \z^*$ endowed with the metric
\begin{equation*}
(g_f)_{(w,z,z^*)}=2 \,dz^*\left( dz - \sum_{\alpha=1}^{n-2} f^{\alpha}(w_{\alpha})^2 dz^* \right) + 
\sum_{\alpha=1}^{n-2} (dw_{\alpha})^2
\end{equation*}
\item[(iii)] The center of the transvection group is 
\begin{equation*}
Z(\G_f)=
\begin{cases}
\z & \text{if} \; \exists \, \alpha \text{ s.t. } f^{\alpha} > 0 \; \text{ or if } 
\; \exists \, \alpha, \beta \text{ s.t } \frac{f^{\alpha}}{f^{\beta}}\not\in \mathbb{Q}^2  \\
\z \times \{\lambda \Z \cdot Z^*\} & \text{else,} \quad  
\text{ with } \lambda = \min \{l > 0 \;|\; l \sqrt{f^{\alpha}} \in 2 \pi i \Z, \; \forall \alpha\}
\end{cases}
\end{equation*}
\end{enumerate}
\end{prop}

\begin{proof} (i): The first derivative of $\g_f$ is spanned by 
$\{X_{\alpha},W_{\beta},Z\}_{1 \le \alpha,\beta \le n-2}$. The only non-zero
commutators are $[X_{\alpha},W_{\alpha}]=-f^{\alpha}\cdot Z$. This Lie algebra is known as
$(2n-3)$-dimensional Heisenberg algebra with nil-index $2$. Thus the Campbell-Hausdorff formula
reduces to 
\begin{equation*}
\big( x,w,z \big) * \big( \ol{x},\ol{w},\ol{z} \big) = 
\big( x + \ol{x}, w + \ol{w}, z + \ol{z} + \frac{1}{2}[x+w,\ol{x}+\ol{w}]\big)
\end{equation*} 
Then as explained above $\G_f=LG(\g_f) \simeq [\g_f,\g_f] \oplus \R \cdot Z^* = \g_f$ with the
usual multiplication of a semi-direct product:
\begin{align*}
\big(x + w,z,z^* \big) \circ \big( \ol{x} + \ol{w},\ol{z},\ol{z^*} \big)
&=\Big( (x+w,z)*\gamma(z^*)(\ol{x} + \ol{w},\ol{z}),z^* + \ol{z^*} \Big) \\
&=\Big( e^{\ad z^*}(\ol{x} + \ol{w}) + x + w,   z + \ol{z} + 
\frac{1}{2} \big[x + w,e^{\ad z^*}(\ol{x} + \ol{w}) \big] , z^* + \ol{z^*}\Big)
\end{align*}
Instead of the multiplication ``$\circ$'' we will consider that one given in the proposition
since the two groups $(\g_f,\circ)$ and $(\g_f,\cdot)$ are isomorphic under the map
$\phi(g^{-1})= -g$ (inversion with respect to ``$\circ$'').
This is exactly the group cited in \cite[thm. 5b]{CahenW70}.

(ii): The connected subgroup $\H=\langle exp(\h) \rangle \subset \g_f$ is $\h$. It acts from
the right on $\G_f$ by
\[
\big(x + w,z,z^* \big)\cdot \big( \ol{x} ,0,0 \big)= 
\Big( x + \ol{x} + w ,   z  + 
\frac{1}{2} \big[w,\ol{x}\big] , z^* \Big)
\]
In order to keep this action simple the isomorphism $\phi$ has been introduced. Thus 
\begin{align*}
s: \quad \G_f/\H  \quad &\rightarrow \quad \G_f \\
(x+w,z,z^*) \,\H \quad &\mapsto \quad (w,z-\frac{1}{2} [w,x],z^* )
\end{align*}  
is a global section of the submersion $\pi:\G_f \rightarrow \G_f/\H$.
Now we extend the scalar product $B$ on $\g_f \simeq T_0 \G_f$ to a bi-invariant metric on $\G_f$.
Then  $\pi$ being a pseudo-Riemannian submersion yields the metric $g_f$ on $\G_f/\H$ 
and $\big( \G_f/\H , g_f \big)$ is isometric to $\big(W  \oplus \z \oplus \z^* ,(g_f) \big)$ with
$(g_f)_{s  \pi (k)} := (s \circ \pi \circ l_k)_* B|_{\m \times \m}$.

The identification $\R^{n} \simeq \m$ by means of the basis 
$\beta_0=\{W_{\alpha},Z,Z^*\}$ yields a global chart of $\m$. The corresponding
basis in $T_p\m$ $(p \in \m)$ induced by the coordinate vector fields will be denoted by $\beta_p$.
Set $k=(x+w,z,z^*) \in \G_f$ and $p=(s \circ \pi) (k)$. Then  
\[
Mat_{\beta_p}\big((g_f)_p \big)= Mat_{\beta_0,\beta_p} \big( d(s \circ \pi \circ l_k)|_{\m} \big)^{-t}\cdot 
Mat_{\beta_0}(B|_{\m \times \m}) \cdot Mat_{\beta_0,\beta_p} \big( d(s \circ \pi \circ l_k)|_{\m} \big)^{-1}
\]
Concretely
\begin{align*}
(dl_k)_0 (\wh{x}+\wh{w},\wh{z},\wh{z^*}) 
& = \frac{d}{dt}\Big( e^{-\ad t\wh{z^*}}(x + w) + t\wh{x} + t\wh{w}
 z + t\wh{z} + \frac{1}{2} \big[e^{-\ad t\wh{z^*}}(x + w), t\wh{x} + t\wh{w} \big] , z^* + t\wh{z^*}\Big)_{t=0}\\
& =\Big( [x+w,z^*] + \wh{x} + \wh{w}, \wh{z} + \frac{1}{2} [x + w, \wh{x} + \wh{w}],\wh{z^*} \Big)\\
d(s \circ \pi)_k (\wh{x}+\wh{w},\wh{z},\wh{z^*}) 
& = \frac{d}{dt}\Big( s\circ \pi (x+t\wh{x}+w+t\wh{w},z+t\wh{z},z^*+t\wh{z^*}) \Big)_{t=0}\\
& = \frac{d}{dt}\Big( w+t\wh{w},z+t\wh{z}+\frac{1}{2}[w+t\wh{w},x+t\wh{x}],z^*+tz^* \Big)_{t=0}\\
& = \Big(  \wh{w}, \wh{z} - \frac{1}{2}\big([x,\wh{w}]-[w,\wh{x}]\big),\wh{z^*}\Big)
\end{align*}
In order to determine the metric $g$ we may restrict ourself to $k \in \m$ since this
subset of $\G_f$ acts already transitively on $\G_f/\H$:
\begin{align*}
&d(s \circ \pi \circ l_k)|_{\m} (\wh{w},\wh{z},\wh{z^*}) = d(s \circ \pi)_k \Big( [w,z^*] + \wh{w}, \wh{z}, \wh{z^*} \Big)
 = \Big( \wh{w}, \wh{z} + \frac{1}{2}\big[w,[w,z^*] \big],\wh{z^*}\Big)\\
\Rightarrow \quad &Mat_{\beta_0,\beta_p} \big( d(s \circ \pi \circ l_k)|_{\m} \big) 
 = \left(\begin{matrix} I_{n-2} &0&0\\0&1&\frac{1}{2}\sum f^{\alpha}(w_{\alpha})^2\\0&0&1 \end{matrix}\right)
\end{align*}
Since the group isomorphism $\phi$ satisfies $d\phi_0 = id_{\g_f}$ we get further
\[
Mat_{\beta_0}(B|_{\m \times \m})= \left(\begin{matrix} I_{n-2} &0&0\\0&0&1\\0&1&0 \end{matrix}\right)
\]
Now the assertion follows immediately.

(iii): For a connected Lie group its center is given by the kernel of the adjoint representation.
\begin{multline*}
Ad(x+w,z,z^*)(\wh{x}+\wh{w},\wh{z},\wh{z^*}) =\\
\shoveleft{ = \frac{d}{dt} \Big( (x+w,z,z^*)\cdot (t\wh{x}+t\wh{w},t\wh{z},t\wh{z^*}) \cdot 
\underbrace{(-e^{\ad z^*}(x+w),-z,-z^*)}_{\ds =(x+w,z,z^*)^{-1}} \Big)_{t=0} }\\
\shoveleft{ = \frac{d}{dt} \Big( e^{\ad z^*} \big( e^{-\ad t\wh{z^*}}(x+w)+t\wh{x}+t\wh{w} - x -w \big),
 t\wh{z} +  \frac{1}{2}\big[ e^{-\ad t\wh{z^*}}(x+w),t\wh{x}+t\wh{w}\big] - \ldots }\\
\shoveright{ \ldots -\frac{1}{2} e^{\ad z^*} \big[ e^{-\ad t\wh{z^*}}(x+w) + t\wh{x}+t\wh{w},x+w\big],t\wh{z^*}\Big)_{t=0}}\\
\shoveleft{ = \Big( e^{\ad z^*}  \big( -\ad \wh{z^*}(x+w)+\wh{x}+\wh{w}\big),\wh{z} +  \big[ x+w,\wh{x}+\wh{w}\big] 
 + \frac{1}{2}\big[ \ad \wh{z^*}(x+w),x+w\big],\wh{z^*}\Big) }\\
\end{multline*}
Thus $(x,w,z,z^*) \in \ker(Ad)$ if and only if the following conditions hold $\forall (\wh{x}+\wh{w},\wh{z},\wh{z^*})$:
\begin{enumerate}
\item[(i)] $e^{\ad z^*}  \big( -\ad \wh{z^*}(x+w)+\wh{x}+\wh{w}\big)=\wh{x}+\wh{w}$
\item[(ii)] $\big[ x+w,\wh{x}+\wh{w}\big] + \frac{1}{2}\big[ \ad  \wh{z^*}(x+w),x+w\big]=0$ 
\end{enumerate} 
The first condition is equivalent to 
\[
\ad \wh{z^*}(x+w)= \big( id - e^{-\ad z^*} \big)(\wh{x}+\wh{w}),\quad \forall (\wh{x}+\wh{w},\wh{z},\wh{z^*})
\]
Since the left side does not depend on $\wh{x},\wh{w}$ the whole equation has to be identically zero hence:\
\[
(x+w,z,z^*) \in \ker(Ad) \quad \Leftrightarrow \quad x=w=0 \text{ and } e^{\ad(z^*)}=id_{\g_f}
\]
It is easy to see that $\ad(s Z^*)$ is diagonalizable and $spec\big( \ad(s Z^*) \big)=\{0,s\sqrt{f^{\alpha}} \}$.
Thus $e^{\ad(z^*)}=id_{\g_f}$ if all eigenvalues of $\ad(z^*)$ are imaginary and have a common integer multiple. 
\end{proof} 
An element $(0,z,z^*) \in Z(\G_f)$ acts on $\m \simeq \G_f/\H$ simply by $(0,z,z^*) \cdot (\ol{w},\ol{z},\ol{z^*})=(\ol{w},z+\ol{z},z^*+\ol{z^*})$.
Thus the factor space with respect to a discrete subgroup of $Z(\G_f)$ can be represented by a fundamental domain as in the flat case.

The full centralizer $Z_{Iso_f}\G_f$ has been determined for the Lorentzian case in \cite[prop. 5]{CahenK78} (see also \cite{Baum00}).
If $Z(\G_f)=\z$ the centralizer coincides with the center. In the other case $Z_{Iso_f}\G_f$ is generated by $\z$ and the isometry 
\[
\phi \big( w_1 +\ldots + w_{n-2} ,z,z^* \big) =  \big( (-1)^{m_1}w_1 + \ldots + (-1)^{m_{n-2}}w_{n-2},z,z^*+ \frac{\lambda}{2}Z^* \big) 
\]
where $w_{\alpha} \in \R \cdot W_{\alpha}$, $\lambda = \min \{l > 0 \;|\; l \sqrt{f^{\alpha}} \in 2 \pi i \Z, \; \forall \alpha\}$ 
and $m_{\alpha}=\frac{ \lambda\sqrt{f^{\alpha}} }{2 \pi i} \in \Z$.
Note that on the factor-space $\m/\lambda \Z \cdot Z^*$ the isometry induced by $\phi$ is of order two as could be expected from \ref{prop:2Cov}(iii).

\bibliographystyle{amsalpha} 
\bibliography{literature}

\end{document}